\pdfoutput=1      

\documentclass[journal]{new-aiaa}
\usepackage[utf8]{inputenc}
\usepackage{textcomp}

\usepackage{amsmath}
\usepackage[version=4]{mhchem}
\usepackage{siunitx}
\usepackage{longtable,tabularx}
\usepackage{caption}
\usepackage{subcaption}

\DeclareMathOperator*{\argmax}{argmax}
\usepackage{multirow}
\usepackage{makecell}
\usepackage{graphicx}

\title{Safe and Optimal N-Spacecraft Swarm Reconfiguration in Non-Keplerian Cislunar Orbits}

\author{
Yuji Takubo \footnote{Ph.D. Candidate, Department of Aeronautics and Astronautics, 496 Lomita Mall, AIAA Student Member, Corresponding Author (ytakubo@stanford.edu)}, 
Walter Manuel \footnote{Ph.D. Candidate, Department of Aeronautics and Astronautics, 496 Lomita Mall, AIAA Student Member}, 
Ethan Foss \footnote{Ph.D. Student, Department of Aeronautics and Astronautics, 496 Lomita Mall, AIAA Student Member}, and 
Simone D'Amico \footnote{Associate Professor, Department of Aeronautics and Astronautics, 496 Lomita Mall, AIAA Associate Fellow}}
\affil{Stanford University, Stanford, CA, 94305}

\begin{document}

\maketitle

\footnotetext{The Preliminary version of this work has been presented at the 2025 AAS/AIAA Space Flight Mechanics Meeting as manuscript 25-331 \cite{takubo2025passively}.}

\begin{abstract}
This paper presents a novel fuel-optimal guidance and control methodology for spacecraft swarm reconfiguration in Restricted Multi-Body Problems (RMBPs) with a guarantee of passive safety, maintaining miss distance even under abrupt loss of control authority. A new set of constraints exploits a quasi-periodic structure of RMBPs to guarantee passive safety. Particularly, the condition for passive safety is expressed as simple geometric constraints by solving optimal control in Local Toroidal Coordinates, which is based on a local eigenspace of a quasi-periodic motion around the corresponding periodic orbit. The proposed formulation enables a significant simplification of problem structure, which is applicable to large-scale swarm reconfiguration in cislunar orbits. The method is demonstrated in the Circular Restricted Three-Body Problem, the Elliptic Restricted Three-Body Problem, and the Bi-Circular Restricted Four-Body Problem. Furthermore, the optimized control profiles are validated in the full-ephemeris dynamics model. By extending and generalizing well-known concepts of relative orbital elements within the restricted two-body problem to the three- and four-body problems, this paper lays the foundation for practical control schemes of relative motion in cislunar space.
\end{abstract}


\section*{Nomenclature}

{\renewcommand\arraystretch{1.0}
\noindent\begin{longtable*}{@{}l @{\quad=\quad} l@{}}
$m$ & mass \\
$M$ & monodromy matrix \\
$N$ & number of time steps \\
$N_{sc}$ & number of deputy spacecraft \\
$\boldsymbol{r}_{ab}$ & relative position vector from $a$ to $b$ \\
$\boldsymbol{w}$ & eigenvector (oscillatory mode) \\
$\boldsymbol{\zeta}$ & local toroidal coordiantes \\
$\boldsymbol{\theta}$ & toroidal variable  \\
$\boldsymbol{\omega}_{A/B}$ & angular velocity vector of A with respect to B \\
$\mu$ & non-dimensional mass ratio \\
$\Phi$ & state transition matrix \\
$\boldsymbol{\rho}$ & relative position vector \\
\multicolumn{2}{@{}l}{Subscripts}\\
$b_1$ / $\mathcal{B}$ & Earth-Moon barycenter synodic frame \\
$c$ & chief spacecraft \\
$d$ & deputy spacecraft \\
$e$ & Earth \\
$i / \mathcal{I}$ & inertial frame \\
$j$ & deputy's index \\
$l / \mathcal{L}$ & Local-Horizontal-Local-Vertical (LVLH) frame\\
$m / \mathcal{M}$ & moon / Moon synodic frame \\
$v / \mathcal{V}$ & Velocity/Normal/Binormal (VNB) frame \\
$s$ & Sun \\
\end{longtable*}}

\section{Introduction}
\lettrine{D}{istributed} Space Systems (DSS) are rapidly emerging as a transformative paradigm in modern spaceflight, enabling new mission capabilities through spacecraft constellations, swarms, formations, and Rendezvous, Proximity Operations, and Docking (RPOD).
To date, almost all RPOD missions have been executed in orbital regimes of the Perturbed Restricted Two-Body Problem (PR2BP), where Keplerian motion dominates and perturbations are several orders of magnitude weaker (e.g., Earth orbits, low lunar orbits). 
Recently, there has been growing interest in deploying DSS along periodic orbits governed by the Restricted Multi-Body Problem (RMBP).
One such motivation comes from the projected deployment of NASA's Gateway in a Near-Rectilinear Halo Orbit (NRHO) of the Earth-Moon $L_2$ libration point \cite{lee2019gateway} and associated RPOD operations. 
Such a DSS requires Guidance, Navigation, and Control (GNC) algorithms that enable bounded relative motion between multiple spacecraft in the encountered multi-body dynamics. 

Traditionally, dynamics models and state representations for relative motion between spacecraft rely on a fundamental assumption of orbiting around a single primary (e.g., Earth) based on the PR2BP. 
These dynamics are often linearized around the chief orbit within Local-Vertical Local-Horizontal (LVLH) or Radial-Tangential-Normal (RTN) frames \cite{clohessy_terminal_1960, yamanaka_new_2002}.
In the last few decades, dynamics models in Relative Orbital Elements (ROE) \cite{damico2010autonomous, schaub2003analytical}, which are based on the difference in the Keplerian orbital elements, have also been developed as integration constants and fundamental matrix solutions of the aforementioned dynamics models.
The ROE-based GNC algorithms have been demonstrated in a number of flight missions, such as GRACE \cite{montenbruck2006GRACE}, PRISMA \cite{damico2012prisma}, TanDEM-X \cite{ardaens_control_2008}, AVANTI \cite{gaias2016design}, Starling \cite{kruger2023starling}, SWARM-EX \cite{lowe2024concept}, VISORS \cite{Guffanti2023Visors}, and ADRAS-J \cite{Atarashi2024ADRAS-J}. 
Previous studies also include various perturbation sources such as geopotential spherical harmonics, solar radiation pressure, atmospheric drag, and third-body perturbations \cite{sullivan2017comprehensive, koenig2017new}. 
Early research on relative motion between spacecraft on a periodic orbit in the RMBPs begins with the linearization of the nonlinear dynamics in the rotating frame \cite{scheeres2003stabilizing, segerman2002preliminary, ardaens_control_2008}.
More recently, relative dynamics of the Circular Restricted Three-Body Problem (CR3BP) and the Elliptic Restricted Three-Body Problem (ER3BP) have been derived in the LVLH frame \cite{franzini_relative_2019, khoury_relative_2022}. 
This is further extended to the Bi-Circular Restricted Four-Body Problem (BCR4BP) to express a higher-fidelity dynamics model in the Sun-Earth-Moon system \cite{Romero2023cislunar}. 
Furthermore, the relative dynamics of the CR3BP resolved in a co-moving frame based on spacecraft velocity (``NTW frame" in Ref. \citenum{vallado2001fundamentals_5th}, ``velocity-reference frame" in Refs. \cite{albert2023relative, vela2025modal}, ``Velocity/Normal/Conormal (VNC) frame" in Ref. \cite{plice2019helioswarm}, and ``Velocity/Normal/Binormal (VNB) frame" in NASA General Mission Analysis Tool (GMAT) \cite{hughes2016general}) are derived in Refs. \cite{takubo2025passively, vela2025modal} contemporaneously.

Stable solutions in the CR3BP provide crucial insights into the design of bounded relative motion in PRMBPs.
Besides five libration points and periodic orbit families surrounding them, equilibria exist in the form of Quasi-Periodic Orbits (QPOs) that naturally trace the surface of a Quasi-Periodic Invariant Torus (QPIT) and envelop the nearby frequency-matching periodic orbit without an intersection \cite{barden1998formation, olikara2016computation, baresi2017spacecraft}.
Thus, a spacecraft on a QPO performs not only a quasi-periodic bounded motion to an attractor (e.g., Moon) in the rotating frame but also a quasi-periodic bounded relative motion with respect to a spacecraft in its co-moving orbital frame.
For this reason, applications of QPOs to long-term formation-flight and proximity operations in the RMBPs have been thoroughly investigated \cite{barden1998formation, simanjuntak2012design, henry2024phd}. 
The characteristics of relative orbital motion are generalized by Lyapunov-Floquet theory \cite{Yakubovich1975linear} with analysis of the eigensystem of the dynamics. 
In the PR2BP, six unit eigenvalues of the monodromy matrix reflect the six first integrals (Integration Constants (ICs) \cite{guffanti2018integration, manuel2022optimal}) that lead to Relative Orbital Elements (ROEs) \cite{damico2010autonomous, guffanti2022phd, delurgio2024closed}. 
Meanwhile, in the PMBPs, the monodromy matrix and Floquet matrix are used to analyze the stability of the (quasi-)periodic orbits to analyze the fundamental structure of the quasi-periodic solutions \cite{olikara2016computation, baresi2017spacecraft, down2023relative}.
Originally devised for different dynamical regimes, the two approaches converge within the methods of fundamental modal solutions that represent relative motion between spacecraft as a superposition of modes whose coefficients are prescribed by the ICs \cite{burnett2022spacecraft}.
The ICs are derived for the relative motion in CR3BP, and the impulsive control scheme is solved within this state space \cite{vela2025modal}. 
Although this method provides a mathematically convenient structure for discussing the natural perturbation, a set of ICs (ROEs) lacks geometric insight into the spacecraft formation.
Recently, a new state representation, Local Toroidal Coordinates (LTC) \cite{elliott2022phd, elliott2022describing}, has been proposed to relate the dynamical motion to the geometry of the torus and offer convenient first-order conditions for bounded relative motion between spacecraft. 
Previous works adopted the LTC to perform a long-term relative station-keeping of multi-spacecraft swarming in the CR3BP, which is then applied to the full-ephemeris dynamics of the Sun-Earth system \cite{elliott2022describing}. 
However, multi-revolution formation reconfiguration in cislunar orbits becomes more complicated due to a larger dynamics modeling error (e.g., eccentricity of Earth-Moon orbit, Sun's perturbations, etc.), as well as the fundamental complexity of the constrained trajectory planning compared to the unconstrained formation-keeping. 

Rigorous safety criteria are one of the most critical considerations for bounded relative motion.
Particularly, passive safety, which guarantees sufficient miss distance between spacecraft, even for free-drift trajectories that emanate from a sudden loss of control authority, is critical for practical operations in a chaotic dynamical system. 
However, this has received little attention in the relative motion control in multi-body dynamics. 
Recent works discuss the safety of quasi-periodic relative orbits by computing the minimum separation between a QPO and the corresponding periodic orbit in the CR3BP using the LTC \cite{elliott2022phd} as well as the fundamental modal solution space \cite{vela2025modal}.
Particularly, Ref.~\cite{vela2025modal} performs its analysis in the velocity-referenced frame, noting that the local-horizontal axis of the LVLH frame departs from the velocity direction on highly elliptic orbits such as NRHOs. 
Aligning the frame with the velocity direction keeps the navigation-error ellipsoid naturally elongated along its axis, which simplifies the construction of geometric safety regions of Keep-Out Zones (KOZs).
Besides, formation-keeping control strategies that can target specific regions of the local eigenspace using both impulsive control \cite{elango2022local, henry2024phd} and continuous control \cite{segerman2002preliminary, scheeres2003stabilizing, marchand2005control, bando2015formation} are developed to realize a specific geometry with some separations. 
In contrast, some works bypass the natural dynamics structure and perform reconfigurations between safe QPROs with a Linear Quadratic Regulator (LQR) controller in the BCR4BP \cite{khoury_relative_2024} or Model Predictive Control (MPC) with Cartesian collision-avoidance constraints in the CR3BP \cite{sekine2024qpt, capannolo2021dynamics}. 
Nonetheless, the major limitation persists as a design that maintains passive safety throughout the transfer is more convoluted than safe orbit designs. 
Passive safety in the relative motion in the PR2BP has been mainly addressed by careful relative orbit design that ensures planar separation \cite{montenbruck2006GRACE, damico2010autonomous}, or forward/backward propagation of drift trajectories over a fixed duration \cite{breger2008safe, marsillach2020fail, guffanti2023passively, takubo2024art, elango2024fast}. 
In RMBPs, the former is difficult due to distinct dynamical structures. 
The latter approach has been applied to cislunar rendezvous \cite{innocenti_failure_2022}, although this is also limited by the need for a long safety horizon due to nonlinear sensitivity in PRMBPs \cite{innocenti_failure_2022} and by the quadratic growth of constraints with the number of spacecraft.
Consequently, there is a need for scalable safety-preserving constraints that extend beyond the CR3BP to mid- and high-fidelity relative dynamics models in RMBPs. 

With these key limitations in mind, this paper proposes a passively safe, fuel-optimal reconfiguration of a spacecraft swarm in RMBPs, bridging the research gap between well-developed studies of bounded relative motion and control in PR2BP and those in RMBPs. 
Five main contributions are presented in this paper.
First, the relative dynamics in the RMBPs and the corresponding torus structures in the co-moving frame are compared between the LVLH frame and the velocity-reference frame. 
The analysis presents that the QPO resolved in the co-moving frame (termed as a Quasi-Periodic Relative Orbit (QPRO)) yields a compact bounded motion, containing large oscillatory behavior only in the velocity direction. 
This provides an easier understanding of relative motion in the RMBPs compared to a description in the LVLH frame. 
Secondly, the geometric relationship between QPROs in the RMBPs and periodic relative orbits in the PR2BP is harmonized in this paper. 
It is shown that the local eigenspace associated with the periodic orbit in the Earth-Moon rotating frame defines the geometric configuration of the quasi-periodic relative motion. 
Moreover, insights into the minimum separation in relevant projections of the co-moving local frames are provided.
Thirdly, the LTC is constructed on a sidereal-resonant periodic orbit in the ER3BP and a synodic-resonant periodic orbit of the BCR4BP, extended from earlier works on relative motion in the CR3BP \cite{elliott2022phd, takubo2025passively} to a higher-fidelity dynamics model. 
Fourthly, an Optimal Control Problem (OCP) is solved in the LTC for the first time. 
Leveraging this state representation realizes a significantly simpler and more intuitive form of safety constraints. 
Furthermore, the proposed safety constraints are entirely decentralized for a given pair of relative orbits, which facilitates an extension to swarm control involving tens of spacecraft while maintaining simplicity in onboard implementation.
Finally, the optimized controls are deployed in the full-ephemeris model of the Sun-Earth-Moon system with an MPC scheme, providing a promising statistical result on the applicability of the proposed method in the full-ephemeris dynamics via Monte Carlo simulation. 
Overall, the proposed method enables the swarm reconfiguration in non-Keplerian cislunar orbits that is scalable to tens of spacecraft, which could open new possibilities for mission designs of cislunar RPOD or swarming.

The remainder of the paper is organized as follows. 
First, nonlinear and linear equations of relative motion for the RMBPs in the velocity-reference frame are derived. 
Secondly, the passive safety of QRPOs is examined in relation to the relative orbit design in the PR2BP, followed by the derivation of the LTC.
Third, using the invariance property of the LTC, an OCP is formulated as a nonconvex optimization problem and solved via Sequential Convex Programming (SCP) \cite{malyuta2022convex}. 
The control performances of the optimized trajectories are examined in the three dynamic models. 
Furthermore, the results of the OCPs are validated with full-ephemeris dynamics. 
Finally, conclusions and future research directions are provided.

\section{Dynamics Models}

This paper considers a spacecraft swarm comprising a (virtual) chief and $N_{sc}$ deputies, operating near the Moon and subject to gravitational forces of multiple attractors, including the Earth, Moon, and Sun.
In particular, RMBPs assume that the mass of each spacecraft is negligible compared to that of these celestial bodies. 

\subsection{Notations}
The following notations are used in this paper. 
The positions of the chief and deputy with respect to the Moon's center of mass are denoted as $\boldsymbol{r}$ and $\boldsymbol{r}_d$, respectively. 
The relative position of the deputy with respect to the chief is expressed as $\boldsymbol{\rho} = \boldsymbol{r}_d - \boldsymbol{r}$.
For a vector $\boldsymbol{\zeta}$, $\zeta = \| \boldsymbol{\zeta} \|$ denotes its $l_2$-norm. 
A vector with a $\hat{\boldsymbol{\zeta}}$ indicates a unit vector. 
$[\boldsymbol{\zeta}]^\times$ denotes a skew-symmetric matrix that provides a cross product ($\boldsymbol{\zeta} \times$).  
The first, second, and third-order time derivatives resolved in frame $\mathcal{A}$ are denoted as 
$[\dot{\boldsymbol{\zeta}}]_\mathcal{A} = [\partial \boldsymbol{\zeta}/\partial t]_\mathcal{A}$, 
$[\ddot{\boldsymbol{\zeta}}]_\mathcal{A} = [\partial ^2 \boldsymbol{\zeta}/\partial t^2]_\mathcal{A}$, and 
$[\dddot{\boldsymbol{\zeta}}]_\mathcal{A} = [\partial ^3 \boldsymbol{\zeta}/\partial t^3]_\mathcal{A}$, respectively.
A vector expressed in frame $\mathcal{A}$ is denoted as $[\boldsymbol{\zeta}]_\mathcal{A}$. 
$\text{diag}([\zeta_1, ..., \zeta_n])$ represents a diagonal matrix with diagonal entries corresponding to the matrix representation of a vector $[\zeta_1, ..., \zeta_n]$.
The operator $\because$ is a shorthand sign for ``because".
The $d$-dimensional real space is denoted as $\mathbb{R}^{d}$, whereas the $d$-dimensional toroidal space is denoted as $\mathbb{T}^{d}$.

\subsection{Coordinate Frames}

\begin{figure}[ht!]
     \centering
     \begin{subfigure}[b]{0.4\textwidth}
         \centering
         \includegraphics[width=\textwidth]{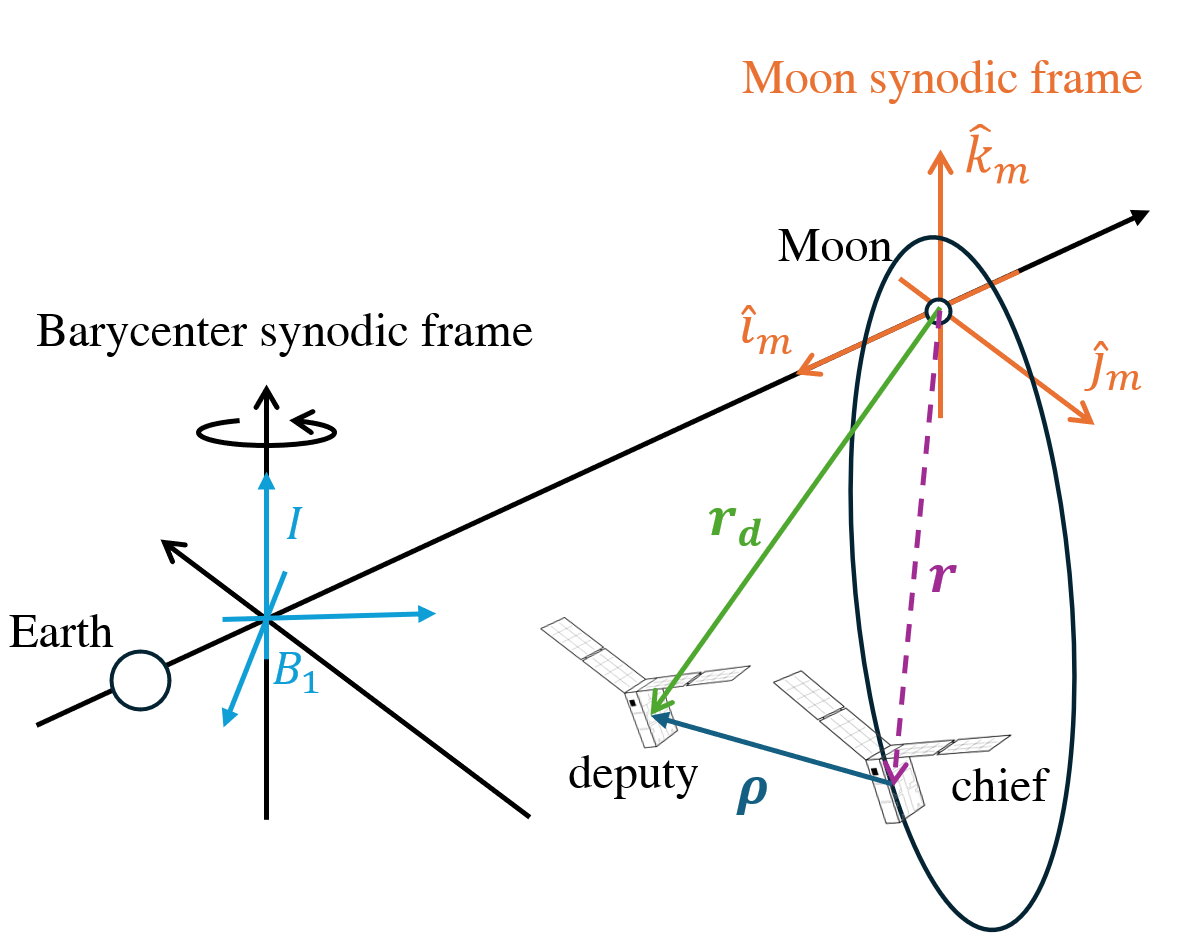}
         \caption{Synodic and inertial frames}
         \label{fig:coord_moon}
     \end{subfigure}
    \begin{subfigure}[b]{0.4\textwidth}
         \centering
         \includegraphics[width=\textwidth]{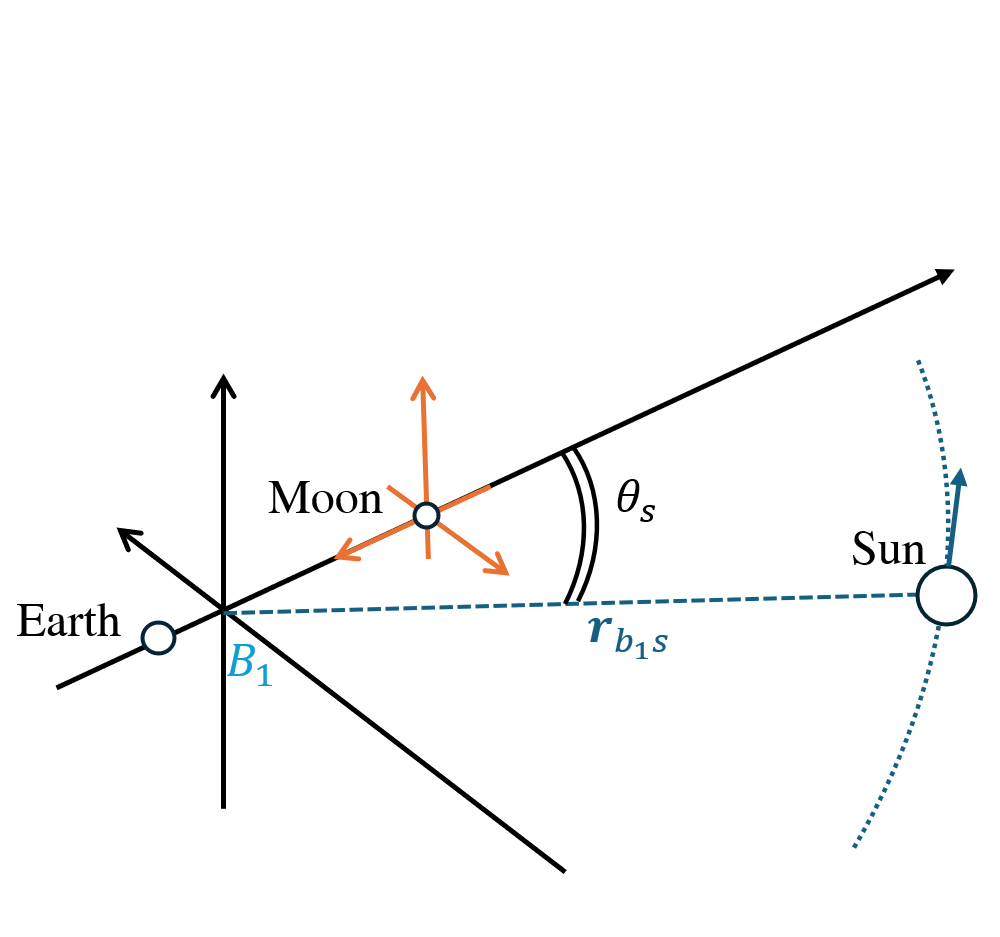}
         \caption{Sun's position}
         \label{fig:coord_EMS}
     \end{subfigure} 
     \\
     \begin{subfigure}[b]{0.4\textwidth}
         \centering
         \includegraphics[width=\textwidth]{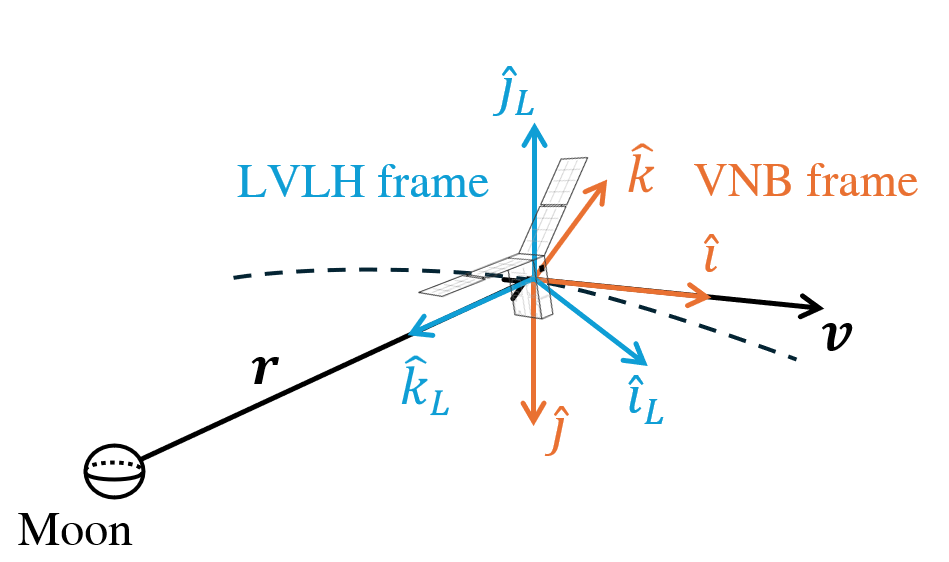}
         \caption{Local coordinate frames}
         \label{fig:coord_tnw_and_lvlh}
     \end{subfigure}  
    \caption{Coordinate frames and vectors used in this work. }
    \label{fig:coord_frame} 
\end{figure}

The coordinate frames used in this article is summarized in Fig. \ref{fig:coord_frame}.
First, Fig. \ref{fig:coord_moon} illustrates coordinate frames within the Earth-Moon system: inertial frame $\mathcal{I}$, barycenter synodic frame $\mathcal{B}$, and Moon synodic frame $\mathcal{M}$. 
The barycenter synodic frame is a rotating frame about the Earth-Moon barycenter $B_1$. 
The Moon synodic frame is attached to the barycenter synodic frame, while the origin of the frame is anchored to the Moon's center of mass.

To describe the relative position and velocity between spacecraft near a (quasi-)periodic orbit (called reference orbit) in the synodic frame, two co-moving frames are employed: LVLH frame $\mathcal L=\{\hat{\imath}_L,\hat{\jmath}_L,\hat{k}_L\}$ and VNB frame $\mathcal{V}:\{\hat{\boldsymbol{\imath}}, \hat{\boldsymbol{\jmath}}, \hat{\boldsymbol{k}}\}$, as shown in Fig. \ref{fig:coord_tnw_and_lvlh}. 
The LVLH frame is defined based on the position of the spacecraft as 
\begin{align}
    \hat{\boldsymbol{\imath}}_L = \hat{\boldsymbol{\jmath}}_L\times\hat{\boldsymbol{k}}_L, \quad 
\hat{\boldsymbol{\jmath}}_L = -(\boldsymbol{r}\times \boldsymbol{v})/ \|\boldsymbol{r}\times \boldsymbol{v}\|, \quad 
\hat{\boldsymbol{k}}_L = -{\boldsymbol{r}}/{r},
\end{align}
whereas the VNB frame is defined based on the velocity of the spacecraft as 
\begin{equation} \label{eq:tnw_basis}
    \boldsymbol{\hat{\imath}} = \boldsymbol{v} / {v}, \quad 
    \boldsymbol{\hat{\jmath}} = (\boldsymbol{r}\times \boldsymbol{v})/ \|\boldsymbol{r}\times \boldsymbol{v}\|, \quad 
    \boldsymbol{\hat{k}} = \boldsymbol{\hat{\imath}} \times \boldsymbol{\hat{\jmath}}.
\end{equation}
Here, $\boldsymbol{v}$ is the velocity vector resolved in the synodic frame.
The LVLH axes are labeled as X, Y, and Z axes, forming the XY, YZ, and ZX planes.
Similarly, the VNB frame uses the V, N, and B axes, which define the BV, VN, and NB planes.

\subsection{Dynamics of RMBPs}

Three dynamical models of the RMBPs are discussed in this paper: the Circular Restricted Three-Body Problem (CR3BP), Elliptic Restricted Three-Body Problem (ER3BP), and Bicircular Restricted Four-Body Problem (BCR4BP). 
The CR3BP assumes that the two primary bodies (i.e., Earth and Moon) execute circular orbits about their barycenter. In contrast, the ER3BP relaxes this assumption by permitting elliptical motions for the primaries, which provides a more realistic representation of their motion. 
Furthermore, the BCR4BP extends the CR3BP formulation by incorporating an additional gravitational source. 
In this paper, the tertiary body (i.e., Sun) is assumed to orbit the barycenter of the primary system circularly, as shown in Fig. \ref{fig:coord_EMS}.

By convention, the equations of motion in the Earth-Moon system are non-dimensionalized by the (a) Earth-Moon distance (length), (b) the sum of the primary mass (mass), and (c) angular rate of primaries around the Earth-Moon barycenter (time), so that the gravitational constant of the Earth-Moon system becomes unity.  
For all these dynamical models, the equations of motion resolved in the barycenter synodic frame are written as follows: 
\begin{align}
    [\ddot{\boldsymbol{r}}_{b_1c}]_{\mathcal{B}} = 
    - 2\boldsymbol{\omega}_{b/i} \times [\dot{\boldsymbol{r}}_{b_1c}]_{\mathcal{B}} 
    - [\dot{\boldsymbol{\omega}}_{b/i}]_{\mathcal{B}} \times \boldsymbol{r}_{b_1c} 
    + \nabla \Upsilon,
\end{align}
where $\Upsilon$ is a pseudo-potential function; 
$\boldsymbol{\omega}_{b/i}$ is the angular velocity of the barycenter synodic frame with respect to the inertial frame, whose components in the barycenter synodic frame are given as  $[\omega^x_{b/i}, \omega^y_{b/i}, \omega^z_{b/i}]$; 
and $\boldsymbol{r}_{b_1c}$ is the position of the chief with respect to $B_1$, whose components in the barycenter synodic frame are given as $[x,y,z]$.
For example, the pseudo-potential of the CR3BP is given by  
\begin{align}
    \Upsilon_{\text{CR3BP}} = \frac{1}{2}\omega_{b/i}^{z^2}(x^2 + y^2) + \frac{\mu}{r_{ec}} + \frac{1-\mu}{r_{mc}},
\end{align}
where $\mu = \dfrac{m_m}{m_e+m_m}$ is the mass ratio of the two primaries.
The CR3BP assumption leads to $[\omega^x_{b/i}, \omega^y_{b/i}, \omega^z_{b/i}] = [0,0,1]$ (i.e., $[\dot{\boldsymbol{\omega}}_{b/i}]_\mathcal{B}=\boldsymbol{0}$).
The pseudo-potenitals of ER3BP and BCR4BP are derived in a similar manner and are shown in Appendix A.
For the ER3BP, the angular velocity and the Earth-Moon distance are the time-varying parameters. 
The BCR4BP parameterizes the Sun's position in the barycenter synodic frame $\boldsymbol{r}_{b_1s}$ as a function of $\theta_s$, which is the angle between the Earth-Moon line and the Sun-$B_1$ line, as shown in Fig. \ref{fig:coord_EMS}.
Within the barycenter synodic frame, the components of $\boldsymbol{r}_{b_1s}$ are described as $[x_s, y_x, z_s] = [r_{b_1s} \cos\theta_s, r_{b_1s} \sin\theta_s, 0]$, which rotates around the Earth-Moon barycenter at a constant angular velocity $\boldsymbol{\omega}_s$ with its non-dimensionalized mass $\mu_s = \dfrac{m_s}{m_e+m_m}$.

\subsection{Relative Dynamics in the VNB} frame

Following the derivations of the relative dynamics in the CR3BP, ER3BP, and BCR4BP within the LVLH frame \cite{franzini_relative_2019, Romero2023cislunar}, this subsection develops nonlinear and linear relative dynamics in R3BP in the VNB frame. 
While this derivation assumes that the spacecraft swarm is anchored to the Moon, the extension to orbits around other anchoring points (e.g., libration points) is readily available \cite{khoury_relative_2022}. 

A key difference in the derivation in the VNB frame to that in the LVLH frame \cite{franzini_relative_2019} is the analytical expression of the angular velocity of the VNB frame with respect to the Moon synodic frame $\boldsymbol{\omega}_{v/m}$ and angular acceleration $[\dot{\boldsymbol{\omega}}_{v/m}]_{\mathcal{V}}$, which are derived as  
\begin{subequations}
\begin{align}
    \boldsymbol{\omega}_{v/m} 
    & = {\omega}^{x}_{v/m} \hat{\boldsymbol{\imath}} + {\omega}^{y}_{v/m} \hat{\boldsymbol{\jmath}} + {\omega}^{z}_{v/m} \hat{\boldsymbol{k}} \\
    & = \dfrac{1}{h} ( ( \boldsymbol{r} \times [\ddot{\boldsymbol{r}}]_{\mathcal{M}}) \cdot {\hat{\boldsymbol{k}}})  \hat{\boldsymbol{\imath}} 
    - \dfrac{1}{v} ([\ddot{\boldsymbol{r}}]_\mathcal{M} \cdot \hat{\boldsymbol{k}})  \hat{\boldsymbol{\jmath}} 
    + \dfrac{1}{hv} ( [\ddot{\boldsymbol{r}}]_\mathcal{M} \cdot \boldsymbol{h} ) \hat{\boldsymbol{k}} 
\end{align}
\end{subequations}
and  
\begin{subequations}
\begin{align}
    [\dot{\boldsymbol{\omega}}_{v/m}]_{\mathcal{V}} & = \dot{\omega}^{x}_{v/m} \hat{\boldsymbol{\imath}} + \dot{\omega}^{y}_{v/m} \hat{\boldsymbol{\jmath}} + \dot{\omega}^{z}_{v/m} \hat{\boldsymbol{k}}, \\
    \dot{\omega}^{x}_{v/m} & = 
    - \dfrac{\dot{h}}{h^2} ((\boldsymbol{r} \times [\ddot{\boldsymbol{r}}]_{\mathcal{M}}) \cdot \hat{\boldsymbol{k}} ) 
    + \dfrac{1}{h} \left(  
    ([\dot{\boldsymbol{r}}]_{\mathcal{M}} \times [\ddot{\boldsymbol{r}}]_{\mathcal{M}} ) \cdot \hat{\boldsymbol{k}} 
    + (\boldsymbol{r} \times [\dddot{\boldsymbol{r}}]_{\mathcal{M}} ) \cdot \hat{\boldsymbol{k}} 
    + (\boldsymbol{r} \times [\ddot{\boldsymbol{r}}]_{\mathcal{M}}) \cdot [\dot{\hat{\boldsymbol{k}}}]_{\mathcal{M}}  \right) \\
    \dot{\omega}^{y}_{v/m} & = \dfrac{\dot{v}}{v^2} ([\ddot{\boldsymbol{r}}]_{\mathcal{M}} \cdot \hat{\boldsymbol{k}}) 
    - 
    \frac{1}{v} ([\dddot{\boldsymbol{r}}]_{\mathcal{M}} \cdot \hat{\boldsymbol{k}} +[\ddot{\boldsymbol{r}}]_{\mathcal{M}} \cdot [\dot{\hat{\boldsymbol{k}}}]_{\mathcal{M}}  )
    \\ 
    \dot{\omega}^{z}_{v/m} & = \dfrac{-\dot{h}v - h\dot{v}}{h^2 v^2} ([\ddot{\boldsymbol{r}}]_{\mathcal{M}} \cdot \boldsymbol{h}) 
    +
    \dfrac{1}{hv} ( [\ddot{\boldsymbol{r}}]_{\mathcal{M}} \cdot \boldsymbol{h}).
\end{align}
\end{subequations}
The derivation is summarized in Appendix B. 
Furthermore, the acceleration and the jerk of each dynamics model resolved in the Moon synodic frame are shown in Appendix C.
The generalized nonlinear equation of relative motion within the Sun-Earth-Moon system is computed as follows \cite{franzini_relative_2019, Romero2023cislunar}: 
\small
\begin{align}
\begin{split}
    [\ddot{\boldsymbol{\rho}}]_{\mathcal{V}}
    = & -2 [\boldsymbol{\omega}_{v/i}]_{\mathcal{V}}^{\times} [\dot{\boldsymbol{\rho}}]_{\mathcal{V}}
    - \left( 
        [\dot{\boldsymbol{\omega}}_{v/i}]_{\mathcal{V}}^{\times} 
        + \left( [{\boldsymbol{\omega}}_{v/i}]_{\mathcal{V}}^{\times} \right)^2 
    \right) \boldsymbol{\boldsymbol{\rho}} 
    \\ & 
    + \mu \left( \frac{\boldsymbol{r}}{r^3} - \frac{\boldsymbol{r} + \boldsymbol{\rho}}{\| \boldsymbol{r} + \boldsymbol{\rho} \|^3} \right) 
    + (1-\mu) \left( \frac{\boldsymbol{r}_{em} + \boldsymbol{r} }{\| \boldsymbol{r}_{em} + \boldsymbol{r} \|^3} - \frac{\boldsymbol{r}_{em} + \boldsymbol{r} + \boldsymbol{\rho}}{\| \boldsymbol{r}_{em} + \boldsymbol{r} + \boldsymbol{\rho} \|^3} \right) 
    + \mu_s \left( \frac{\boldsymbol{r}_{sm} + \boldsymbol{r} }{\| \boldsymbol{r}_{sm} + \boldsymbol{r} \|^3} - \frac{\boldsymbol{r}_{sm} + \boldsymbol{r} + \boldsymbol{\rho}}{\| \boldsymbol{r}_{sm} + \boldsymbol{r} + \boldsymbol{\rho} \|^3} \right).
    \label{eq:CNERM}
\end{split}
\end{align}
\normalsize
Linearization of Eq. \eqref{eq:CNERM} at the chief's position yields the follwoing state-space form of linearized dynamics: 
\begin{subequations} \label{eq:plant_mat_TNW}    
\begin{align} 
    & \dot{\boldsymbol{x}} = \boldsymbol{A}(\boldsymbol{r}(t)) \boldsymbol{x} 
    = \begin{bmatrix}
        \boldsymbol{0}_3 & \boldsymbol{I}_3 \\
        \boldsymbol{A}_{\dot{\boldsymbol{\rho}}\boldsymbol{\rho}} & -2 [\boldsymbol{\omega}_{v/i}]^{\times} 
    \end{bmatrix} \boldsymbol{x}, 
    \quad \boldsymbol{x} = \begin{bmatrix}
        \boldsymbol{\rho} \\ [\dot{\boldsymbol{\rho}}]_{\mathcal{V}}
    \end{bmatrix},
    \\
    & 
\boldsymbol{A}_{\dot{\boldsymbol{\rho}}\boldsymbol{\rho}} 
    = - [\dot{\boldsymbol{\omega}}_{v/i}]_{\mathcal{V}}^{\times} 
    - \left( [\boldsymbol{\omega}_{v/i}]_{\mathcal{V}}^{\times} \right)^2 
    - \mu \frac{\partial}{\partial \boldsymbol{r}}\left[\frac{\boldsymbol{r}}{r^3}\right] 
    - (1-\mu) \frac{\partial}{\partial \boldsymbol{r}}\left[\frac{\boldsymbol{r}+\boldsymbol{r}_{em}}{\left\|\boldsymbol{r} + \boldsymbol{r}_{em}\right\|^3}\right]
    - \mu_s \frac{\partial}{\partial \boldsymbol{r}}\left[\frac{\boldsymbol{r}+\boldsymbol{r}_{sm}}{\left\|\boldsymbol{r} + \boldsymbol{r}_{sm}\right\|^3}\right].
\end{align}
\end{subequations}
\normalsize
Note that ${\boldsymbol{\omega}}_{v/i} = {\boldsymbol{\omega}}_{v/m} + {\boldsymbol{\omega}}_{m/i}$ 
and
$[\dot{\boldsymbol{\omega}}_{v/i}]_{\mathcal{V}} = [\dot{\boldsymbol{\omega}}_{v/m}]_{\mathcal{V}} + [\dot{\boldsymbol{\omega}}_{m/i}]_{\mathcal{M}} - \boldsymbol{\omega}_{v/m} \times \boldsymbol{\omega}_{m/i}$ \cite{franzini_relative_2019}.

\section{Bounded Relative Motion in Multi-body Problems} \label{sec:bounded_relative_motion}

In this section, a thorough comparison is presented between natural bounded relative motions in the RP2BP and the RMBPs. 
Due to nonlinear and chaotic nature of dynamics, the computation of bounded relative motions in the RMBPs has heavily relied on numerical continuation methods \cite{olikara2016computation, baresi2017spacecraft}.
However, by analyzing the linearized dynamical system and the fundamental solution matrix \cite{guffanti2022phd} offers a geometric insight into the relative motion in RMBPs.
Available dynamical modes in two dynamical systems are first presented, followed by the introduction of the LTC \cite{elliott2022describing} that leverages the dynamical mode that excites the bounded relative motion in the RMBPs. 
To analyze the dynamics, this section focuses on the derivation in the VNB frame as a specific choice of the local frame. However, this can be readily replaced with the LVLH frame or other local frames. 

\subsection{Eigensystem in the Vicinity of Periodic Orbits}

Consider a periodic solution in a dynamical system $\boldsymbol{f}$ with a period of $T$. 
Here, $T$ is the time it takes for the state of the system to return exactly to its initial condition under the chosen dynamics model, potentially over multiple revolutions of orbits. 
The first-order dynamics in the vicinity of a periodic orbit can be characterized by its monodromy matrix
\begin{align} \label{eq:monodromy}
    M(t_0) := \Phi(t_0, t_0+T) = 
    \Phi(t_0, t_0) + \int_0^{T} A(\tau) \Phi(t_0, t_0+\tau) d\tau, \quad \Phi(t_0, t_0) = \boldsymbol{I}_6,
\end{align}
where a plant matrix $A(\tau) = \dfrac{\partial f}{\partial \tau}$ is the linearized dynamics of a system evaluated at the reference orbit. 
The eigensystem of the monodromy matrix defines the linear stability of the perturbed states near the reference orbit.

Also, consider a linear mapping from the Cartesian state perturbation $\boldsymbol{x} = [\boldsymbol{\rho}^\top, [\dot{\boldsymbol{\rho}}]_{\mathcal{V}}^\top]^\top$ from the reference orbit to a set of new state representation $\boldsymbol{\kappa}(t)\in \mathbb{R}^6$ as follows: 
\begin{align} \label{eq:dx_to_integ_const}
    \boldsymbol{\kappa}(t) = \Psi^{-1}(t) \boldsymbol{x}. 
\end{align}
Then, the spectrum of the monodromy matrix is retained in the $\boldsymbol{\kappa}(t)$-space due to a similarity transform %
\begin{align} \label{eq:monod_transform}
M_{\boldsymbol{\kappa}}(t_0) = \Psi^{-1}(t_0+T)  M(t_0)  \Psi(t_0) = \Psi^{-1}(t_0) M(t_0) \Psi(t_0). 
\end{align}
Note that $\Psi(t_0+T) = \Psi(t_0) \forall t_0$ because of the periodicity of the reference orbit. 

\subsubsection{Eigensystem in the R2BP} \label{sec:eigen_2BP}

\begin{figure}[ht!]
     \centering
     \begin{subfigure}[b]{0.49\textwidth}
         \centering
         \includegraphics[width=\textwidth]{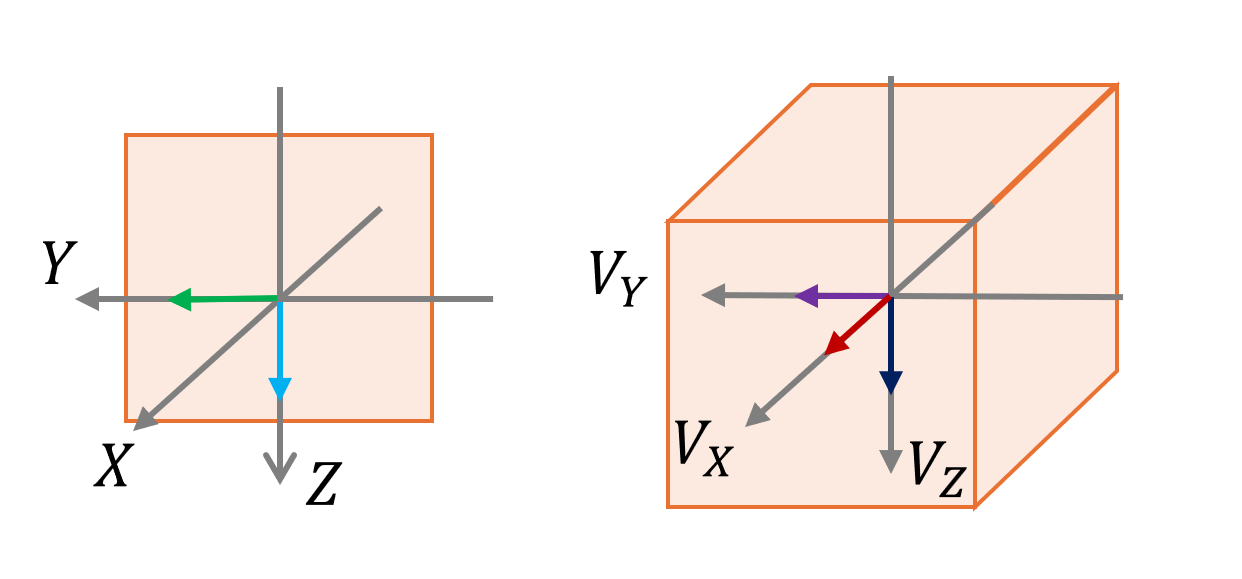}
         \caption{ Central mode in the relative motion of the R2BP. }
         \label{fig:eigvec_periodic}
     \end{subfigure}
     \begin{subfigure}[b]{0.49\textwidth}
         \centering
         \includegraphics[width=\textwidth]{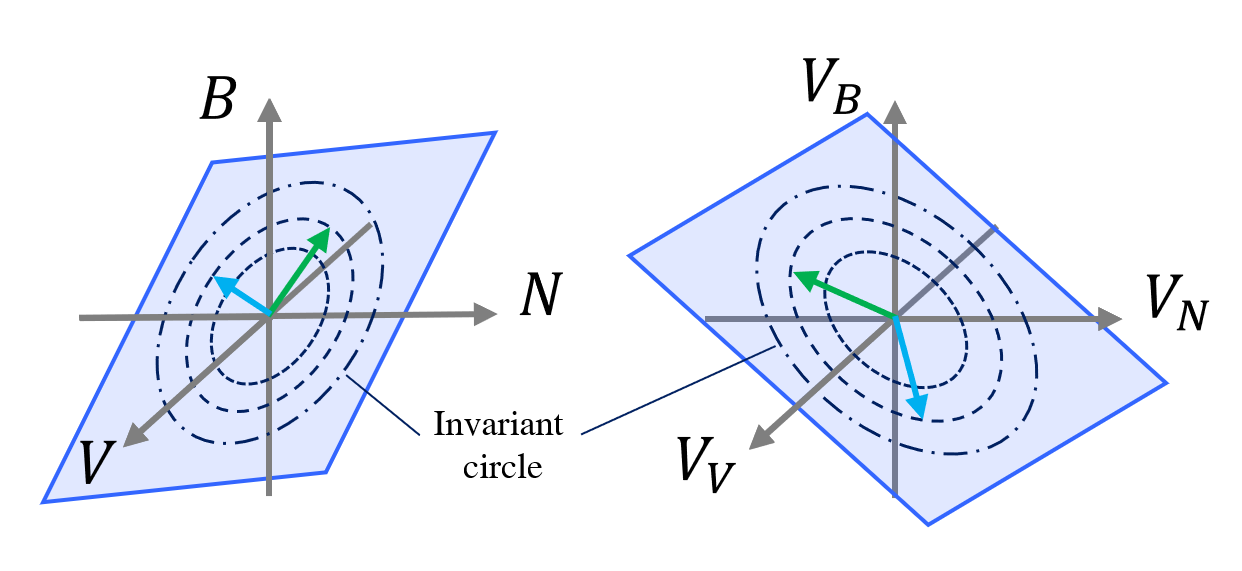}
         \caption{ Oscillatory mode in the relative motion of the CR3BP. }
         \label{fig:eigvec_quasiperiodic}
     \end{subfigure}  
    \caption{ Local eigenspace of the monodromy matrix: central mode in the R2BP and the oscillatory mode in the CR3BP. }
    \label{fig:eigvecs} 
\end{figure}

Periodic orbits in the Keplerian dynamics are integrable and possess six first integrals, which implies that all six eigenvalues of the monodromy matrix are one, with five linearly independent eigenvectors
\cite{vallado2001fundamentals_5th, schaub2003analytical}. 
The dynamical mode associated with the unity eigenvalue is called a central mode.
The five eigenvectors are unit vectors of all six directions except for a mode in the along-track direction, as shown in Fig. \ref{fig:eigvec_periodic}. 
This eigensystem induces periodic relative motion from any initial condition, except for the along-track direction, where a difference in specific mechanical energy induces drift \cite{damico2010autonomous}. 
Therefore, based on  Eq.~\eqref{eq:monod_transform}, $M_{\boldsymbol{\kappa}}$ retains six unit eigenvalues (i.e., periodicity of $\boldsymbol{\kappa}(t)$) for an arbitrary linear invertible map $\Psi(t)$. 
Thus, a set of invariant integration constants $\boldsymbol{\kappa}(t) = \boldsymbol{\kappa}$ can be defined with $\Psi(t)$ being a fundamental matrix solution to the variational equation, which leads to the well-known ROE theory \cite{damico2010autonomous} based on two periodic orbits, as shown in Fig. \ref{fig:rel_motion_2bp}.
Note that this eigensystem provides five degrees of freedom in the design space of bounded relative orbits. 

When additional perturbations are considered (i.e., PR2BP), some eigenvalues of the monodromy matrix perturb from unity. 
Due to the additional conditions required for periodicity, the solution set for periodic relative motion is reduced. 
Previous works address bounded relative orbit designs that cancel a secular drift caused by the $J_2$-perturbation \cite{damico2010autonomous, koenig2018robust}. 

\subsubsection{Eigensystem in the RMBP}

The local linearization of the dynamics of RMBPs contains time-varying components \cite{franzini_relative_2019}, preventing a closed-form expression of the monodromy matrix. 
Instead, the monodromy matrix can be numerically integrated (cf. Eq. \eqref{eq:monodromy}) concurrently with the periodic solution (cf. Eqs. \eqref{eq:rddot_M_cr3bp}, \eqref{eq:rddot_M_er3bp}, and \eqref{eq:rddot_M_bcr4bp}) for the time interval $[t_0, t_0+T]$, using the plant matrix $\boldsymbol{A}(t)$.
This system becomes a 42-dimensional ordinary differential equation. 

Generally, not all eigenvalues of the monodromy matrix computed along a periodic solution of RMBP are unity.
Therefore, the entire solution set of bounded relative motion around a periodic orbit within the RMBPs is severely constrained.
One trivial solution is a set of relative motions between spacecraft on the same periodic orbit but with phase offsets.
Besides, a periodic orbit in the RMBPs may admit up to two oscillatory modes in the corresponding monodromy matrix, characterized by a conjugate pair of complex unimodular eigenvalues \cite{scheeres2003stabilizing, elliott2022describing}. 
When the oscillatory mode is excited, an initial relative state within the corresponding eigenspace evolves into a quasi-periodic motion; after one period $T$, the relative position and velocity do not repeat but remain confined within a rotating ellipse, known as an invariant circle. 
This structure is illustrated in Fig. ~\ref{fig:eigvec_quasiperiodic} as dotted ellipses on the eigenspace.
\begin{figure}[ht!]
     \centering
     \begin{subfigure}[b]{0.45\textwidth}
         \centering
         \includegraphics[width=\textwidth]{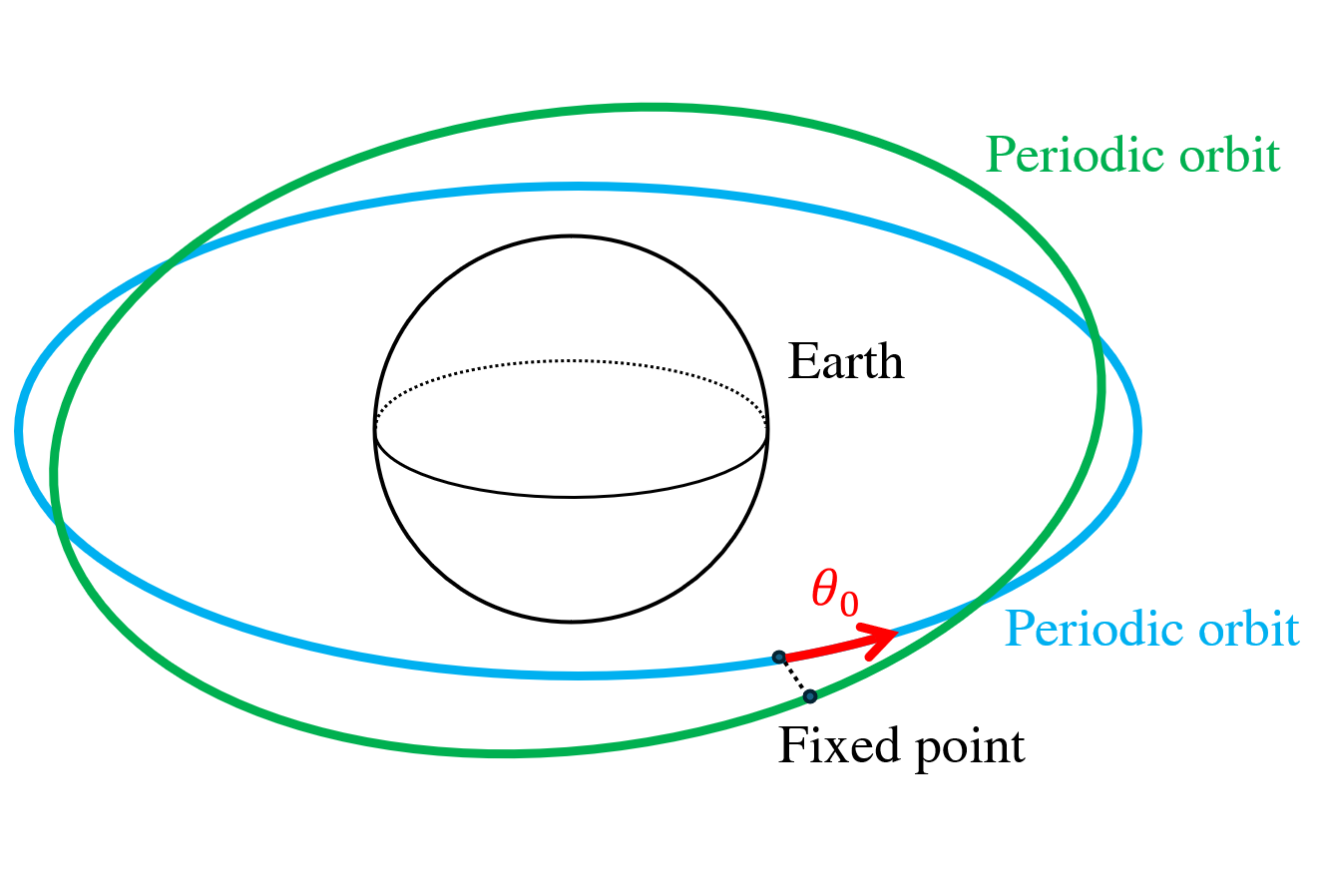}
         \caption{Periodic orbits around Earth characterized by one-dimensional torus $\mathbb{T}$}
         \label{fig:rel_motion_2bp}
     \end{subfigure}
     \begin{subfigure}[b]{0.45\textwidth}
         \centering
         \includegraphics[width=\textwidth]{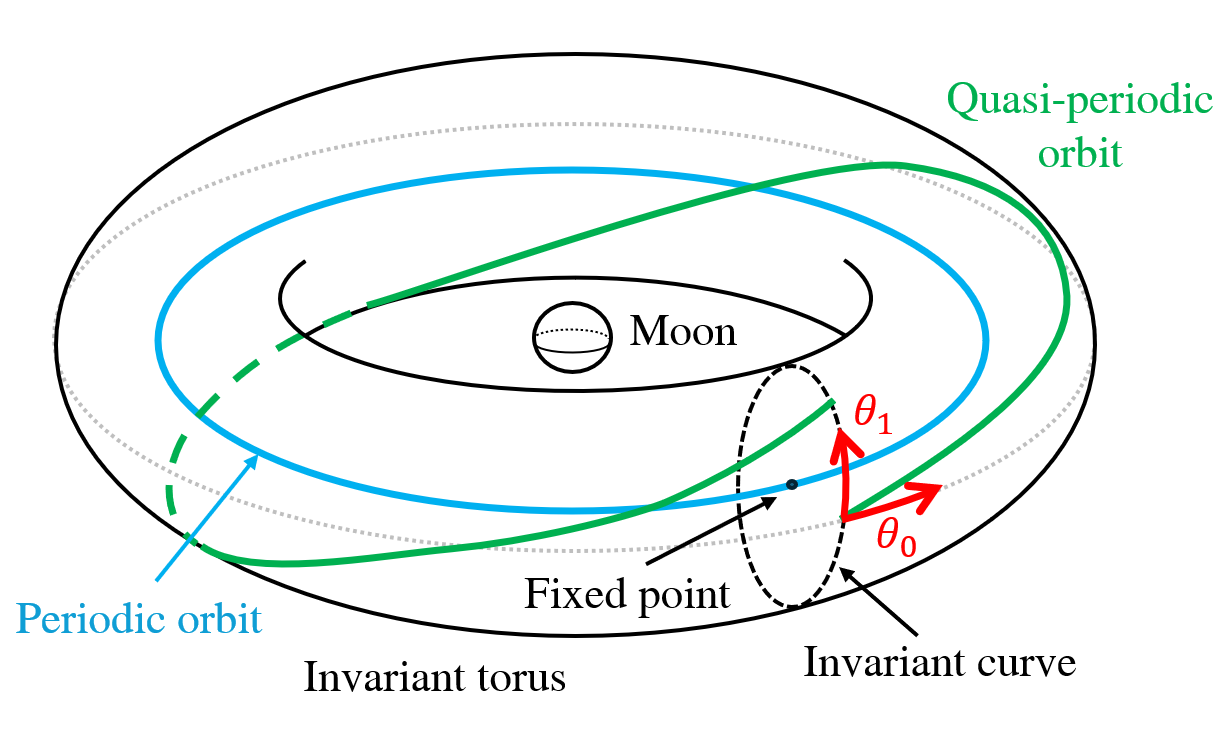}
         \caption{Two-dimensional ($\mathbb{T}^2$) quasi-periodic invariant torus and a quasi-periodic orbit}
         \label{fig:qpt}
     \end{subfigure}  
    \caption{Periodic and quasi-periodic motion in the R2BP and RMBP.}
    \label{fig:relative_2bp_and_3bp} 
\end{figure}

A QPRO constructed based on the oscillatory mode in the linearized dynamics is first-order equivalent to a QPO in the synodic frame that envelopes a reference orbit under the original nonlinear dynamics.
The QPO evolves along the surface of a two-dimensional Quasi-Periodic Invariant Torus (QPIT), which constitutes a nonlinearly stable geometric structure in the phase space of the RMBPs. 
This torus exhibits bounded perturbation responses in transverse directions and preserves its quasi-periodic motion under small disturbances. 
This enables solutions to evolve along it without altering their overall shape, as shown in Fig.~\ref{fig:qpt}.
This bounded motion is ensured because each point on the torus surface follows the vector field of the system's equations of motion. 
Any position and velocity ($\mathbb{R}^6$) on the surface of a QPIT can be characterized by a two-dimensional phase vector $\boldsymbol{\theta} \in \mathbb{T}^2$. 
Consequently, a motion initiated from a point on the torus can be viewed as an evolution of $\boldsymbol{\theta}$. 
The progression along a periodic orbit of the RMBPs is characterized by $\theta_0$, whereas $\theta_1$ represents the phase within an invariant curve, defined as a closed curve where quasi-periodic solutions are confined under the $T$-stroboscopic map.
In this paper, a QPIT resolved in the local frame around the chief is termed Quasi-Periodic Relative Invariant Torus (QPRIT). 
Compared to the QPIT around the attractor, the structure of the QPRIT around the chief spacecraft provides a fundamental geometry of a QPRO, enabling a simple characterization of the relative motion as shown in the following sections.
Thus, an invariant circle is equivalent to the first-order approximation of the invariant curve (cf. Fig. \ref{fig:qpt}) in a local frame \cite{olikara2016computation}, leading to two degrees of freedom in the QPRO design in each oscillatory mode. 

\subsection{Local Toroidal Coordinates} 

Local Toroidal Coordinates (LTC) \cite{elliott2022describing} are coordinate systems on a local frame along a reference orbit that project the dynamics onto the oscillatory eigenspace of the RMBPs. 
Instead of a full linear mapping to integration constants, a mode-isolating transformation is employed to emphasize the dynamical mode that provides a quasi-periodic bounded relative motion. 
In this paper, it is assumed that the chief is on a periodic orbit of RMBPs that contains at least one oscillatory mode.
 
\begin{figure}[ht!]
    \centering    \includegraphics[width=0.5\textwidth]{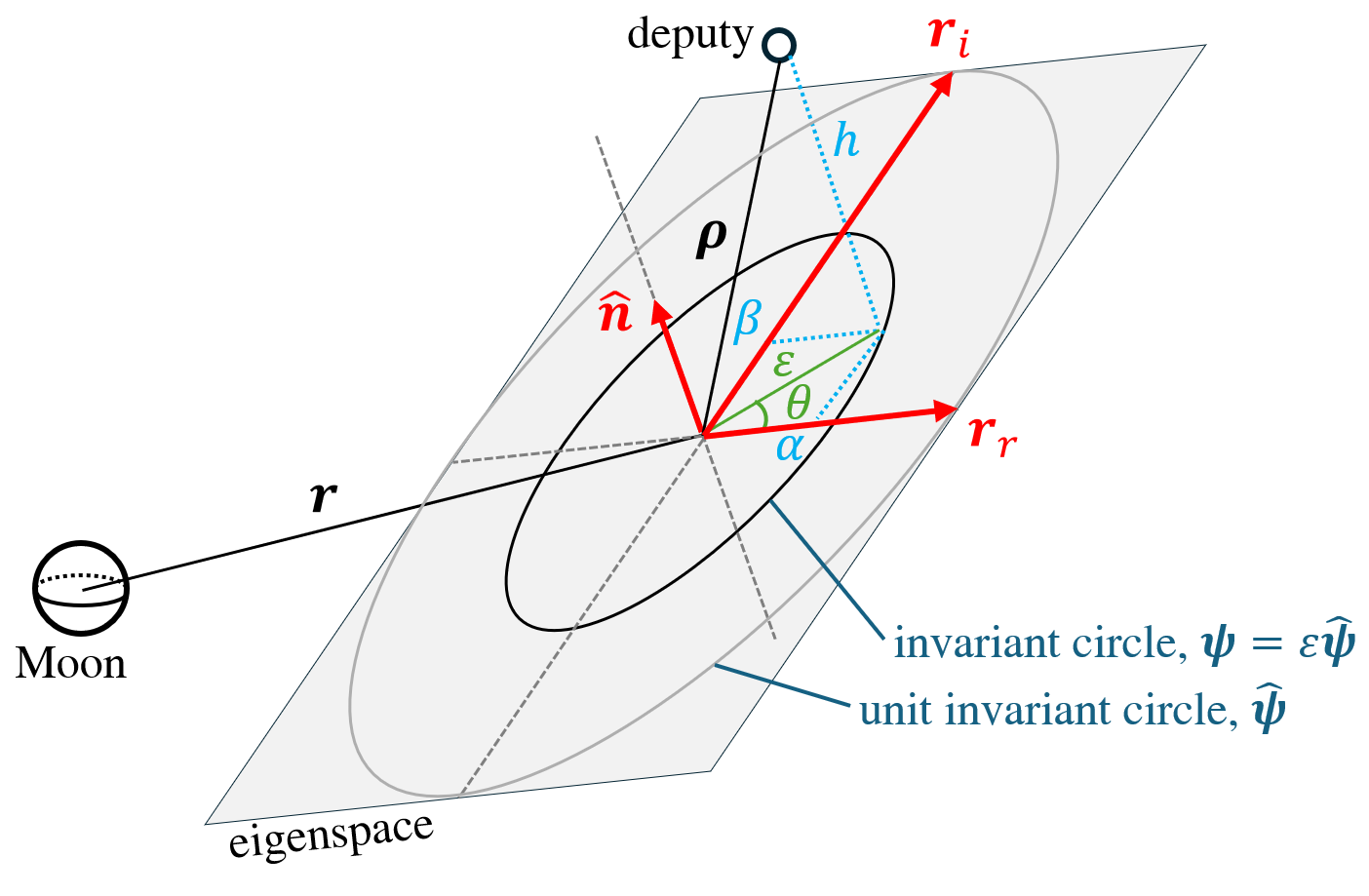}
    \caption{Nonsingular local toroidal coordinates $[ \alpha, \beta, h]$ and geometric local toroidal coordinates $[ \varepsilon, \theta, h]$}
    \label{fig:toroidal_frame}
\end{figure}

A local eigenspace for an oscillatory mode at a fixed point in a reference orbit is illustrated in Fig. \ref{fig:toroidal_frame}.
An osculating invariant circle is expressed as \cite{olikara2016computation}
\begin{align} \label{eq:inv_curve_first_order_approx}
    \boldsymbol{\psi}(t) = \varepsilon \left(\boldsymbol{w}_r(t) \cos \theta + \boldsymbol{w}_i(t) \sin \theta \right), \quad \boldsymbol{w}_r = \begin{bmatrix}
        \boldsymbol{r}_r(t) \\
        [\dot{\boldsymbol{r}}_r(t)]_\mathcal{V} 
    \end{bmatrix} 
    , \boldsymbol{w}_i = 
    \begin{bmatrix}
        \boldsymbol{r}_i(t) \\
        [\dot{\boldsymbol{r}}_i(t)]_\mathcal{V}  
    \end{bmatrix}, 
\end{align}
where $\boldsymbol{w}(t) = \boldsymbol{w}_r(t) + \boldsymbol{w}_i(t) i$ is the eigenvector corresponding to the oscillatory mode at the fixed point. 
A scaling parameter $\varepsilon$ defines the size of the invariant circle, and $\theta$ (cf. equivalent to $\theta_1$ in Fig.  \ref{fig:qpt}) defines a phase within the invariant circle. 
Note that although the eigenvalues of the monodromy matrix are invariant along a periodic orbit, the associated eigenvectors depend on the selected fixed point.

With this in mind, the local toroidal frame is defined based on a time-varying non-orthogonal and non-unimodular basis triad $\{\boldsymbol{r}_r(t), \boldsymbol{r}_i(t), \hat{\boldsymbol{n}}(t)\}$, where $\hat{\boldsymbol{n}}(t)$ is a unimodular vector parallel to $ \boldsymbol{r}_r(t) \times \boldsymbol{r}_i (t)$.
In order to achieve a unique and consistent representation of an eigenvector and to construct a triad from a local eigensystem, basis vectors at the fixed point at the initial epoch $\{ \boldsymbol{r}_r(t_0), \boldsymbol{r}_i(t_0) \}$ are normalized and rotated via singular value decomposition (SVD) so that these vectors are orthogonal and $\boldsymbol{r}_r(t_0)$ becomes unimodular. 
For details of the normalization process, see Ref. \citenum{elliott2022phd}.
Once normalized at the fixed point, the evolution of $\boldsymbol{w}$ can be propagated by the linear dynamics 
\begin{align} \label{eq:w_dot_dyn}
    [\dot{\boldsymbol{w}}(t)]_\mathcal{V} = \boldsymbol{A}(t) \boldsymbol{w}(t),
\end{align}
where the plant matrix is obtained from Eq. \eqref{eq:plant_mat_TNW}, and the time derivative is taken in the local frame. 
Therefore, the basis vectors $\{\boldsymbol{r}_r(t), \boldsymbol{r}_i(t), \hat{\boldsymbol{n}}(t)\}$ are generally neither orthogonal (nor unimodular for $\{\boldsymbol{r}_r, \boldsymbol{r}_i\}$) except for the initial triad at the fixed point. 

Within the local toroidal frame, the relative position of the deputy to the chief is written as $\boldsymbol{\rho} = \alpha \boldsymbol{r}_r +  \beta \boldsymbol{r}_i + h \hat{\boldsymbol{n}}$.
Based on this, nonsingular LTC is defined as $\boldsymbol{\zeta}_{\text{ns}} := [ \alpha, \beta, h, \dot{\alpha}, \dot{\beta}, \dot{h} ]^\top$, where time derivatives are taken in the local (toroidal) frame.  
Using Eq. \eqref{eq:dx_to_integ_const}, the mapping from the nonsingular LTC $\boldsymbol{\zeta}_{\text{ns}}$ to the relative position and velocity resolved in the local frame $\boldsymbol{x}$ is related as
\begin{align} \label{eq:map_LTW_to_cart}
    \boldsymbol{\zeta}_{\text{ns}}
    = T(t)^{-1} 
    \boldsymbol{x}, \quad T(t) := 
    \begin{bmatrix}
        R(t) & \boldsymbol{0}_3 \\
        R'(t) & R(t) 
    \end{bmatrix}  \Leftrightarrow 
    T^{-1} (t) = \begin{bmatrix}
        R(t)^{-1} & \boldsymbol{0}_3 \\
        - R(t)^{-1} R'(t) R(t)^{-1} & R(t)^{-1}
    \end{bmatrix}, 
\end{align}
where $R(t) := [\boldsymbol{r}_r(t), \boldsymbol{r}_i(t), \hat{\boldsymbol{n}}(t)]$ and 
$R'(t) := [[\dot{\boldsymbol{r}}_r(t)]_\mathcal{V}, [\dot{\boldsymbol{r}}_i(t)]_\mathcal{V}, [\dot{\hat{\boldsymbol{n}}}(t)]_\mathcal{V}]$. 
The time derivative of the unit normal vector $\hat{\boldsymbol{n}}(t)$ is computed as 
\begin{align}
    [\dot{\hat{\boldsymbol{n}}}(t)]_\mathcal{V} = [\dot{\boldsymbol{r}_r}(t)]_\mathcal{V} \times \boldsymbol{r}_i(t) + \boldsymbol{r}_r(t) \times [\dot{\boldsymbol{r}_i}(t)]_\mathcal{V}  - (\boldsymbol{n}(t)^\top ([\dot{\boldsymbol{r}_r}(t)]_\mathcal{V}  \times \boldsymbol{r}_i(t) + \boldsymbol{r}_r(t) \times [\dot{\boldsymbol{r}_i}(t)]_\mathcal{V} )) \boldsymbol{n}(t).
\end{align}
Note that $T(x)$ serves as a fundamental matrix solution to the variational equation in the relative motion of RMBPs, similarly to Eq. \eqref{eq:dx_to_integ_const}.
Moreover, in the same fashion to Eq. \eqref{eq:monod_transform}, the State Transition Matrix (STM) of the LTC is expressed as 
\begin{align} 
\Phi_{\boldsymbol{\zeta}_{\text{ns}}}(t+\Delta t, t) = T^{-1}(t +\Delta t) \Phi(t+\Delta t, t) T(t), 
\end{align}

One can also define LTC by expressing the eigenspace component in polar coordinates (cf. Eq. \eqref{eq:inv_curve_first_order_approx}) as $\boldsymbol{\zeta}_{\text{geo}} := [\varepsilon, \theta, h, \dot{\varepsilon}, \dot{\theta}, \dot{h}]^\top$, which is referred to as geometric LTC.
This approximation holds for the relative motion around any periodic orbit with oscillatory mode thanks to the linearization adopted in Ref.~\cite{franzini_relative_2019}, and the validity of its accuracy is directly affected by the size and direction of the invariant circle (i.e., eigenspace) at each phase of the reference orbit.
A similar mapping from the geometric LTC to Cartesian coordinates is available \cite{elliott2022phd}, although the transformation is no longer a function of only the states on a reference orbit but also a time derivative of nonsingular LTC. 
Nonetheless, the geometric LTC provides a convenient characterization of an invariant circle's size, thereby facilitating the specification of boundary conditions and waypoints for trajectory design.

The first-order condition to be on a QPRIT (QPRO) is to stay on the eigenspace spanned by $\boldsymbol{w}_r(t)$ and $\boldsymbol{w}_i(t)$. 
This leads to the condition $h = \dot{\alpha} = \dot{\beta} = \dot{h} = 0$, or using geometric LTC, $h = \dot{\varepsilon} = \dot{\theta} = \dot{h} = 0$ \cite{elliott2022describing}. 
While other coordinates can still describe the deputy’s relative motion, the (first-order) quasi-periodicity linked to the oscillatory mode specific to the toroidal coordinate set cannot be maintained.

The quasi-periodic condition allows for arbitrary choices of the variables $(\alpha, \beta)$ or $(\varepsilon, \theta)$, defining both the size of the invariant circle and its phase.
However, nonzero values in $(\dot{\alpha}, \dot{\beta})$  or $(\dot{\varepsilon}, \dot{\theta})$ result in a drift in the $h$-component, breaking the quasi-periodicity. 

\subsection{Safety evaluation of QROs}

This subsection demonstrates the connection between the ROE and LTC based on the geometric configuration of the osculating relative orbit and the minimum separation. 
Although the approach is general, the application case of this section is on the Earth-Moon $L_2$ 9:2 South NRHO, which is the orbit of the Lunar Gateway \cite{lee2019gateway, nakamura2023rendezvous} and has an orbital period of 6.56 days. 

\subsubsection{Comparison between the LVLH and VNB} frame

\begin{figure}[ht!]
     \centering
     \begin{subfigure}[b]{0.45\textwidth}
         \centering
         \includegraphics[width=\textwidth]{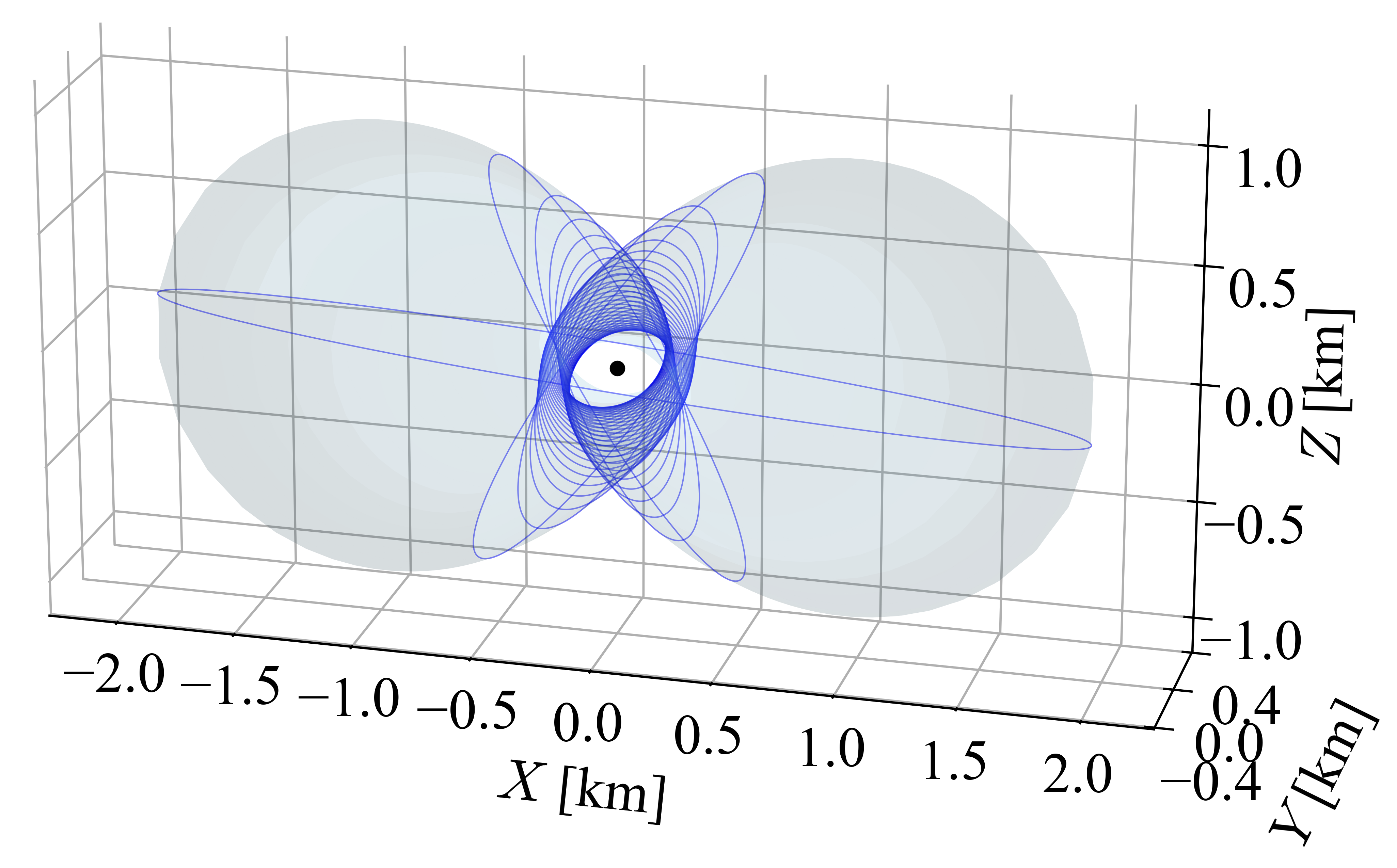}
         \caption{LVLH frame}
         \label{fig:qpt_lvlh}
     \end{subfigure}
     \begin{subfigure}[b]{0.45\textwidth}
         \centering
         \includegraphics[width=\textwidth]{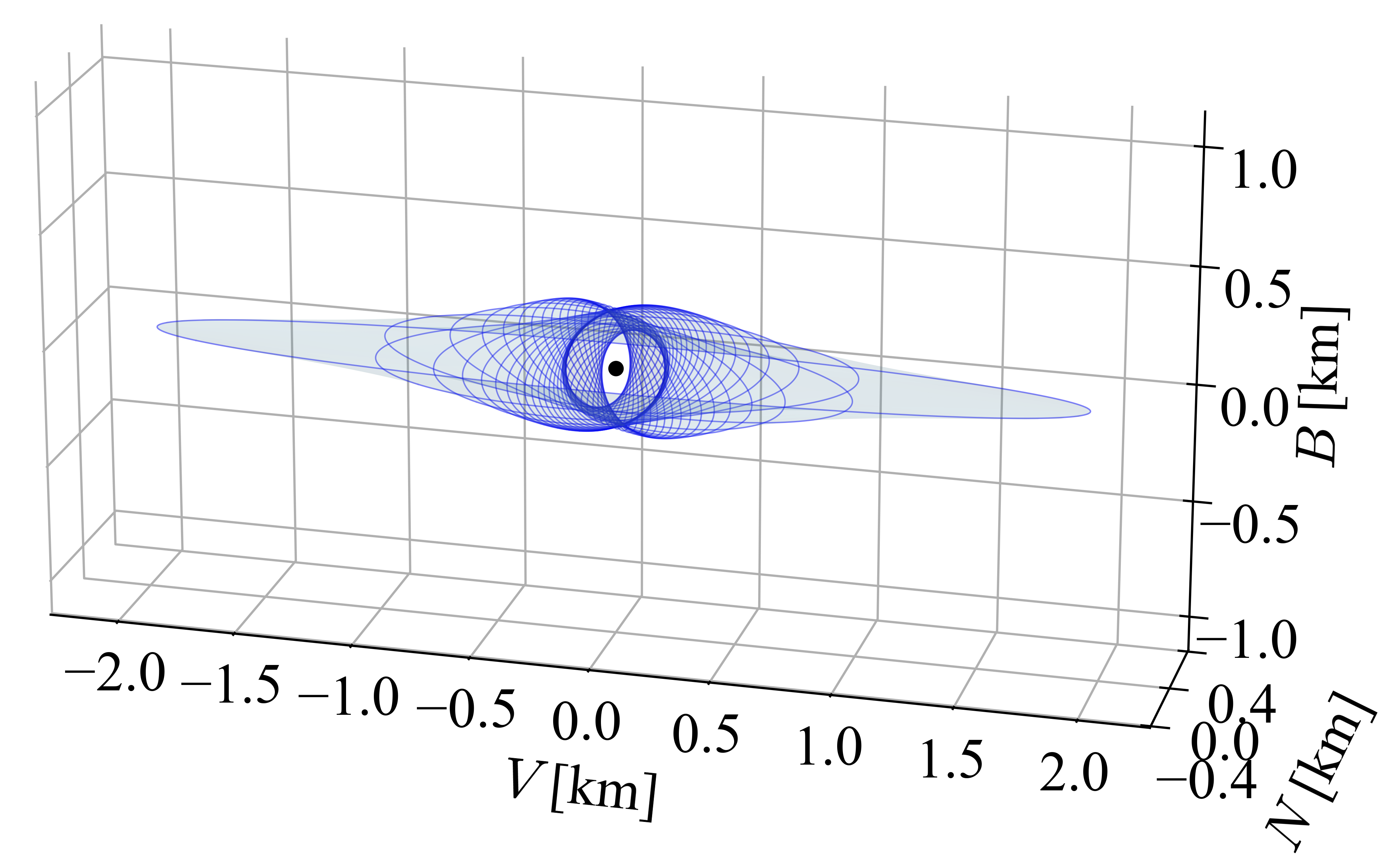}
         \caption{VNB} frame
         \label{fig:qpt_tnw}
     \end{subfigure}  
    \caption{First-order QPRIT (gray surface) and 50 invariant circles (blue ellipses) of CR3BP $L_2$ South 9:2 (synodic) NRHO in co-moving frames ($\varepsilon$=200 m, fixed point at apolune).}
    \label{fig:qpt_3d_comparison} 
\end{figure}

The first-order approximation of a QPRIT resolved in the LVLH and VNB frames is illustrated in Figs. \ref{fig:qpt_lvlh} and \ref{fig:qpt_tnw}, respectively. 
Both figures depict the same QPRIT, where the fixed point is anchored at apolune with a scaling factor $\varepsilon = 200$ m. 
The surfaces of the QPRIT are generated by the continuation of invariant circles (cf. Eq. \eqref{eq:inv_curve_first_order_approx}) along the reference orbit, shown as a gray surface. 
Additionally, fifty osculating invariant circles, sampled uniformly over one orbital period, are shown as blue ellipses. 
In the VNB frame, the semimajor axes of the invariant circles align with the chief’s velocity direction, producing a more compact structure. 
In contrast, for highly eccentric orbits like NRHOs, the angular offset of the QPRIT in the LVLH frame is evident.
Previous works \cite{elliott2022describing, elliott2022phd, down2023relative} report a QPRIT structure of NRHO characterized by two disk-like surfaces for a wide range of period, a pattern arising from positional displacement relative to the spacecraft velocity in the synodic frame. 
The reduced density of invariant circles within these disks reflects higher relative velocities and greater radial uncertainty.

\subsubsection{Geometric Interpretation of Osculating Relative Invariant Circle}

\begin{figure}[ht!]
     \centering
     \begin{subfigure}[b]{0.7\textwidth}
         \centering
         \includegraphics[width=\textwidth]{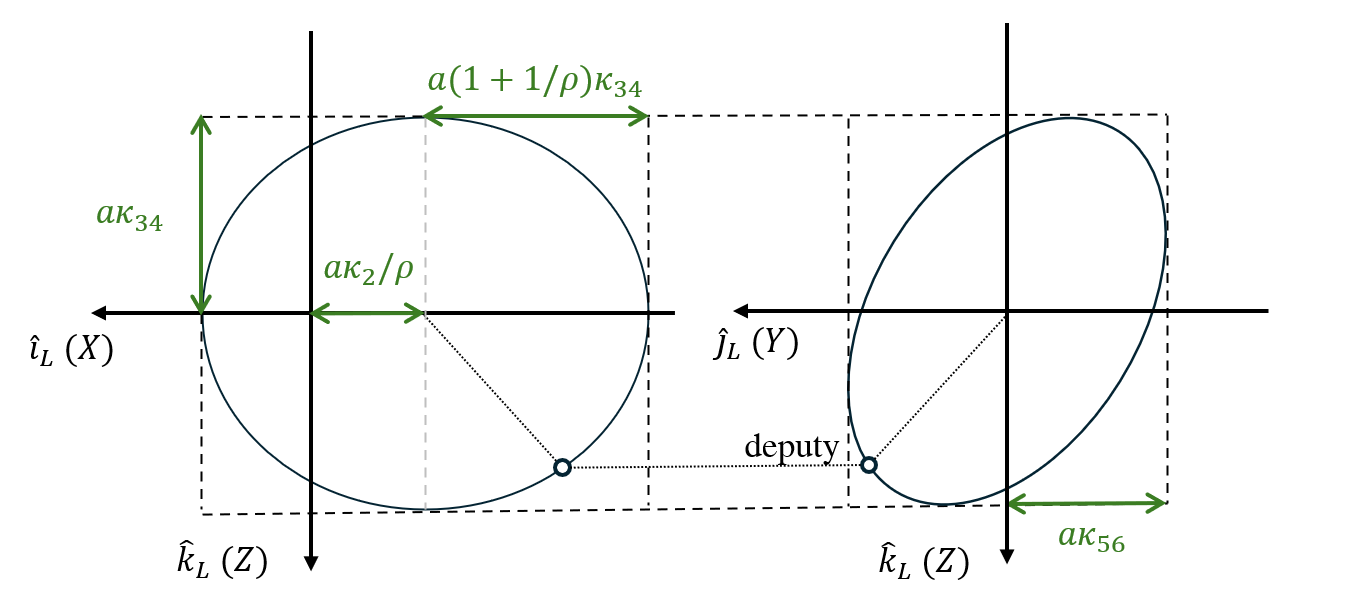}
         \caption{Osculating relative ellipse of a relative orbit (PR2BP) in the LVLH frame ($\kappa_1$=0).}
         \label{fig:ellipse_lvlh}
     \end{subfigure}
     \\
     \begin{subfigure}[b]{0.7\textwidth}
         \centering
         \includegraphics[width=\textwidth]{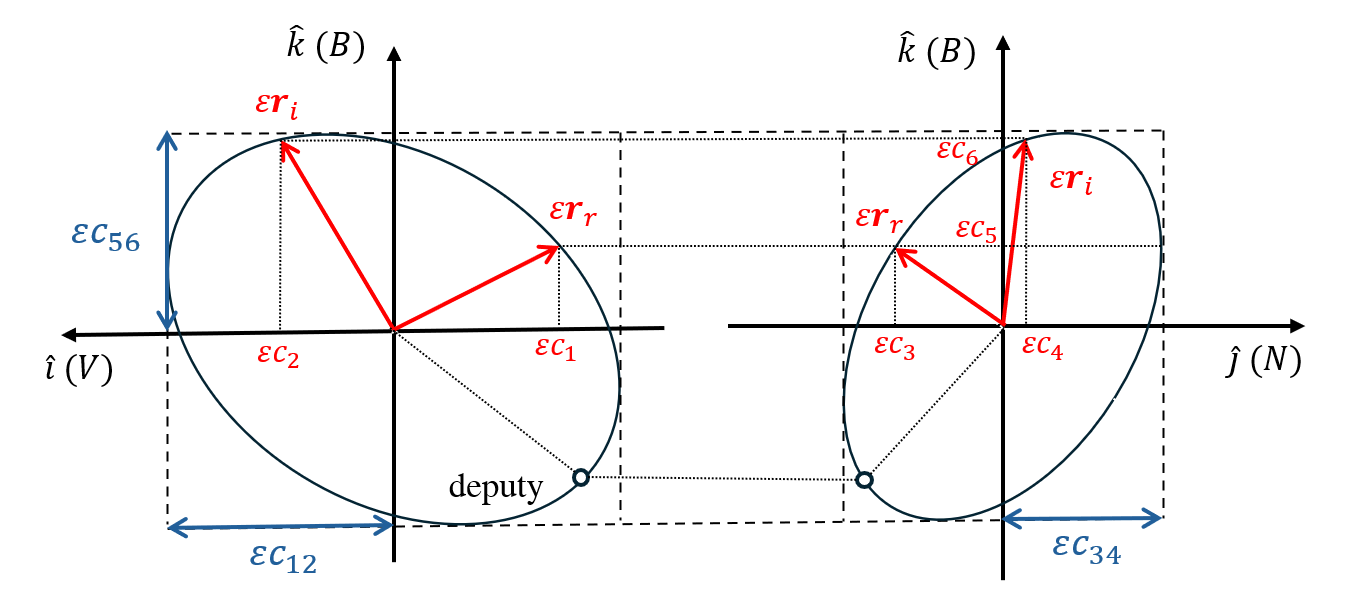}
         \caption{Osculating invariant circle of a QPRO (RMBP) in the VNB} frame based on the parametrization in Eq.~\eqref{eq:psi_amplitude}.
         \label{fig:ellipse_tnw}
     \end{subfigure}  
    \caption{Comparison of the geometric parametrization of osculating relative orbit in the PR2BP and the osculating invariant circle in the RMBPs.}
    \label{fig:ellipse_projection} 
\end{figure}

Two-dimensional projections of an osculating relative orbit in the PR2BP and an osculating invariant circle in the RMBPs are shown in Figs. \ref{fig:ellipse_lvlh} and \ref{fig:ellipse_tnw}, respectively.
For a relative motion around a Keplerian orbit, a six-dimensional variable $\boldsymbol{\kappa} = \in \mathbb{R}^6$ parameterizes a relative orbit of the deputy with respect to the chief (cf. Eq.~\eqref{eq:dx_to_integ_const}), with additional independent variable $\rho(t) = 1 + e\cos \nu(t)$, where $e$ and $\nu(t)$ denote the eccentricity and true anomaly of the chief, respectively. 
In this paper, quasi-nonsingular ROE in near-circular orbits \cite{damico2010autonomous} and its extension to elliptic orbits \cite{delurgio2024closed} are chosen for $\boldsymbol{\kappa}$, which correspond to integration constants of the linearized equations of motion \cite{clohessy_terminal_1960, yamanaka_new_2002}.
In Fig.~\ref{fig:ellipse_lvlh}, $\kappa_1 = 0$ is enforced to realize a non-drifting bounded relative motion.

The position-components of the osculating invariant circle in Eq. \eqref{eq:inv_curve_first_order_approx} can be parameterized within the local frame as 
\small
\begin{align} \label{eq:psi_amplitude}
    & \boldsymbol{\psi}_p(t) = \varepsilon 
        \left(
        \sqrt{c_1^2(t) + c_2^2(t)} \sin (\theta + \phi_x(t)) \hat{\boldsymbol{\imath}}
        +
        \sqrt{c_3^2(t) + c_4^2(t)} \sin (\theta + \phi_y(t)) \hat{\boldsymbol{\jmath}}
        +
        \sqrt{c_5^2(t) + c_6^2(t)} \sin (\theta + \phi_z(t)) \hat{\boldsymbol{k}}
    \right) 
    \\
    & 
    \tan \phi_x(t) = c_1(t) / c_2(t), \quad \tan \phi_y(t) = c_3(t) / c_4(t), \quad \tan \phi_z(t) = c_5(t) / c_6(t)
    \\
    & \boldsymbol{r}_r(t) = 
        c_1(t)\hat{\boldsymbol{\imath}} +  c_3(t)\hat{\boldsymbol{\jmath}} + c_5(t) \hat{\boldsymbol{k}},
    \quad 
    \boldsymbol{r}_i(t) = 
        c_2(t)\hat{\boldsymbol{\imath}} + c_4(t)\hat{\boldsymbol{j}} + c_6(t)\hat{\boldsymbol{k}}    
\end{align}
\normalsize
\noindent
Vectors $[c_1(t), c_2(t)]$, $[c_3(t), c_4(t)]$, and $[c_5(t), c_6(t)]$ provide the amplitudes of the projected ellipse in each plane as well as the oscillation phase. 
As stated in Sec.~\ref{sec:eigen_2BP}, $\boldsymbol{\kappa}$ is time-invariant\footnote{When perturbations (e.g., spherical harmonics, drag, etc.) are added to the Keplerian dynamics, the osculating ROE also experiences a drift, although it does not induce any quasi-periodic motion \cite{damico2010autonomous, schaub2003analytical, koenig2017new, delurgio2024closed}}. 
Nonetheless, the magnitude of the perturbation on the eigensystem is negligible (i.e., the eigenvectors of the monodromy matrix are almost 1) and can be treated as a secular drift in an averaged motion over multiple revolutions in a phase space.  
in the R2BP, allowing the explicit expression of the minimum planar separation in each projected plane \cite{damico2010autonomous}.
In contrast, vector $\boldsymbol{c}(t) \in \mathbb{R}^6$ is time-varying in RMBPs, where a combined analysis of the magnitudes and phases of this vector is necessary to assess the passive safety of a QPRO.

\begin{figure}[ht!]
     \centering
     \begin{subfigure}[b]{0.8\textwidth}
         \centering
         \includegraphics[width=\textwidth]{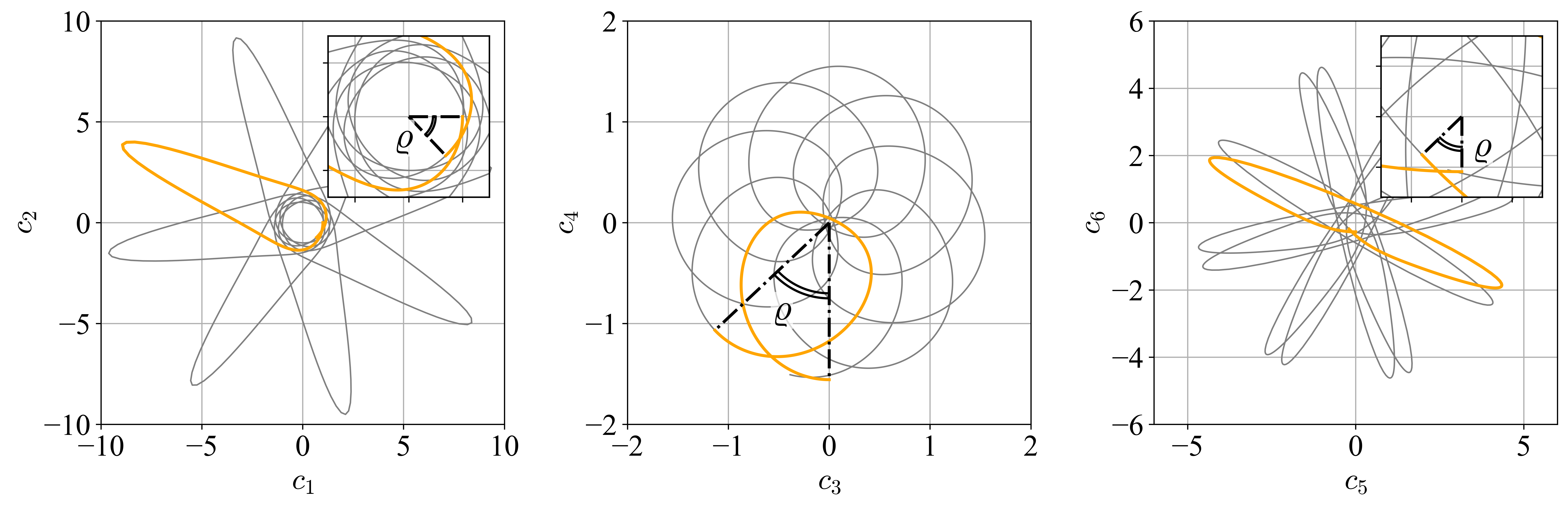}
         \caption{LVLH frame}
         \label{fig:inv_circ_lvlh}
     \end{subfigure}
    \\
     \begin{subfigure}[b]{0.8\textwidth}
         \centering
         \includegraphics[width=\textwidth]{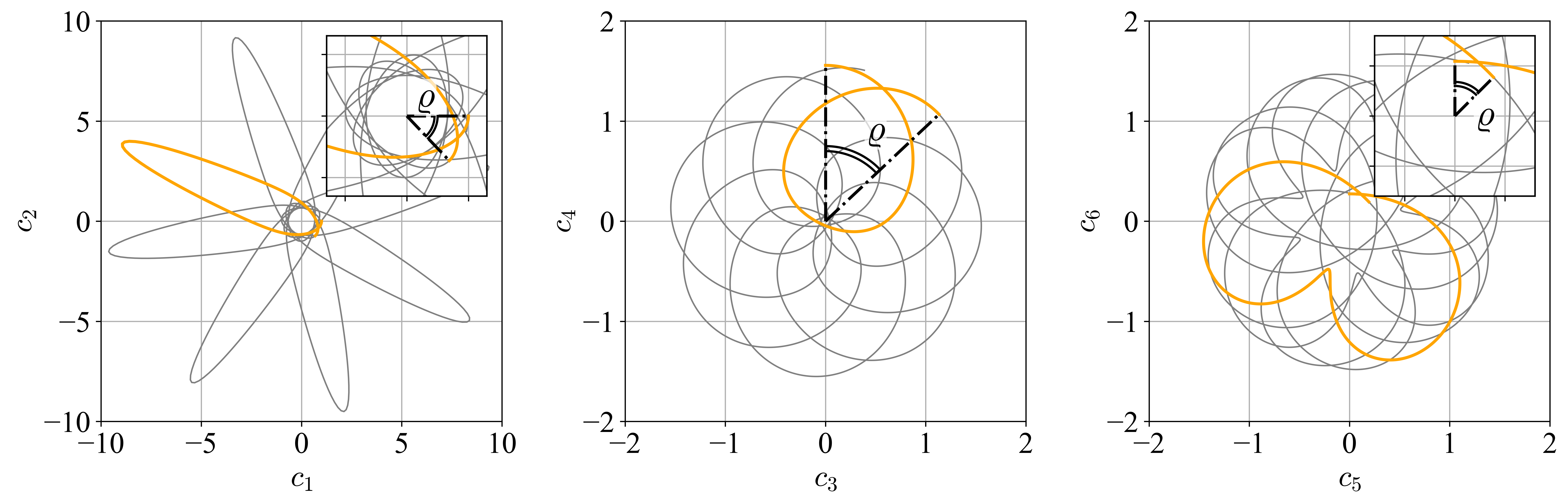}
         \caption{VNB frame}
         \label{fig:inv_circ_tnw}
     \end{subfigure}  
    \caption{ Time history of the components of the eigenvectors $\boldsymbol{r}_r(t)$ and $\boldsymbol{r}_i(t)$ corresponding to the unit osculating invariant circle. Orange curves correspond to the first revolution, where each vector rotates by the rotation number $\varrho$.} 
    \label{fig:inv_circ_amplitude} 
\end{figure}

Fig. \ref{fig:inv_circ_amplitude} presents the evolution of $\boldsymbol{c}(t)$ over eight revolutions, resolved in both the LVLH and VNB frames. 
The history is computed using a constant size $\varepsilon=1$ for the osculating invariant circle, ensuring preservation of the geometric structure.
The norm of each vector field represents the projected length of the unit osculating invariant circle along each axis of the co-moving frame.
Fig.~\ref{fig:inv_circ_amplitude} is analogous to the plots of the vector fields of relative eccentricity vector $[\kappa_3, \kappa_4]$ and relative inclination vector $[\kappa_5, \kappa_6]$ in the quasi-nonsingular ROE \cite{damico2010autonomous, guffanti2022phd} because the magnitude of each vector represents the amplitude of osculating invariant circle in each direction of the axis ($c_{12}, c_{34}$, and $c_{56}$), and the magnitudes of the two vectors and the phase difference between them determine the minimum separation in the corresponding two-dimensional projected plane.
However, note again that the vectors $[c_1(t), c_2(t)]$, $[c_3(t), c_4(t)]$, and $[c_5(t), c_6(t)]$ for a unit invariant circle are determined by the local eigensystem, whereas the relative eccentricity and inclination vectors are free design variables of the deputy.

From each vector field illustrated in Fig. \ref{fig:inv_circ_amplitude}, a quasi-periodic behavior is confirmed with a drastic change in its amplitude. 
Note that the equivalent vectors in relative motion for the PR2BP (i.e., ROE) would be stationary (for a near-circular orbit) or oscillatory about a line segment (for an elliptic orbit), with possible slow secular drifts due to perturbations.
The evolution of the eigenvector can be characterized by two angles $ [\theta_0, \theta_1]$ by separating the magnitude and the phase as $\boldsymbol{w}(t) = w_r(\theta_0(t)) \hat{\boldsymbol{w}}_r(\theta_1(t)) + w_i(\theta_0(t)) \hat{\boldsymbol{w}}_i(\theta_1(t)) i $.
The first angle, $\theta_0$, determines the magnitude of the vector corresponding to the fixed point on the reference orbit. 
In contrast, the second angle, $\theta_1$, captures the secular drift of the vector over multiple orbital revolutions, reflecting the geometric structure of the two-dimensional QPIT and QPRIT (cf. Fig.~\ref{fig:qpt}).

The angular rate of this secular rotation $\dot{\theta}_1$ remains constant over the three vector fields and corresponds to the angle of the complex eigenvalue of the oscillatory mode, known as rotation number $\varrho$ \cite{olikara2016computation, baresi2017spacecraft}. 
For example, an eigenvalue associated with the oscillatory mode of the monodromy matrix of the 9:2 South NRHO is $\lambda \simeq 0.6845-0.7290i$. 
This leads to the rotation number of $\varrho \simeq -46.80^\circ$, which coincides with the angular displacement between the initial and terminal points of the orange segment in Fig. \ref{fig:inv_circ_amplitude}.
For both Figs. \ref{fig:inv_circ_lvlh} and \ref{fig:inv_circ_tnw}, the minimum amplitude in the N-direction collapses to nearly zero, as shown in the history of $[c_3(t), c_4(t)]$. 
A significant difference between Figs. \ref{fig:inv_circ_lvlh} and \ref{fig:inv_circ_tnw} lies in the vector $[c_5(t), c_6(t)]$, which represents the amplitude of the invariant circle in the R(W)-direction. 
This is illustrated by the contraction of the amplitude in the W-direction in the VNB frame compared to the R-direction in the LVLH frame.
\begin{figure}[ht!]
     \centering
     \begin{subfigure}[b]{0.49\textwidth}
         \centering
         \includegraphics[width=\textwidth]{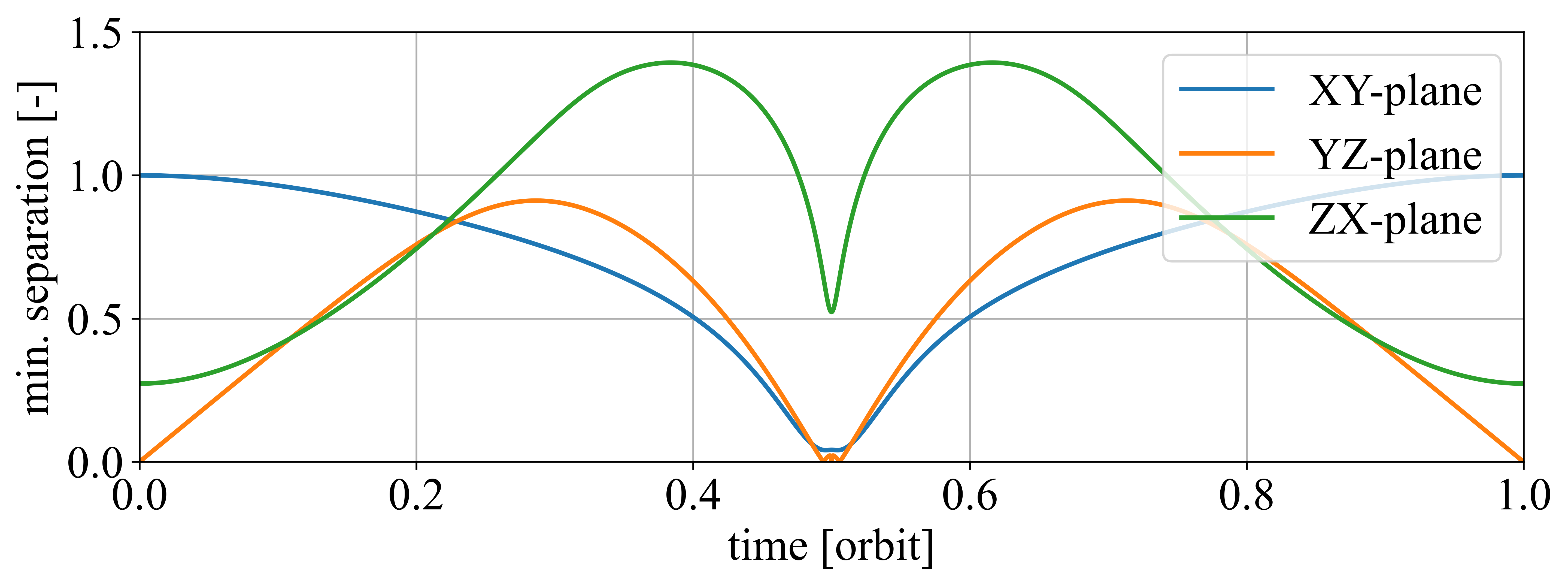}
         \caption{LVLH frame}
         \label{fig:minsep_LVLH}
     \end{subfigure}
     \begin{subfigure}[b]{0.49\textwidth}
         \centering
         \includegraphics[width=\textwidth]{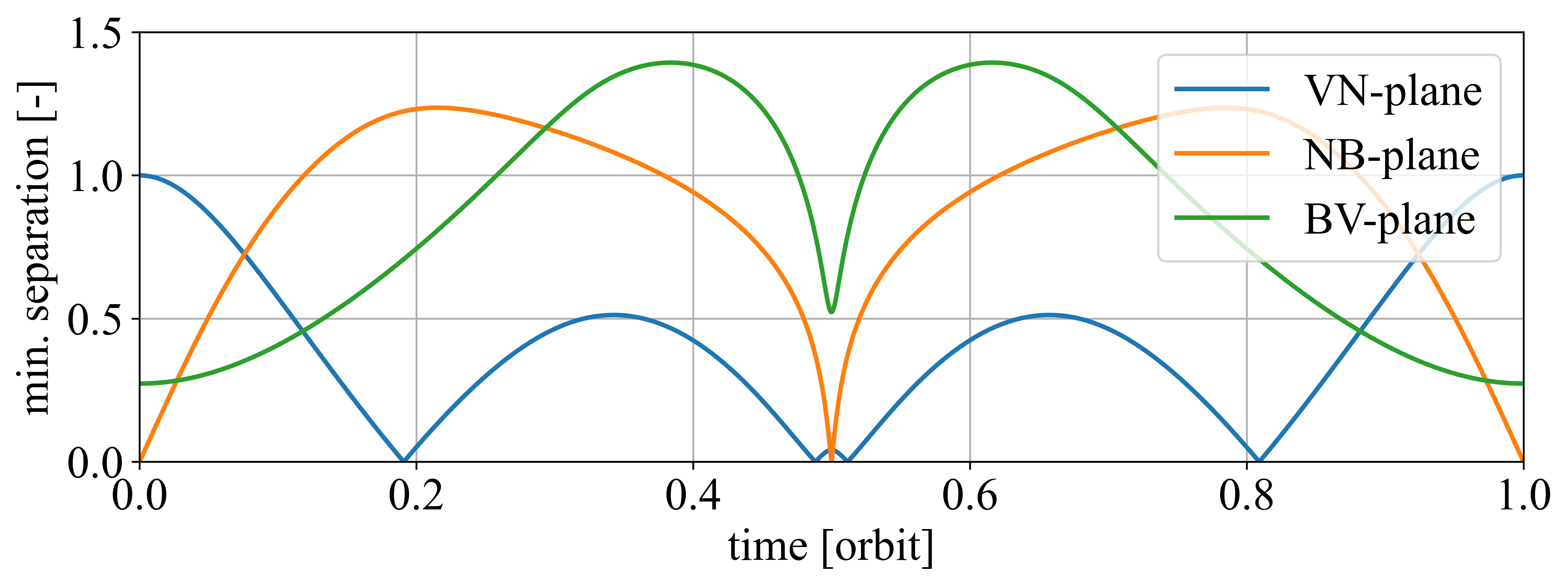}
         \caption{VNB frame}
         \label{fig:minsep_TNW}
     \end{subfigure}  
    \caption{Evolution of the minimum planar separation of a unit invariant circle. The initial point is set to apolune. }
    \label{fig:minsep_planar} 
\end{figure}

Previous studies provide analysis on the minimum positional separation by performing an SVD \cite{elliott2022describing} or modal decomposition \cite{ vela2025modal}. 
Applying SVD to the eigenvectors projected onto the two-dimensional co-moving frame allows for the evaluation of the minimum planar separation, which offers additional insights into the design of the safety region. 
Fig.~\ref{fig:minsep_planar} illustrates the evolution of the minimum separation in each plane of the LVLH frame and the VNB frame along the full orbit, starting from the apolune. 
Only the ZX/BV-plane has a nonzero separation for both frames, whereas the separations in the XY/VN- and YZ/NB-planes collapse to zero at perilune (0.5 orbit).
Furthermore, the separation in the NB plane increases by resolving the QPRIT in the VNB frame, whereas additional singularities are generated in the VN- plane separation. 
Note that although the separation increases with a larger QPIT radius ($\varepsilon$), this trend is specific to the orbit. 
In the PR2BPs, an appropriate relative orbit configuration can always secure planar separation in the YZ plane, which is known as relative eccentricity/inclination-vector separation \cite{montenbruck2006GRACE,damico2010autonomous}. 
This also motivates a unique passive safety formulation for the relative motion in the RMBPs.

\section{Passively Safe Reconfiguration using Local Toroidal Coordinates} \label{sec:passive_safety_ocp}

This section presents a passively safe, fuel-optimal swarm reconfiguration method in cislunar relative motion. 
Assuming that the chief is on a reference orbit, the focus is on transferring an $N_{sc}$ deputy spacecraft between two passively safe QPROs, where the initial and terminal conditions are defined by $\boldsymbol{\zeta}_{\text{ns},i}^{(j)}$ and $\boldsymbol{\zeta}_{\text{ns},f}^{(j)}$ for each $j = 1, \dots, N_{sc}$, with a fixed transfer time $t_f$ while ensuring passive safety throughout the transfer. 
Such safe transfers are crucial for missions involving on-orbit inspection or formation flying.

To facilitate numerical optimization and enable efficient trajectory planning under mission constraints, the problem is cast in a discrete-time framework. 
This approach also aligns with the use of linearized dynamics models and impulsive control inputs.
With a sufficiently fine time discretization, this problem can be formulated as a discrete-time OCP using the linearized dynamics on the nonsingular LTC within a local frame as follows: 
\small
\begin{subequations}  \label{eq:tpbvp}
\begin{align}
    \min_{ \{\boldsymbol{\zeta}_{\text{ns},k}^{(j)}\}_{k=0}^N, \{\boldsymbol{u}_k^{(j)}\}_{k=0}^N} \ & 
     \sum_{j=1}^{N_{sc}} \mathcal{J}_j = \sum_{j=1}^{N_{sc}} \sum_{k=0}^{N} \| \boldsymbol{u}_k^{(j)} \| \\
    \text{subject to} \quad & \boldsymbol{\zeta}_{\text{ns},i}^{(j)} = \boldsymbol{\zeta}_{\text{ns},0}^{(j)}, \quad \boldsymbol{\zeta}_{\text{ns},f}^{(j)} = \boldsymbol{\zeta}_{\text{ns},N}^{(j)} + B_N\boldsymbol{u}_N, & \forall j  \\
    & \boldsymbol{\zeta}_{\text{ns},k+1}^{(j)} = A_k \left(\boldsymbol{\zeta}_{\text{ns},k}^{(j)} + B_{k} \boldsymbol{u}_{k}^{(j)}\right) & k = 0, ..., N-1, \forall j \\
    & g_{ps}\left(\boldsymbol{\zeta}_{\text{ns},k}^{(j)}, \boldsymbol{\zeta}_{\text{ns},k}^{(j')}, \right) \leq 0 , &  k = 0, ..., N, j = 0,..., N_{sc}, j' = 0, ..., N_{sc}, j\neq j' \label{eq:ps_constr_general} \\
    & \boldsymbol{u}_k^{(j)} = \boldsymbol{0}, & k \in [k_{p}^{-},k_{p}^{+}], \forall j  \label{eq:no_control} \\
    \text{where} \quad & A_k = T_{k+1}^{-1} \Phi(t_{k+1}, t_k) T_k, 
    \quad 
    B_k = \begin{bmatrix}
        \boldsymbol{0}_3 \\
        R^{-1}_k
    \end{bmatrix} \label{eq:stm_cim}.
\end{align}
\end{subequations}
\normalsize
The control input $\boldsymbol{u}_k^{(j)}$ is an impulsive velocity change added to the $j$-th deputy within the VNB frame.
The minimum-fuel problem is equivalent to the minimization of the sum of the $l_2$-norm of impulsive maneuvers.
Using Eq. \eqref{eq:map_LTW_to_cart}, the control input matrix $B_k$ is derived from the last three columns of $T_k^{-1}$.
Eq. \eqref{eq:ps_constr_general} presents the passive safety constraint that is applied to two spacecraft $\boldsymbol{\zeta}_{\text{ns},k}^{(j)}$ and $\boldsymbol{\zeta}_{\text{ns},k}^{(j')}$, where $j=0$ denotes the chief. 
Two passive safety formulations are discussed in the following subsections.

To avoid control near perilune, where system sensitivity is high and the linearized dynamics incur substantial nonlinear errors, Eq.~\eqref{eq:no_control} enforces a no-control condition. 
These no-control windows, denoted by the heuristic intervals $k \in [k_p^-,k_p^+]$, define time segments during which control actions are not allowed.

Because the reference orbit may be highly eccentric (e.g., NRHO), a time regularization scheme is considered to distribute the discrete-time control nodes to reduce the nonlinear errors emanating from the dynamics linearization.
The time discretization $dt$ is regularized to a pseudo-time discretization $d\tau$ by a scaling function of the chief's position and velocity as $g(\boldsymbol{y})$ as $dt = g(\boldsymbol{y}) d\tau  \ \Leftrightarrow \ dt / d\tau = g(\boldsymbol{y})$. 
This leads to the reformulation of the dynamics in the pseudo-time domain as 
\begin{align}
        \dfrac{d\boldsymbol{y}}{d\tau} = \dfrac{d\boldsymbol{y}}{dt} \dfrac{dt}{d\tau}  = g (\boldsymbol{y}) \dot{\boldsymbol{y}}  
        , \quad
        \dfrac{dt}{d\tau} = g(\boldsymbol{y}).
\end{align}
A set of regularized time epochs can be determined by propagating the above dynamics over the pseudo-time interval $[\tau_0, \tau_f]$, where the pseudo-time epochs are distributed uniformly. 
The retrieved time epochs that correspond to the pseudo-time epochs represent a vector of the regularized time epochs. 
In this paper, the Sundman transformation for the two-body system \cite{sundman1913memoire} with respect to the Moon is employed, considering that NRHOs are close to an elliptic orbit. 
The regularization function is defined as $g(\boldsymbol{y}) = r^\alpha$, where $\alpha=1$ is a user-defined parameter.

\subsection{Passive Safety via Drift Trajectories From Controlled Nodes}

One way to achieve passive safety is to enforce the exclusion of the drift trajectory from the keep-out zone for a predefined period \cite{breger2008safe, guffanti2023passively, takubo2024art}. 
The geometric configuration of an ellipsoidal KOZ is characterized by the diagonal matrix $P$. 
Using this matrix, the exclusion from the KOZ is expressed as a quadratic inequality $\sqrt{{\boldsymbol{x}}^\top P \boldsymbol{x}} \geq 1$.
In this paper, the KOZ exclusion constrained is defined within the VNB frame, where the ellipsoid matrix is defined as $\text{diag}([1/x_{e}^2, 1/y_{e}^2, 1/z_{e}^2,0,0,0])$ and $(x_{e}, y_{e}, z_{e})$ are the semi-axes of the ellipsoid in the VNB frame.
Therefore, the passive safety of the $j$-th deputy and the chief is rewritten as \cite{morgan2014model} 
\begin{subequations}
\label{eq:drifty_safety}
\begin{align}
    & \max_{m = 0, ..., N_{\text{safe}}} \sqrt{{\boldsymbol{x}_{km}^{(j)}}^\top P \boldsymbol{x}_{km}^{(j)}}  -{\boldsymbol{x}_{km}^{(j)}}^\top P \boldsymbol{x}_{km}^{(j)} \leq 0, \quad k = 0, ..., N, \label{eq:drifty_safety1} \\
    \text{where} \quad &  \boldsymbol{x}_{km}^{(j)} = \prod_{l=0}^{m-1} \Phi(t_{l+1}, t_l) T_k \boldsymbol{\zeta}_{\text{ns},k}^{(j)}, \label{eq:drifty_safety2}
\end{align}
\end{subequations}
where Eq.~\eqref{eq:drifty_safety2} defines the natural propagation of the deputy after each controlled epoch, and Eq.~\eqref{eq:drifty_safety1} ensures the exclusion of the drift trajectory from the KOZ at each discretized time step. 
In this subsection, for simplicity of notation, $\boldsymbol{x}$ represents the matrix representation of a state vector within the local frame.
Note that this is an equivalent formulation to the well-known quadratic form of the exclusion constraint 
To prevent the quadratic growth of the number of constraints along the drift period, the maximum operator is applied over the interval $m = 0, ..., N_{\text{safe}}$. 
Note that an instantaneous collision avoidance is achieved by considering only $m=0$.
Since Eq. \eqref{eq:drifty_safety1} is nonconvex, the following convexification is performed with respect to the reference orbit resolved in a local frame $\bar{\boldsymbol{x}}_k$ as follows \cite{morgan2014model, guffanti2023passively}:
\small
\begin{subequations} \label{eq:ps_lin}
\begin{align} 
    & -\boldsymbol{a}^\top_{k} T_k \boldsymbol{\zeta}^{(j)}_{\text{ns},k} + b_k \leq 0, \quad k = 0, ..., N \\
    \text{where} \quad & \boldsymbol{a}^\top_k =  {\bar{\boldsymbol{x}}}_k^{{(j)}^\top} \prod_{l=0}^{m^*_k-1} \Phi(t_{l+1}, t_l) D^\top P D \left( \prod_{l=0}^{j^*_k-1} \Phi(t_{l+1}, t_l) \right)^\top, 
    \quad 
    b_k = \sqrt{ \boldsymbol{a}^\top_k \bar{\boldsymbol{x}}_k^{(j)}}, \quad D= [ \boldsymbol{I}_3, \boldsymbol{0}_3], \\
    & m^*_k = \argmax_{m = 0, ..., N_{\text{safe}}} \sqrt{{\bar{\boldsymbol{x}}}_{km}^{{(j)}^\top} P \bar{\boldsymbol{x}}^{(j)}_{km}} - \bar{\boldsymbol{x}}_{km}^{{(j)}^\top} P \bar{\boldsymbol{x}}_{km}^{(j)},
 \end{align}
\end{subequations}
\normalsize
where ${m^*_k}$ is the time epoch that has the worst constraint violation in the drift trajectory emanating from the time step $k$, obtained from the previous iteration of the SCP.
The formulation for the full passive safety of $N_{sc}$ agents with all consideration of different failure timing is available in Ref. \cite{guffanti2022phd}. 

\subsection{Passive Safety via First-Order condition of QPRO}

A new constraint formulation is proposed in this paper to guarantee passive safety throughout the transfer, achieved by applying the first-order condition for staying on a QPRO.
This approach ensures that the deputy always adheres to the surface of a QPRIT of arbitrary size (defined by $\varepsilon$) throughout transfers between QRPOs, thereby achieving passive safety.
The constraint is expressed as
\begin{subequations}
\label{eq:ps_qpt}
\begin{align}
    & \| h_k^{(j)} \| \leq \Delta_h, \quad 
    \| \dot{\alpha}_k^{(j)} \| \leq \Delta_{\dot{\alpha}}, \quad 
    \| \dot{\beta}_k^{(j)} \| \leq \Delta_{\dot{\beta}},  \quad 
    \| \dot{h}_k ^{(j)}\| \leq \Delta_{\dot{h}}, \quad k = 0, ... , N, \forall j \label{eq:scp_qpt_constr1} \\
    & \varepsilon_k^{(j)} = \sqrt{\alpha_k^{{(j)}^2} + \beta_k^{{(j)}^2}} \geq \varepsilon_f^{(j)},  \quad k = 0, ... , N,  \forall j.  \label{eq:scp_qpt_constr2}
\end{align}
\end{subequations}
The deputy's condition to be on a QPRO is $h = \dot{\alpha} = \dot{\beta} = \dot{h} = 0$. 
However, with the strict form of this equality, the problem would be infeasible as LTC becomes stationary (cf. Eq. \eqref{eq:map_LTW_to_cart}). 
Therefore, this constriant is relaxed to $l_2$-norm inequality at each time step as shown in Eq. \eqref{eq:scp_qpt_constr1}, where $\{\Delta_h,\Delta_{\dot{\alpha}},\Delta_{\dot{\beta}},\Delta_{\dot{h}} \}$ are user-defined small parameters.
The choice of these slack parameters depends on the length of the passive safety period and the accuracy of the dynamics modeling. 
Furthermore, Eq. \eqref{eq:scp_qpt_constr2} ensures that the deputy is outside of the QPRIT's surface that corresponds to the terminal QPRO. 
To alleviate nonconvexity, the following convexification using the reference state $\bar{\alpha}_k$ and $\bar{\beta}_k$ is introduced: 
\begin{align}
    \varepsilon_f^{(j)} - \frac{1}{\sqrt{\bar{\alpha}_k^{{(j)}^2} + \bar{\beta}_k^{{(j)}^2} }} \left(\bar{\alpha}_k^{(j)} \alpha_k^{(j)} + \bar{\beta}_k^{(j)} \beta_k^{(j)}\right) \leq 0, \quad k = 0, ..., N, \forall j.
\end{align}

In this proposed approach, safety between agents is achieved solely via passively safe QPRO designs. 
Therefore, the above multi-agent OCP could be solved in a decentralized fashion, once the boundary conditions of each agent are defined. 
If there exists a strict requirement for the minimum separation between agents, a relative orbit design problem must be conducted by iteratively solving the above problem. 
However, the placement of safe relative orbits can be readily done by placing agents in the $\alpha-\beta$ plane of the local toroidal frame, similarly to the spacecraft swarm placement in ROE space \cite{koenig2018robust}, significantly simplifying the QRPO design. 

\section{Results and Analysis}

This section demonstrates the results of trajectory optimizations using the LTC-based nonconvex OCP with the dynamically-informed passive safety constraints.
The primary focus of this section is the application to southern NRHO families around the Earth-Moon $L_2$ point, which are of high practical interest for future cislunar missions, anticipating large-scale swarm deployments and sustained PROD activities over the next decades.

Three primary experiments are performed to examine the performance of the proposed method.
First, passive safety during transfer is analyzed in the CR3BP for a dual-spacecraft system ($N_{sc}=1$) by comparing the two passive safety constraints introduced in Sec. \ref{sec:passive_safety_ocp}.
Secondly, the proposed method is demonstrated in the higher-fidelity dynamical model: ER3BP and BCR4BP. 
Similar dual spacecraft rendezvous scenarios are demonstrated as case studies, and the nonlinear error due to the dynamics linearization is reported.   
Finally, a five-spacecraft swarm reconfiguration strategy ($N_{sc}=4$) is designed in the CR3BP, which is then converted to the full-ephemeris dynamics. 
To address the dynamics modeling error, an MPC scheme is proposed to close the loop. 
Monte Carlo analysis is conducted to validate the proposed methodology. 

Table~\ref{tab:sim_setup} summarizes the scenarios and parameters used for the respective trajectory optimization, such as reference orbits, boundary conditions, and constraint configurations. 
This paper adopts \texttt{SCvx*} \cite{oguri2023successive} as an SCP algorithm to solve a nonconvex optimization. 
Note that because of the usage of a direct method with fixed time discretization, the strict fuel optimality \cite{koenig2017new, hunter2025fast} is not attainable. 
However, as shown in the results, the transfer costs of reconfiguration in this study are on the order of $10^0$–$10^1$ mm/s. 
Therefore, using $N=30$ for time discretization is justified, given that the station-keeping cost of the reference orbit would be significantly higher \cite{davis2022orbit}.
All numerical computation is performed using Intel Core i7-1260P.

\begin{table}[ht!]
\caption{Summary of the Experiment Setup}
\label{tab:sim_setup}
\centering
\scalebox{0.8}{
\begin{tabular}{ll cc c cc c cc}
\hline \hline
\multicolumn{2}{c}{} & \multicolumn{2}{c}{Sec. \ref{sec:safety_comparison_cr3bp}} & & \multicolumn{2}{c}{Sec. \ref{sec:mid_fidelity_demo}} & & \multicolumn{2}{c}{Sec. \ref{sec:hifi_demo}} \\ 
\cline{3-4} \cline{6-7} \cline{9-10}
& Dynamics & \multicolumn{2}{c}{CR3BP} & & ER3BP & BCR4BP & & CR3BP & Full-Ephemeris  \\ 
& \makecell[l]{Orbits\\($L_2$ South NRHO)} & \multicolumn{2}{c}{\makecell{9:2 synodic \\ (Fig. \ref{fig:9_2_nrho})}} & & \makecell{4:1 sidereal\\(Fig. \ref{fig:4_1_er3bp})} & \makecell{3:1 synodic \\ (Fig. \ref{fig:3_1_bcr4bp})} & & \makecell{9:2 synodic \\ (Fig. \ref{fig:9_2_nrho})} & \makecell{ --- \\ (Fig. \ref{fig:9_2_nrho})}  \\ 
& Initial fixed point & \multicolumn{2}{c}{apolune} & & \multicolumn{2}{c}{apolune} && \multicolumn{2}{c}{apolune}\\ 
& $t_f$ & \multicolumn{2}{c}{2 revs.} & & \multicolumn{2}{c}{2 revs.} && \multicolumn{2}{c}{2 revs.} \\ 
& $N$ & \multicolumn{2}{c}{30} && \multicolumn{2}{c}{30} && \multicolumn{2}{c}{30} \\ 
& $[k_p^{-}, k_p^{+}]$ & \multicolumn{2}{c}{[6,9], [21,24]} && \multicolumn{2}{c}{[6,9], [21,24]} && \multicolumn{2}{c}{[6,9], [21,24]} \\ 
& $N_{sc}$ & \multicolumn{2}{c}{1} & & 1 & 1 & & \multicolumn{2}{c}{4} \\ 
& $(\varepsilon_i, \theta_i)$ [km, rad] & \multicolumn{2}{c}{(0.5, 4.2)} & &  (0.5, 4.2) & (1, 3.8) & & \multicolumn{2}{c}{(0.5, 4.2), (0.75, 4.2-$\pi$), (2.5, 4.2), (1.5, 4.2-$\pi$)}  \\ 
& $(\varepsilon_f, \theta_f)$ [km, rad]  & \multicolumn{2}{c}{(0.2, 0)} & & (0.2, 0) & (0.7, 0) & & \multicolumn{2}{c}{(0.2, 0), (0.3, $\pi$), (1, 0), (0.6, $\pi$)}  \\ 
& Safety & SCP-DRIFT & SCP-QRPIT & & \multicolumn{2}{c}{SCP-QRPIT} & & \multicolumn{2}{c}{SCP-QRPIT}\\
\hline
\multirow{2}{*}{QPRIT} & $\Delta_h$ & {---} & 1 m & & \multicolumn{2}{c}{1 m} & & \multicolumn{2}{c}{1 m}\\ 
& $\Delta_{\dot{\alpha}}, \Delta_{\dot{\beta}}, \Delta_{\dot{h}}$ & --- & 50 mm/s & &  \multicolumn{2}{c}{50 mm/s} & &  \multicolumn{2}{c}{50 mm/s}\\ 
\hline
\multirow{2}{*}{DRIFT} & $T_{\text{safe}}$ & 1 rev. & --- & & --- & --- &\multirow{2}{*}{MPC} & $Q$ & \makecell{position: 5e-8 m/s \\ velocity: 5e-10 m/s$^2$} \\
& $N_{\text{safe}}$ & 30 & --- & & --- & --- & & $k_{\text{MPC}}$ & 13, 25 \\ 
\hline \hline
\end{tabular}
}
\end{table}

In this case study, if there are multiple oscillatory modes in a reference orbit, one with a smaller maximum along-track separation of the QPRIT is chosen to construct LTC. 
Such a torus structure is more desirable because a highly eccentric invariant circle requires a larger torus to ensure the minimum separation around the chief.

\subsection{Comparison of Passively Safe Reconfiguration Strategies in the CR3BP} \label{sec:safety_comparison_cr3bp}

\begin{figure}[ht!]
    \centering
     \begin{subfigure}[b]{0.24\textwidth}
         \centering
         \includegraphics[width=\textwidth]{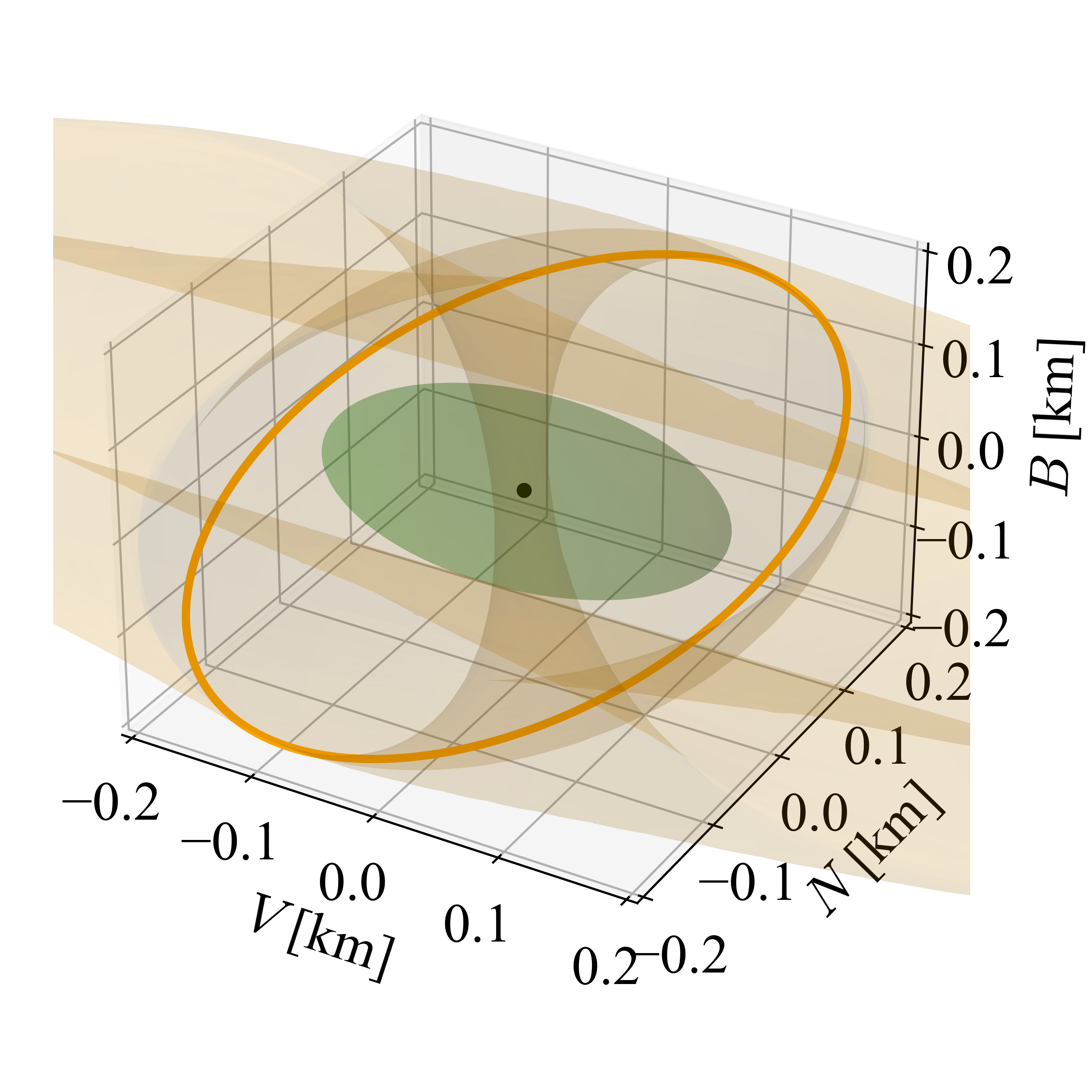}
         \caption{3D view}
         \label{fig:koz_config}
     \end{subfigure}
     \begin{subfigure}[b]{0.24\textwidth}
         \centering
         \includegraphics[width=\textwidth]{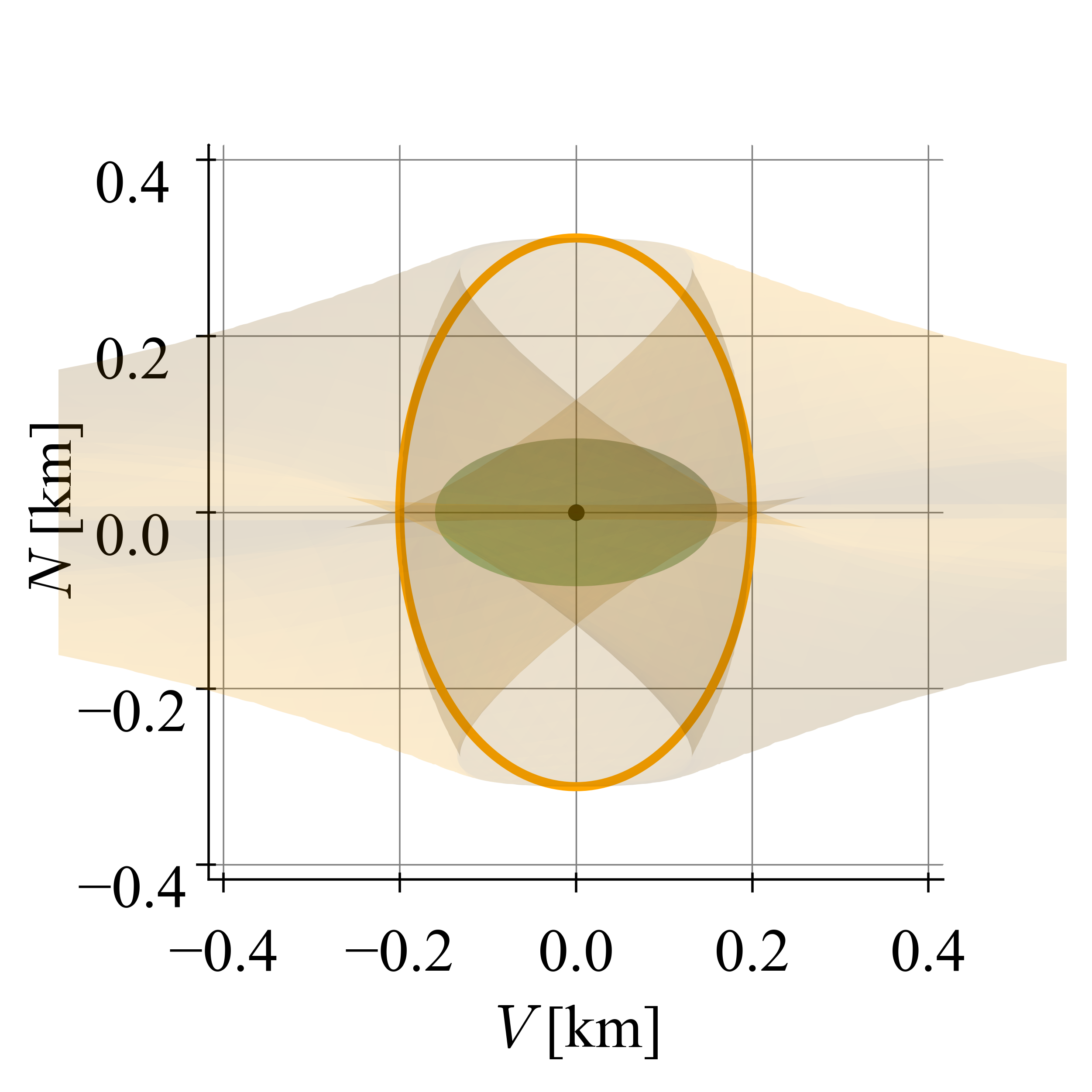}
         \caption{VN plane projection}
         \label{fig:koz_tn}
     \end{subfigure}  
     \begin{subfigure}[b]{0.24\textwidth}
         \centering
         \includegraphics[width=\textwidth]{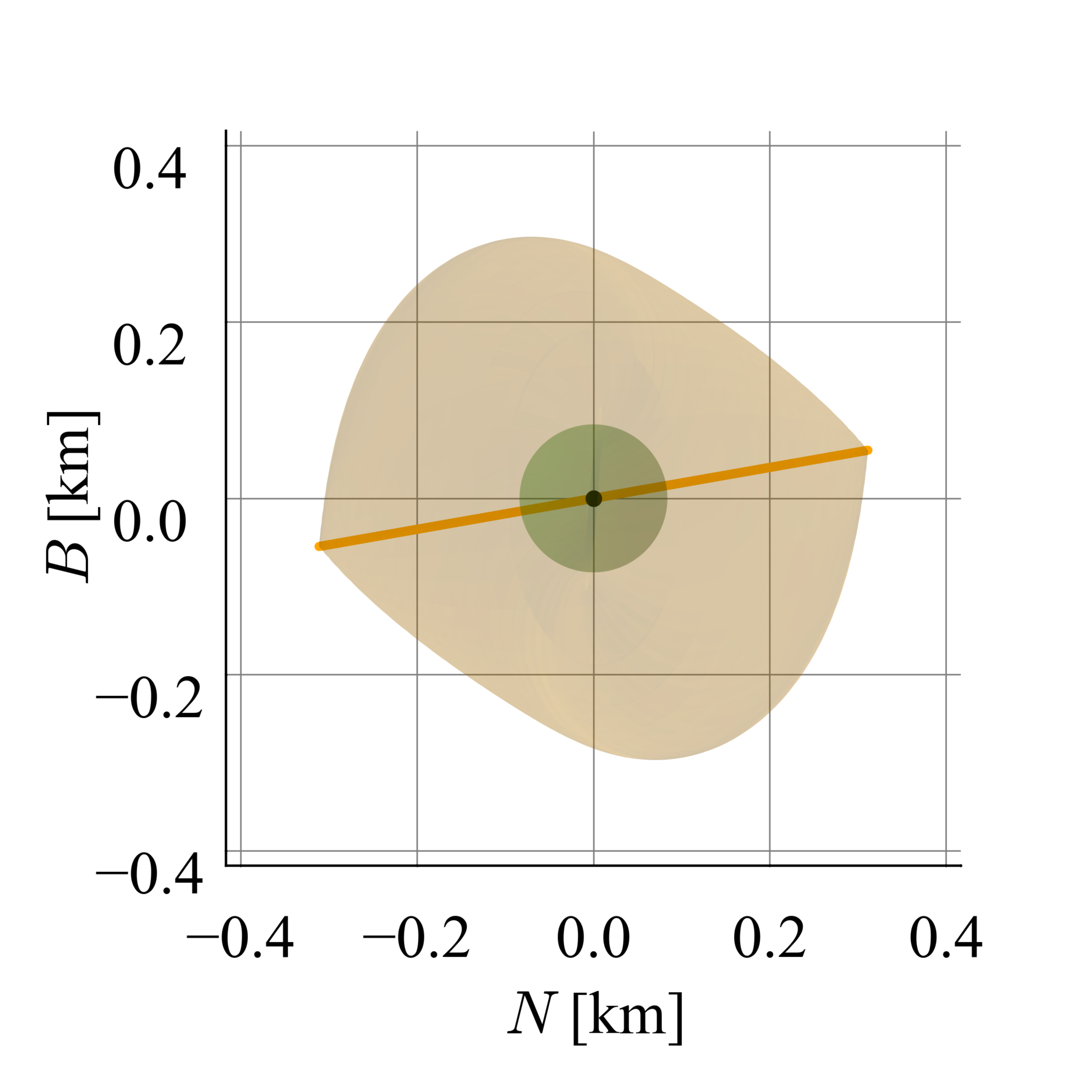}
         \caption{NB plane projection}
         \label{fig:koz_nw}
     \end{subfigure}  
     \begin{subfigure}[b]{0.24\textwidth}
         \centering
         \includegraphics[width=\textwidth]{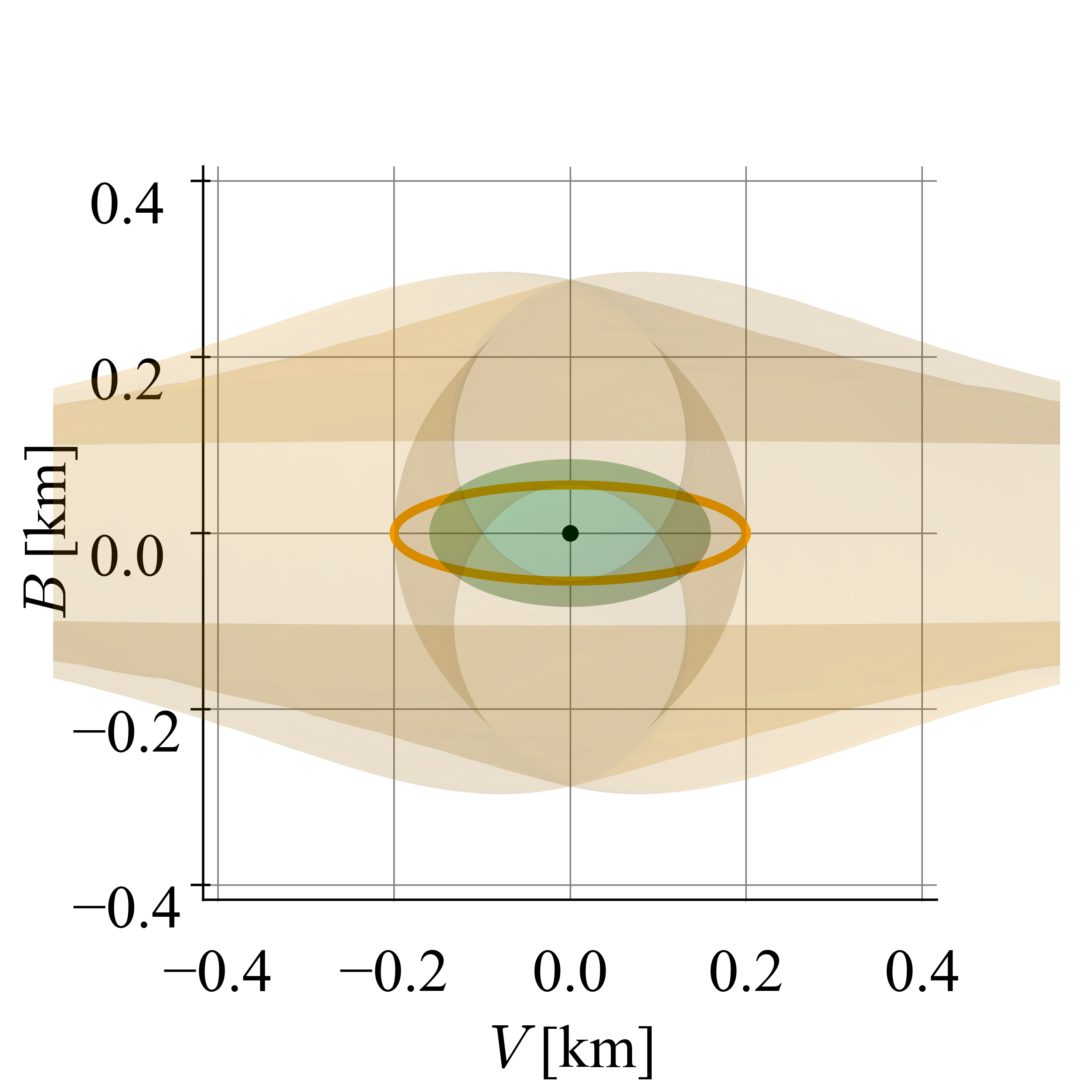}
         \caption{BV plane projection}
         \label{fig:koz_tw}
     \end{subfigure}  
    \caption{ Keep-out zone (green) in the VNB frame and the passively-safe terminal QPRIT (orange).}
    \label{fig:koz} 
\end{figure}

Two passive safety strategies are examined through a dual spacecraft proximity operation mission scenario in the CR3BP $L_2$ South 9:2 synodic resonance NRHO. 
Since it takes approximately 8 orbits for the relative orbits to return to a similar torus state in $\mathbb{T}^2$ in this orbit (cf. Fig. \ref{fig:inv_circ_amplitude}), a relative transfer is considered between particular phases of QPROs, for which a solution of a two-point boundary value problem lacks passive safety. 
The geometry of the KOZ is defined to have a stretch in the velocity direction as $(x_{e}, y_{e}, z_{e}) = (200, 95, 95)$m, from which both the initial and terminal QPROs are passively safe. 
The QPRIT and KOZ are illustrated in Fig. \ref{fig:koz}. 
The orange ellipse corresponds to the osculating invariant circle at the apolune. 
\begin{figure}[ht!]
     \centering
     \begin{subfigure}[b]{0.33\textwidth}
         \centering
         \includegraphics[width=\textwidth]{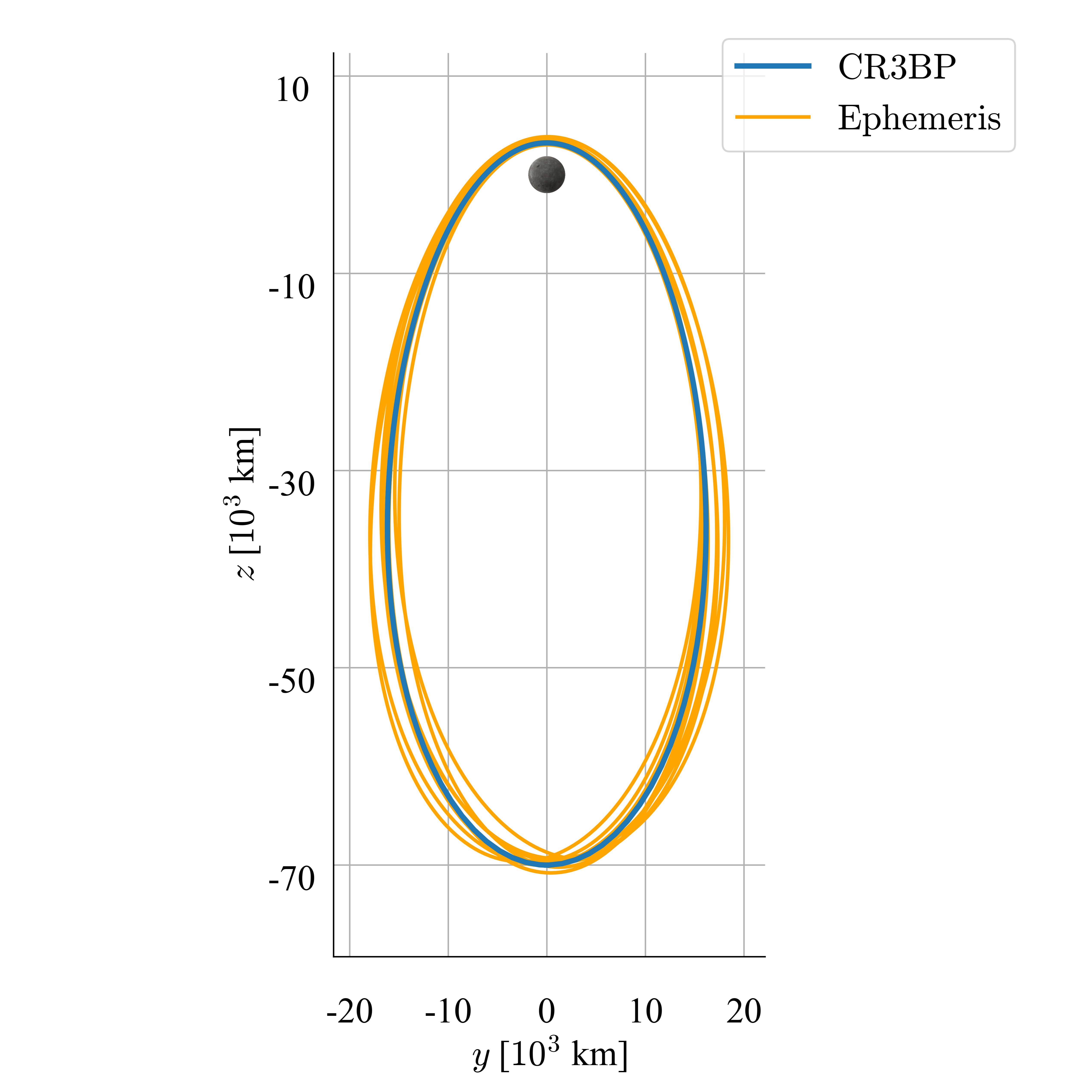}
         \caption{9:2 Synodic}
         \label{fig:9_2_nrho}
     \end{subfigure}
     \begin{subfigure}[b]{0.33\textwidth}
         \centering
         \includegraphics[width=\textwidth]{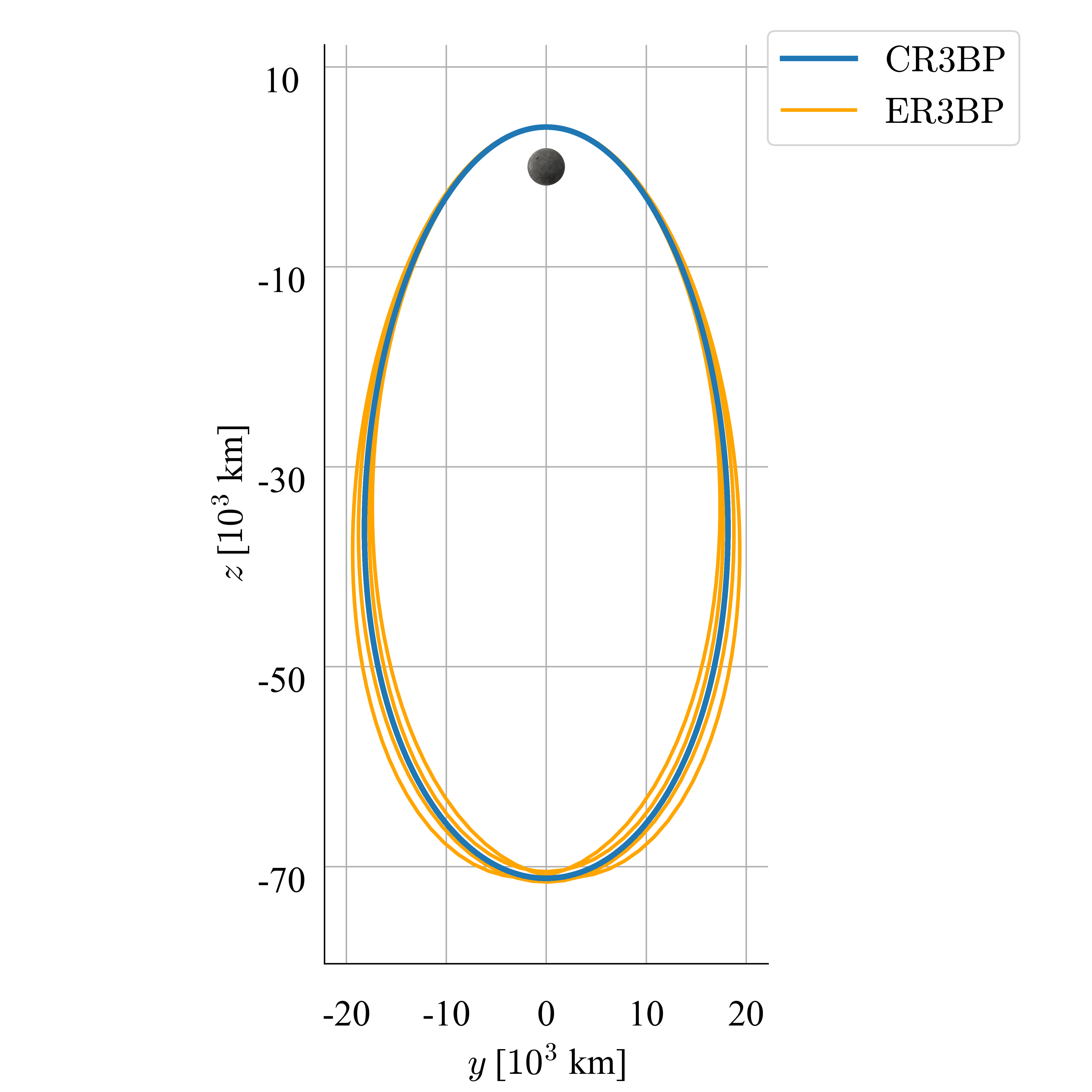}
         \caption{4:1 Sidereal}
         \label{fig:4_1_er3bp}
     \end{subfigure}  
     \begin{subfigure}[b]{0.33\textwidth}
         \centering
         \includegraphics[width=\textwidth]{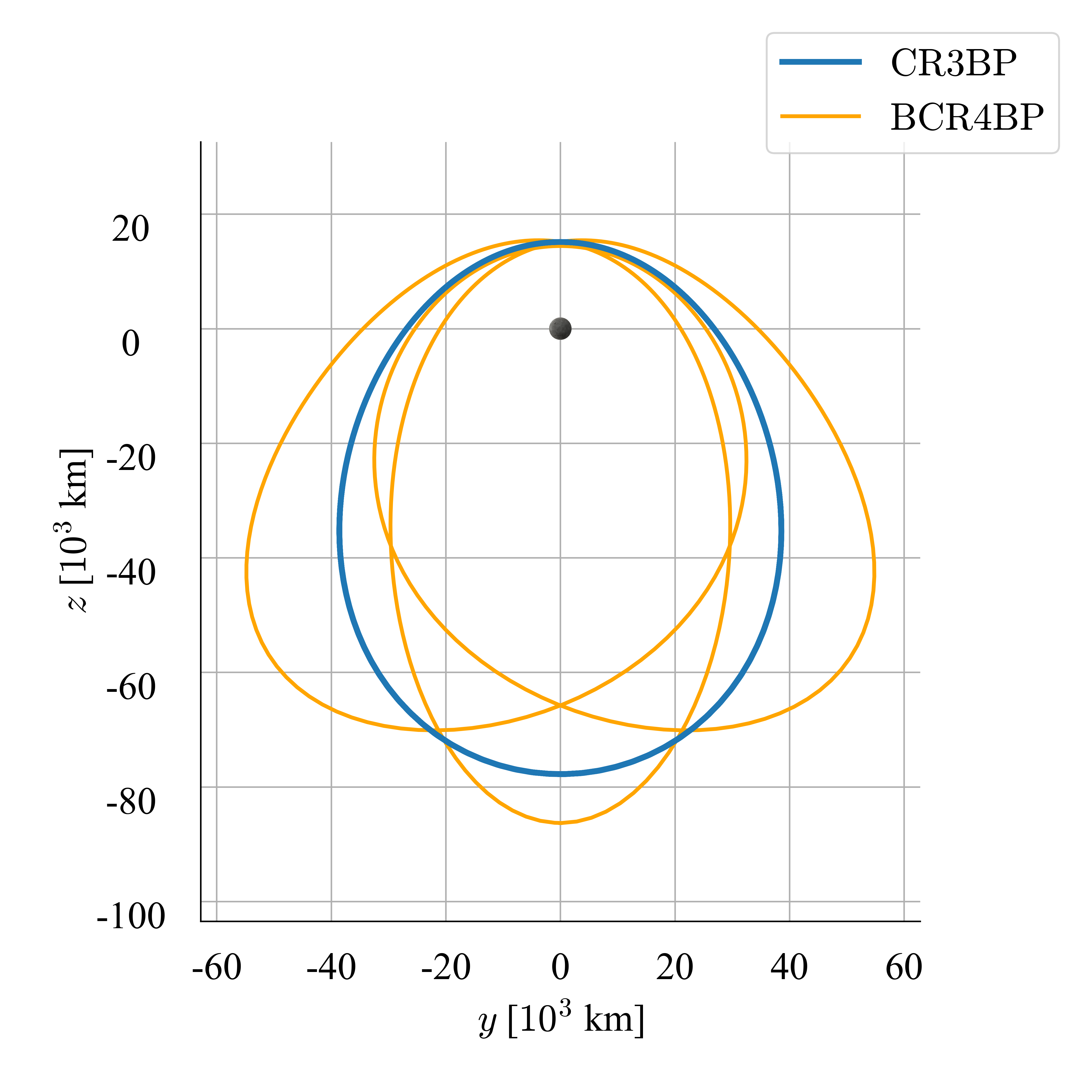}
         \caption{3:1 Synodic}
         \label{fig:3_1_bcr4bp}
     \end{subfigure}  
    \caption{Earth-Moon $L_2$ South NRHOs used for case studies (Moon synodic frame)}
    \label{fig:nrho_ref} 
\end{figure}

The solution to the nonconvex problem with the drift-safe constraint in Eq. \eqref{eq:drifty_safety} and the first-order QPRIT condition in Eq.~\eqref{eq:ps_qpt} is referred to as SCP-DRIFT and SCP-QPRIT, respectively.
The solution of a convex OCP that only considers the boundary conditions and dynamics constraints, denoted as CVX, serves as the initial guess for the SCP-DRIFT, whereas the convex solution of a convex OCP with the $l_2$-norm constraints in Eq. \eqref{eq:scp_qpt_constr1}, denoted as CVX-QPRIT, serves as the initial guess for the SCP-QPRIT.  

\begin{table}[ht]
\centering
\caption{Experiment results of open-loop trajectories in CR3BP (Sec. \ref{sec:safety_comparison_cr3bp}) }
\label{tab:exp_results}
\scalebox{0.8}{
    \begin{tabular}{c cccc}
    \hline\hline
    \multicolumn{1}{c}{Safety strategy} 
    & CVX & SCP-DRIFT & CVX-QPRIT & SCP-QPRIT \\
    \hline
    STM computation [s] & 15.63 & 100.93 & 15.63 & 15.63 \\
    Opt. runtime [s]  & 0.101 & 1.677 & 0.279 & 4.946  \\
    One-rev. passive safety &  & \checkmark &  & \checkmark  \\
    Fuel cost [mm/s]  & 5.157 & 5.220 & 7.287 & 8.510 \\
    \makecell{Nonlinear error \\ (terminal position) [m]} & 0.008 &  0.0167 & 0.280 & 0.501 \\
    \hline\hline
    \end{tabular}
    }
\end{table}

\begin{figure}[ht!]
     \centering
     \begin{subfigure}[b]{0.35\textwidth}
         \centering
         \includegraphics[width=\textwidth]{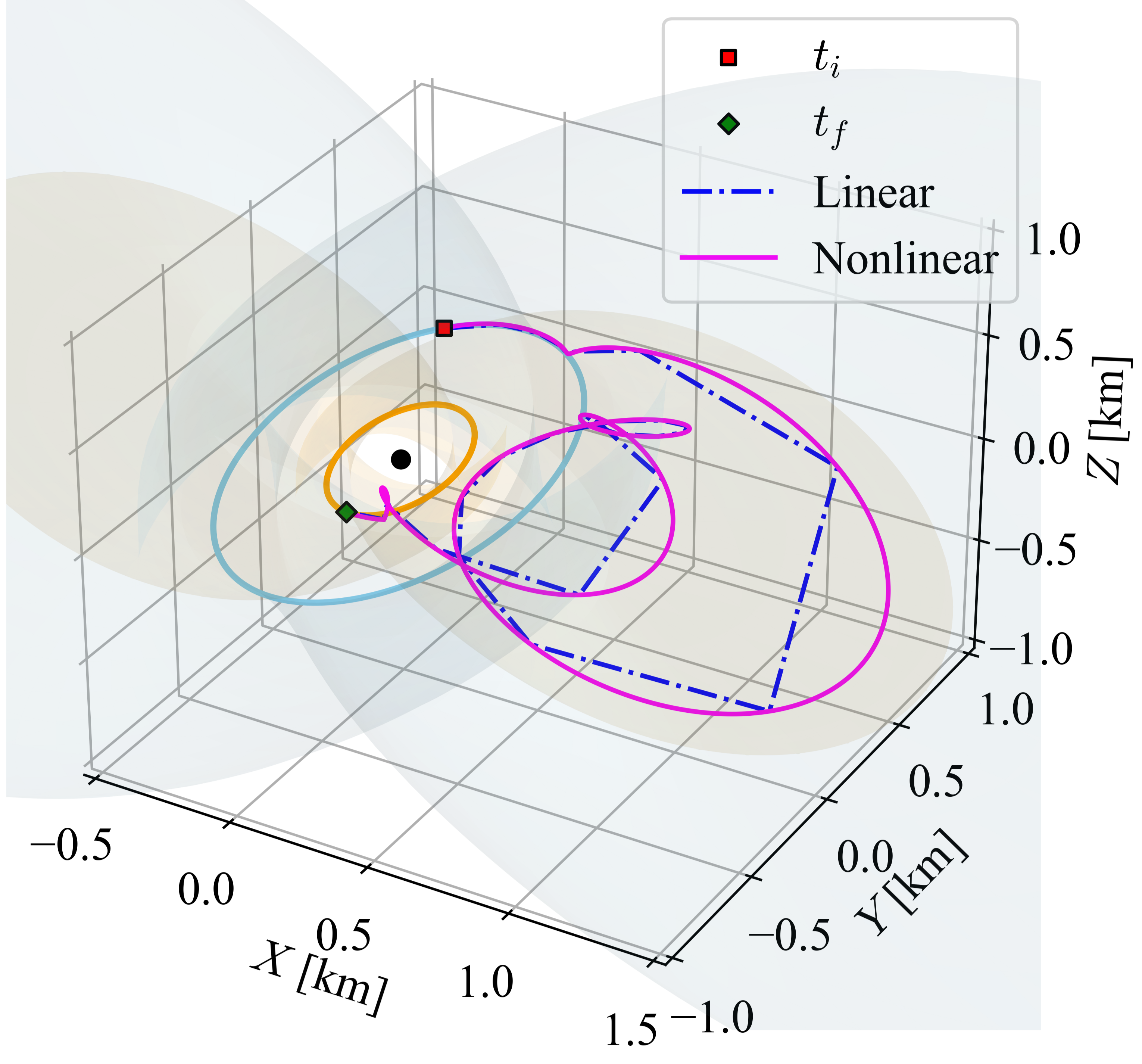}
         \caption{LVLH frame}
         \label{fig:OptTraj_lvlh_scp}
     \end{subfigure}
     \begin{subfigure}[b]{0.35\textwidth}
         \centering
         \includegraphics[width=\textwidth]{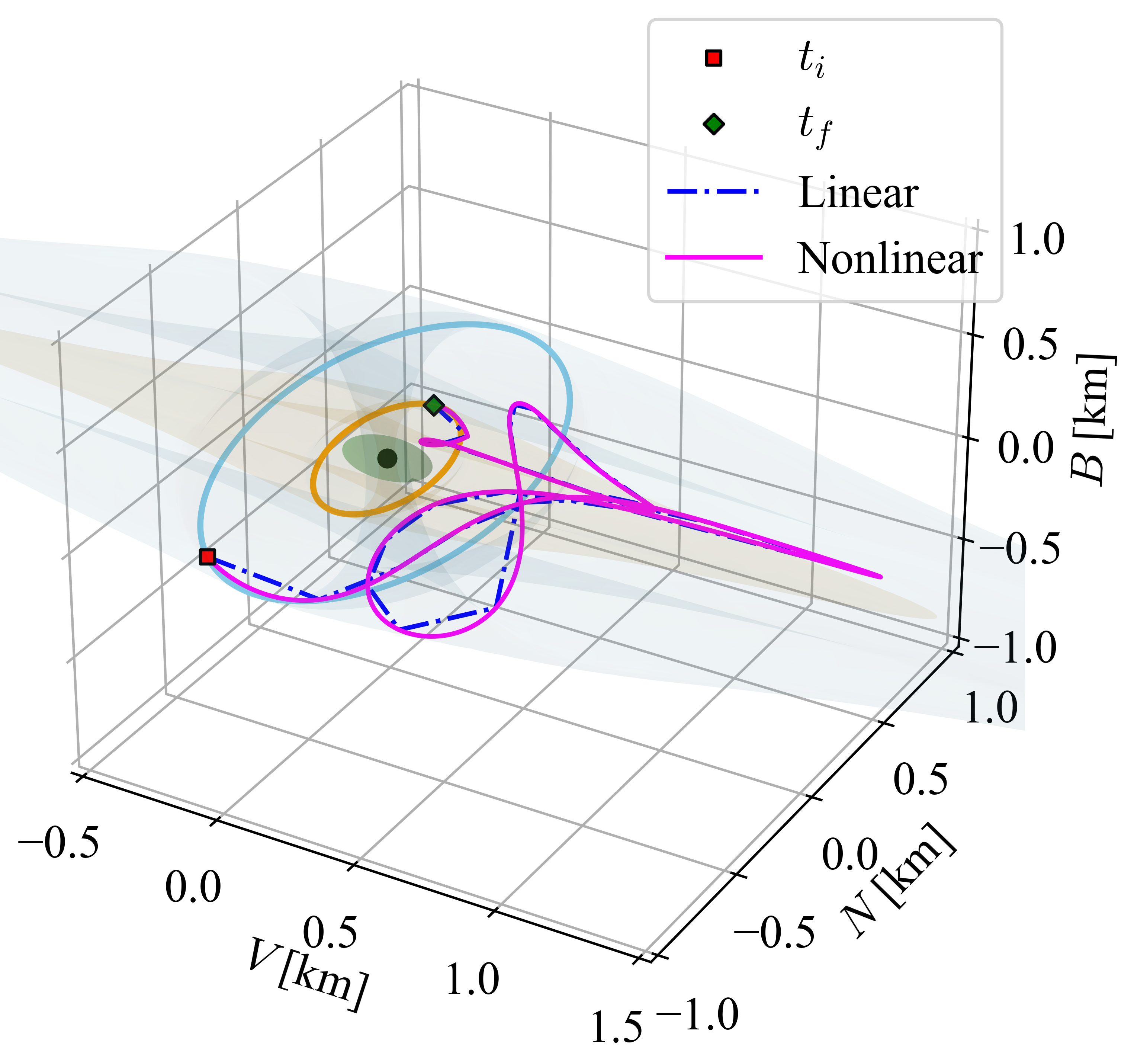}
         \caption{VNB frame}
         \label{fig:OptTraj_tnw_scp}
     \end{subfigure}  
     \begin{subfigure}[b]{0.28\textwidth}
         \centering
         \includegraphics[width=\textwidth]{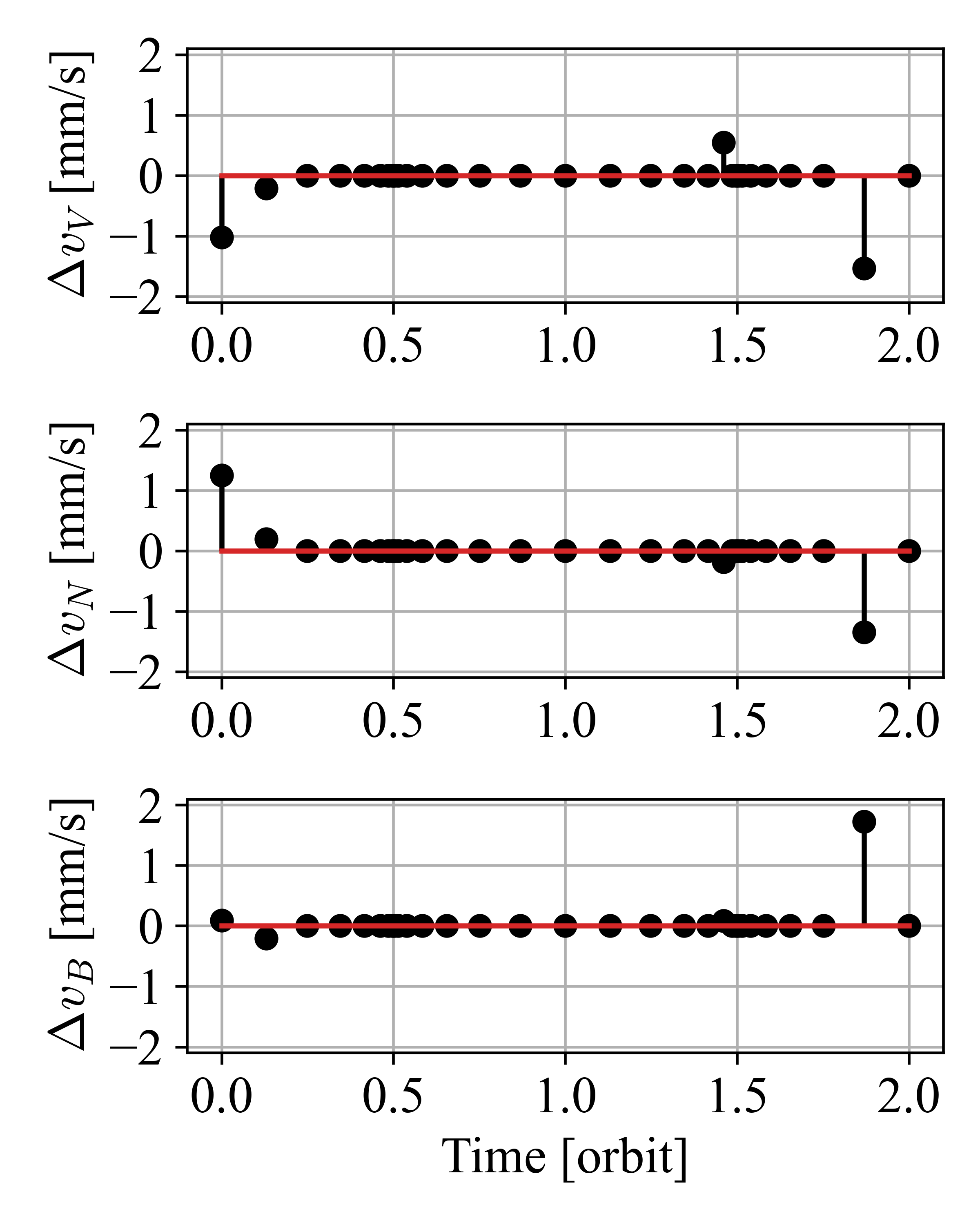}
         \caption{Control history (VNB frame)}
         \label{fig:cntrl_hist_DRIFT}
     \end{subfigure}  
    \caption{Optimized trajectory of the SCP-DRIFT and its nonlinear propagation in the co-moving frames. }
    \label{fig:opt_traj_scp} 
\end{figure}
\begin{figure}[ht!]
     \centering
     \begin{subfigure}[b]{0.35\textwidth}
         \centering
         \includegraphics[width=\textwidth]{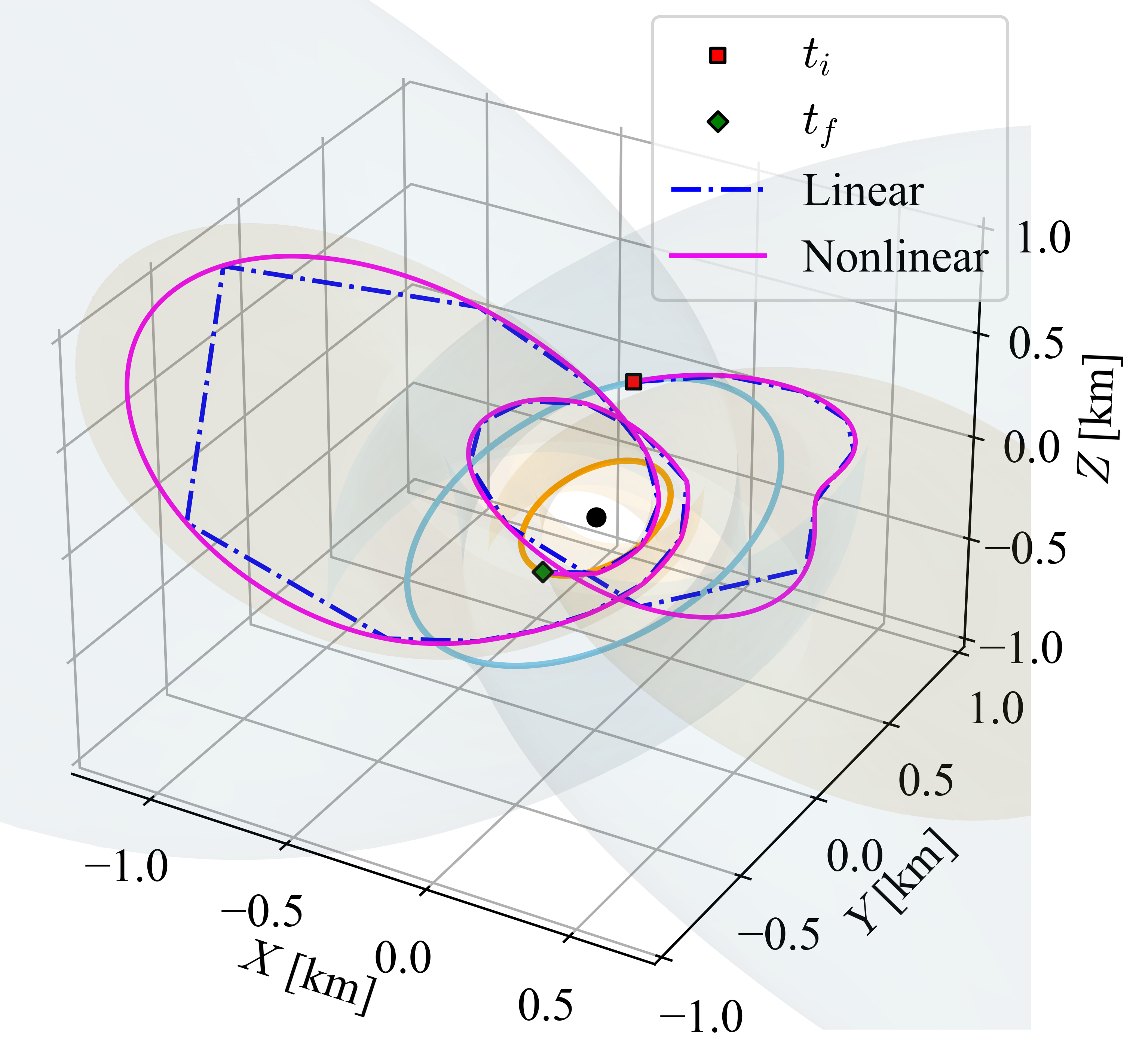}
         \caption{LVLH frame}
         \label{fig:OptTraj_lvlh_qpt}
     \end{subfigure}
     \begin{subfigure}[b]{0.35\textwidth}
         \centering
         \includegraphics[width=\textwidth]{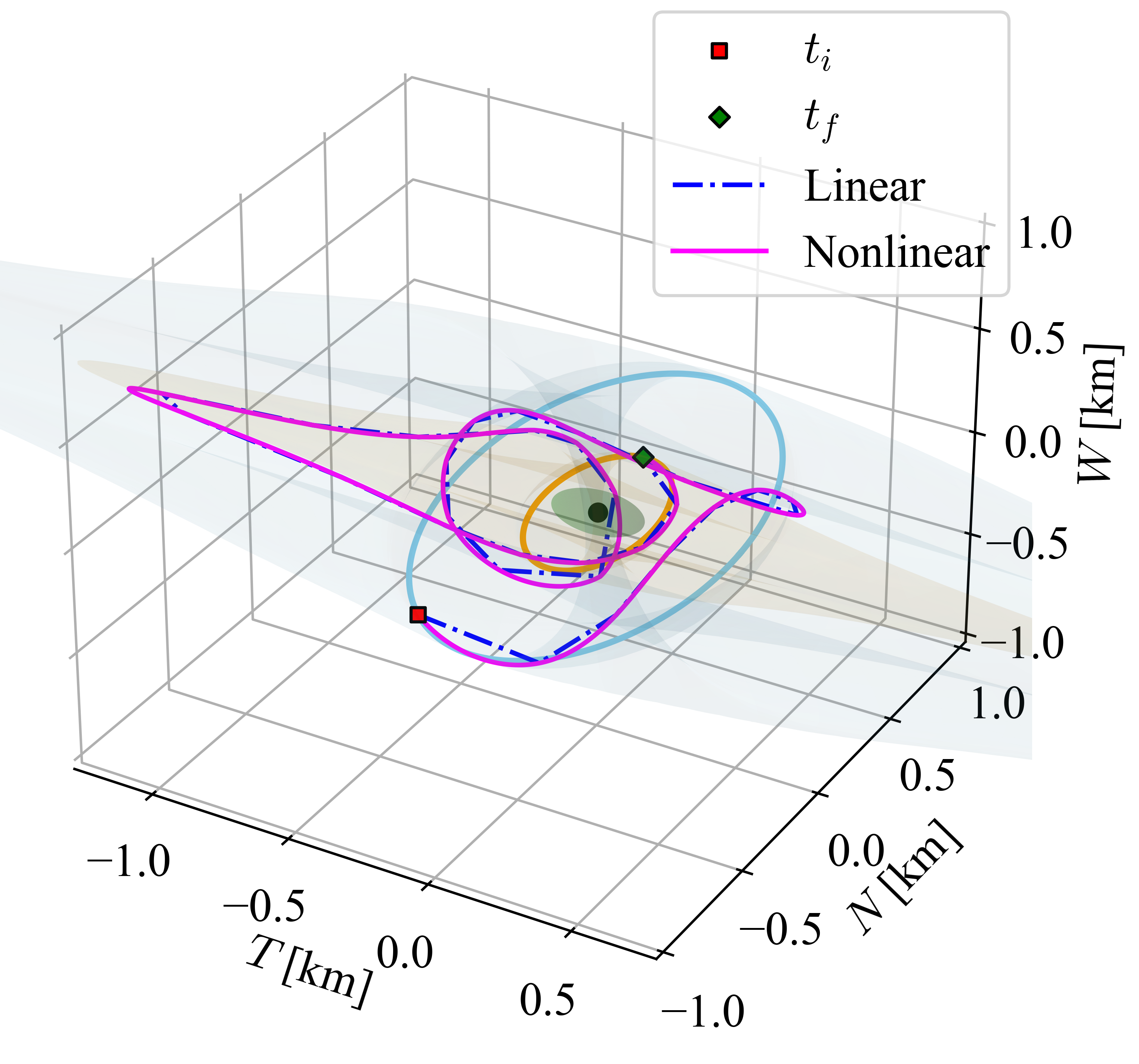}
         \caption{VNB} frame
         \label{fig:OptTraj_tnw_qpt}
     \end{subfigure}  
     \begin{subfigure}[b]{0.28\textwidth}
         \centering
         \includegraphics[width=\textwidth]{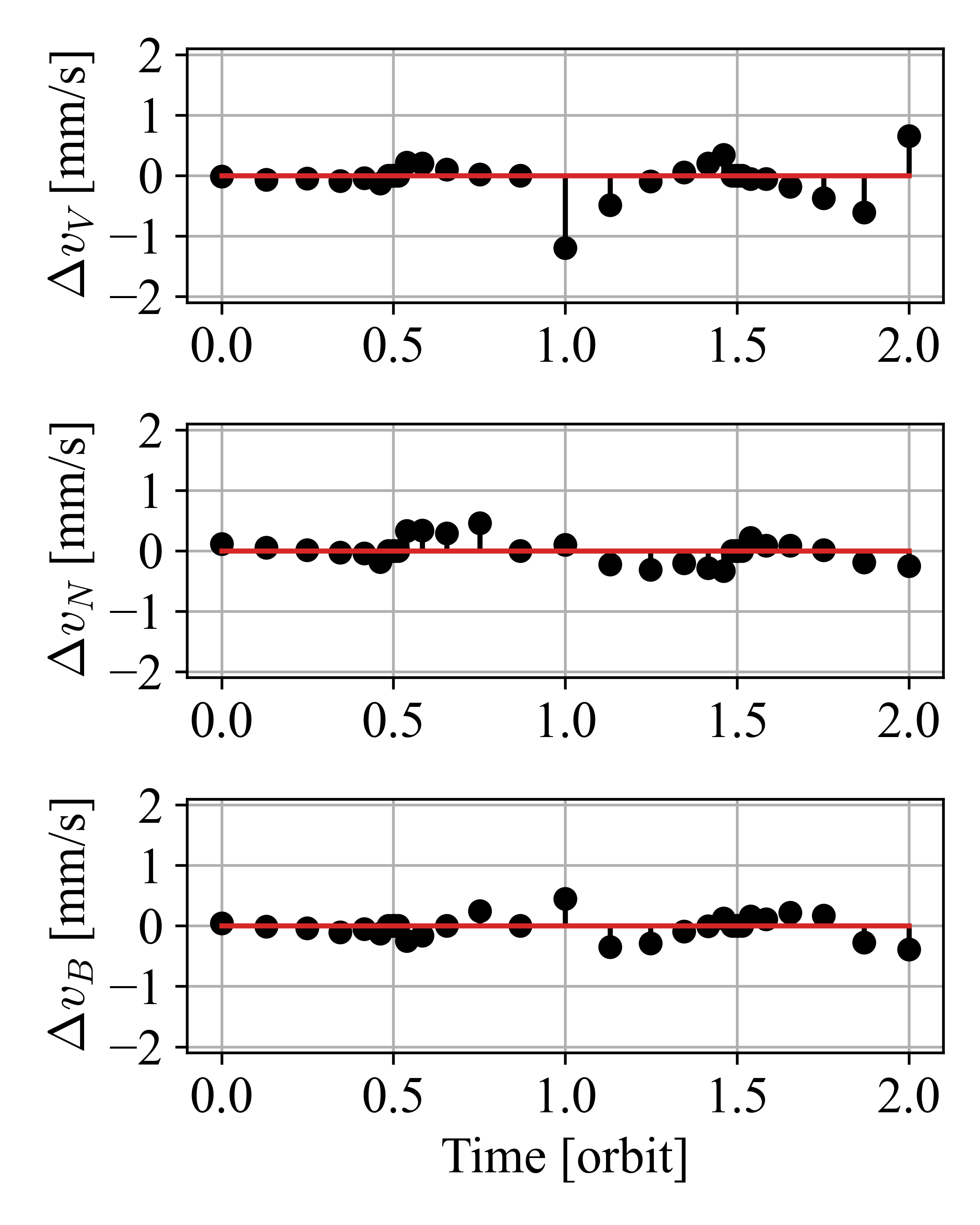}
         \caption{Control History (VNB} frame)
         \label{fig:cntrl_hist_QPRIT}
     \end{subfigure}  
    \caption{Optimized trajectory of the SCP-QPRIT and its nonlinear propagation in the co-moving frames.}
    \label{fig:opt_traj_qpt} 
\end{figure}
\begin{figure}[ht!]
     \centering
     \begin{subfigure}[b]{0.47\textwidth}
         \centering
         \includegraphics[width=\textwidth]{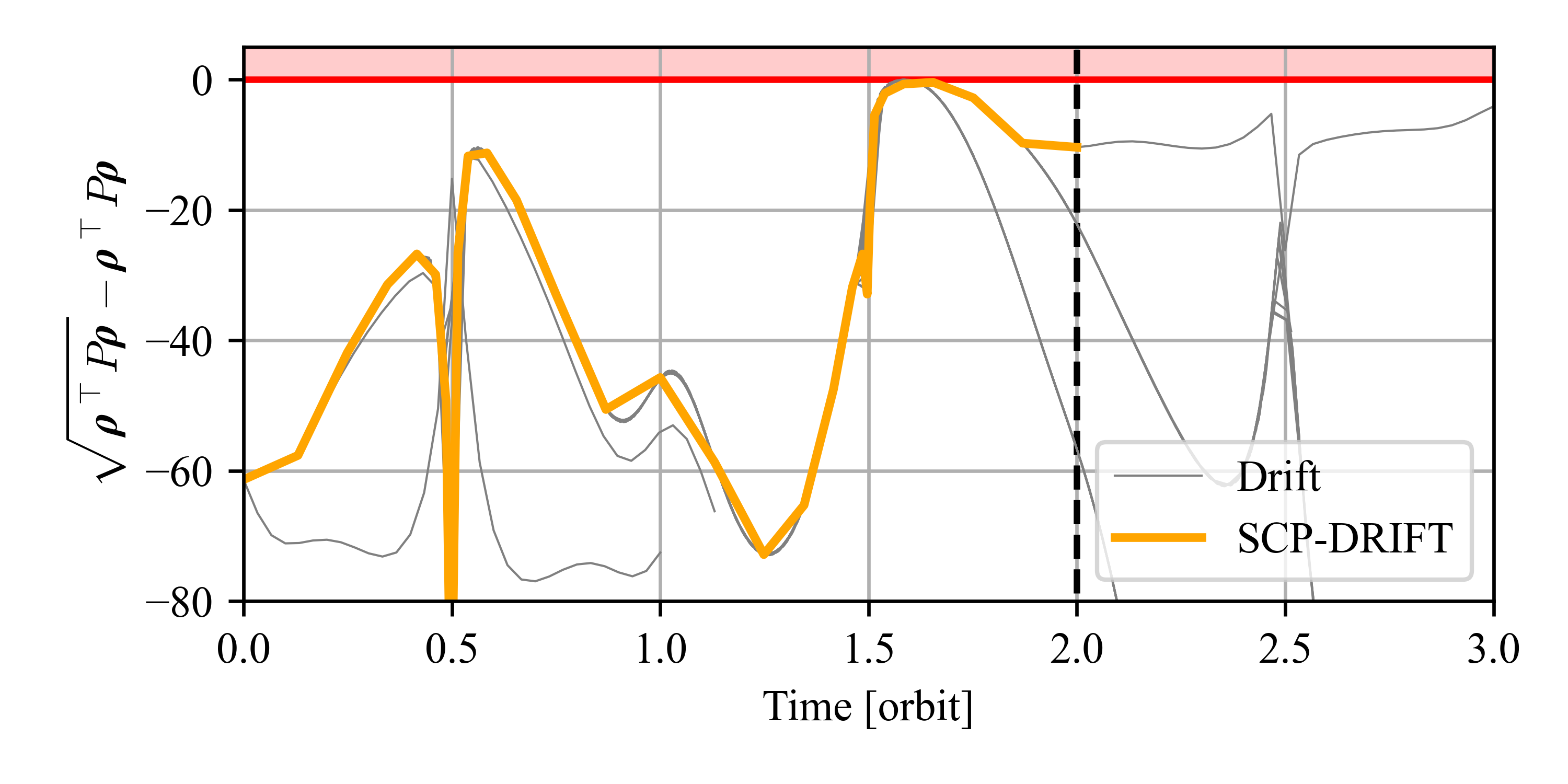}
         \caption{SCP-DRIFT}
         \label{fig:ps_constr_drift}
     \end{subfigure}
     \begin{subfigure}[b]{0.47\textwidth}
         \centering
         \includegraphics[width=\textwidth]{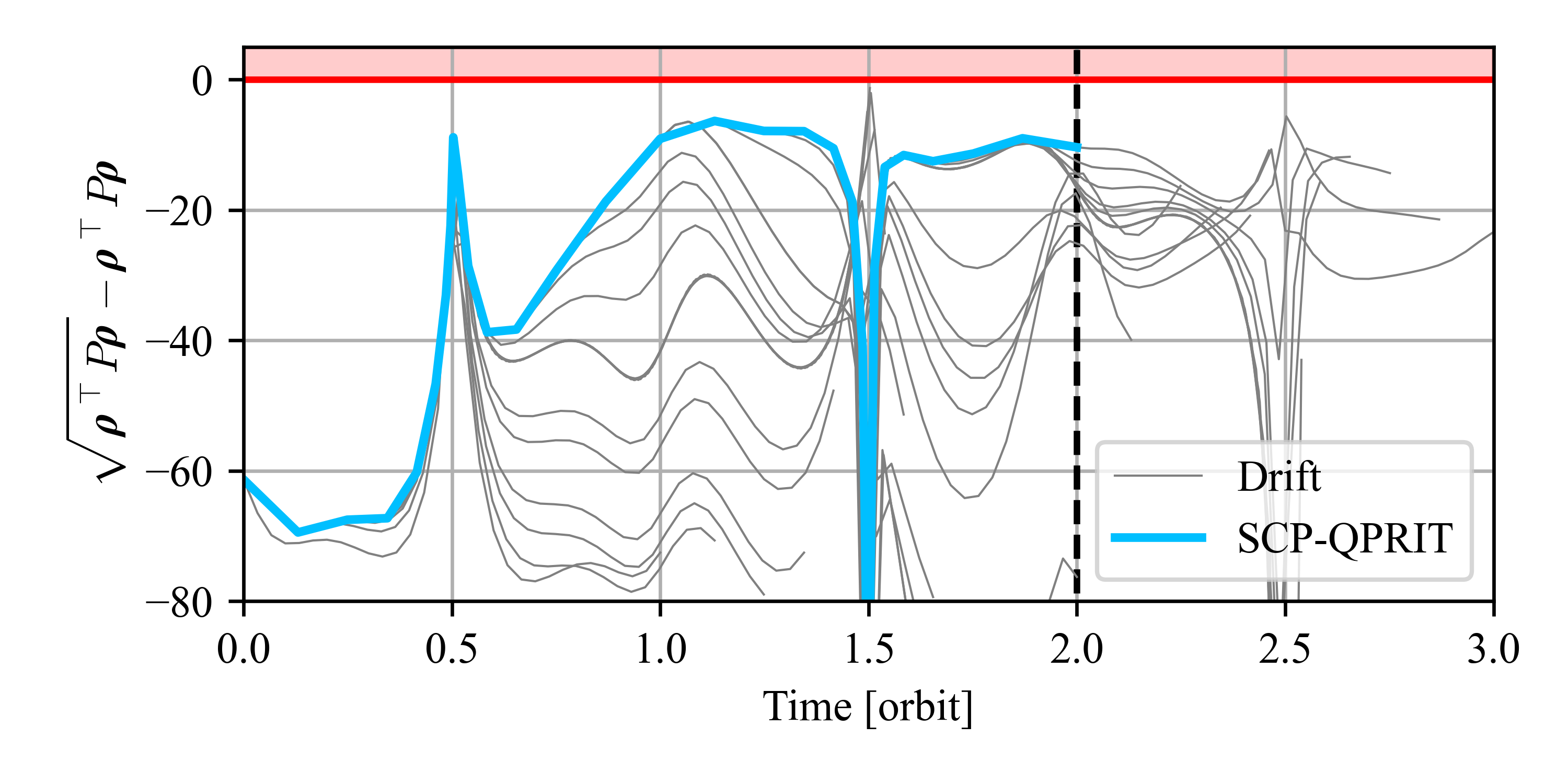}
         \caption{SCP-QPRIT}
         \label{fig:ps_constr_qpt}
     \end{subfigure}  
    \caption{Satisfaction of the passive safety constraint for drift trajectories.}
    \label{fig:ps_constr_hist} 
\end{figure}

The results of the optimization are summarized in Table \ref{tab:exp_results}, which enumerates (i) upfront computation time of monodromy matrix, history of eigenvectors, and STMs, (ii) optimization runtime, (iii) one-revolution passive safety, (iv) fuel cost, and (v) nonlinear error of the optimized trajectory evaluated at its terminal position. 
Note that the computation of STMs and other parameters is shared between CVX, CVX-QPRIT, and SCP-QPRIT, whereas additional computation is required for the SCP-DRIFT to propagate STMs for drift trajectories. 
Both SCPs converged to the local optimum, and the optimized trajectories are illustrated in Figures \ref{fig:opt_traj_scp} and \ref{fig:opt_traj_qpt}, respectively.
Trajectories are resolved in both the VNB frame and the LVLH frame, and the impulsive control profile is also shown in each figure. 
The sky blue and orange surfaces are the initial and terminal QPRIT, respectively. 
Bold ellipses attached to the QPRITs illustrate the invariant circle at the initial and terminal epoch. 
Furthermore, the control profile is propagated along the nonlinear relative dynamics of the CR3BP, and the resulting nonlinear trajectories (solid magenta) are compared to the optimized linear trajectories (dotted blue).
Additionally, the passive safety constraint in Eq. \eqref{eq:drifty_safety1} is evaluated for each optimized trajectory and summarized in Fig. \ref{fig:ps_constr_hist}. 
The state below the red line indicates the satisfaction of the safety constraint. 
It is observed from the figures that the one-orbit passive safety is achieved in both optimal trajectories. 
The SCP-DRIFT converges to a total cost of $\Delta v_{tot} = 5.220$ mm/s, whereas the SCP-QPRIT yields $\Delta v_{tot} = 8.510$ mm/s.
The transfer cost is higher in the SCP-QPRIT due to the conservative nature of the passive safety constraint. 
It is evident from Fig. \ref{fig:OptTraj_tnw_scp} that the controlled trajectory is inside the surface of the QPRIT that corresponds to the terminal state and makes a tangent contact with the surface of the KOZ. 
In contrast, the converged solution of SCP-QPRIT ensures that the deputy not only remains outside of the terminal QPRIT but also leverages the dynamical structure by staying on the surface of QPRIT throughout the transfer. 
Thus, the worst (maximum) constraint value shown in Fig. \ref{fig:ps_constr_qpt} is below the constraint violation threshold and achieves one-revolution passive safety without explicitly constraining the future state. 
A comparison of Figures \ref{fig:cntrl_hist_DRIFT} and \ref{fig:cntrl_hist_QPRIT} highlights the distinct optimal control strategies. 
While the SCP-DRIFT yields a control history resembling a two-impulse solution, the SCP-QPRIT solution spreads small impulses across the transfer, ensuring adherence to the QPRIT surface at each time step and the smooth transition along the QPRIT surface from the initial to the terminal state.
Identifying which thrust profiles are more prone to nonlinear errors is not straightforward. 
SCP-QPRIT applies thrusts near the perilune immediately before and after the no-control windows, where system dynamics are highly sensitive, while SCP-DRIFT uses a small number of larger impulses that can potentially violate linearization assumptions.
However, the nonlinear errors at the terminal positions are overall negligible in the case study as shown in Table \ref{tab:exp_results}, as the total $\Delta v$ values are typically small. 
Given that the evaluation is purely open-loop, the results show that the linearized dynamics provide a sufficiently accurate approximation of the nonlinear dynamics for the problems considered.

\begin{figure}[ht!]
     \centering
     \begin{subfigure}[b]{0.45\textwidth}
         \centering
         \includegraphics[width=\textwidth]{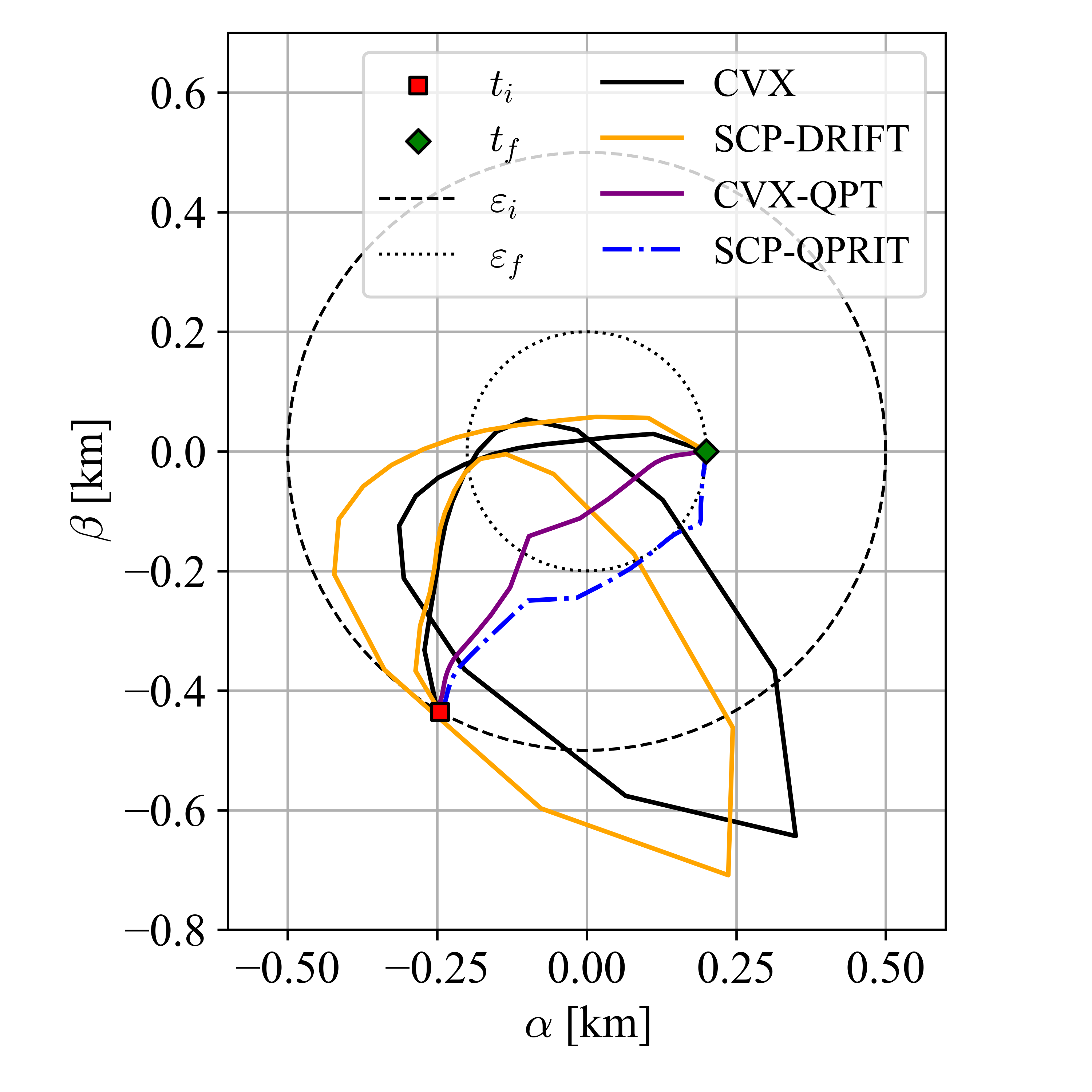}
         \caption{$\alpha-\beta$ plane}
         \label{fig:traj_ltc1}
     \end{subfigure}
     \begin{subfigure}[b]{0.45\textwidth}
         \centering
         \includegraphics[width=\textwidth]{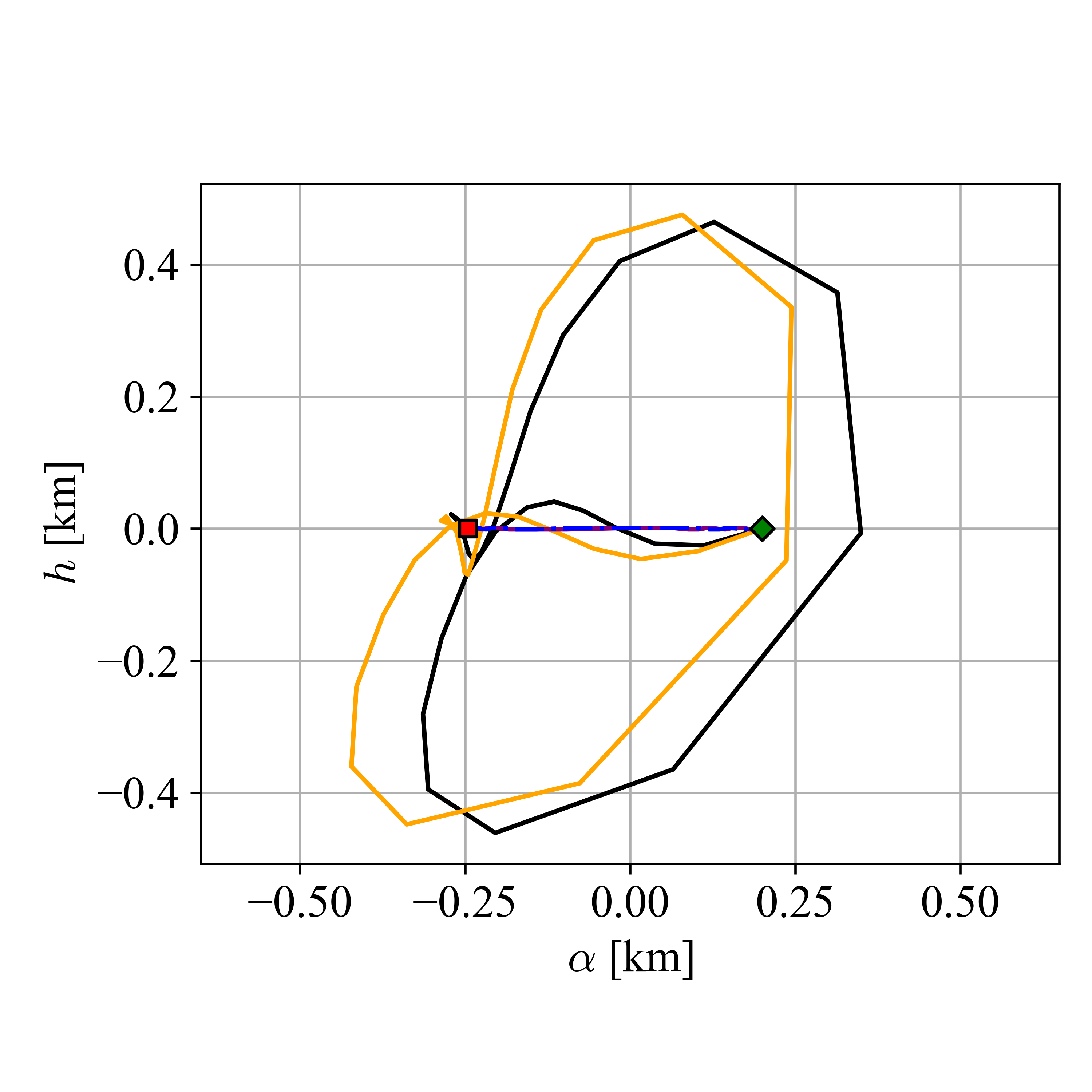}
         \caption{$\alpha-h$ plane}
         \label{fig:traj_ltc2}
     \end{subfigure}  
    \caption{Optimal trajectories in the nonsingular local toroidal coordinates.}
    \label{fig:traj_ltc} 
\end{figure}

Fig. \ref{fig:traj_ltc} illustrates the optimal trajectories projected on the LTC. 
In addition to the two optimal solutions, the solutions of CVX and CVX-QPRIT are also presented.
The solutions of the CVX and SCP-DRIFT exhibit significant variations in the $h$ component, which is orthogonal to the local eigenspace.
This indicates a separation from the QPRIT surface and a loss of quasi-periodicity during the transfer. 
Furthermore, although these solutions include a segment with a radius ($\varepsilon$) smaller than the terminal radius ($\varepsilon_f$), the SCP-DRIFT solution retains one-orbit passive safety.
This is because a passive safety constraint given by Eq. \eqref{eq:drifty_safety1} does not evaluate infinite-time passive safety under natural dynamics.
Conversely, the SCP-QPRIT solution remains constrained within the eigenspace and avoids entering the terminal radius, ensuring passive safety through the evaluation of the LTC.
This is an advantage of using the LTC in the OCP, where the safety condition is geometrically evident and operationally intuitive. 
Additionally, as shown in Table \ref{tab:exp_results}, the computational overhead of propagating an STM and a drifted state at each iteration of SCP leads to a significantly larger computational overhead in SCP-DRIFT, where the linearization of the nonconvex constraint in SCP-QPRIT significantly reduces this burden. 
In the following sections, only the SCP-QPRIT solutions are presented. 

\subsection{LTC-based safe and optimal Reconfiguration in ER3BP and BCR4BP} \label{sec:mid_fidelity_demo}

The original work on LTC \cite{elliott2022phd} only considered the CR3BP, while this can also be extended to mid-fidelity dynamics models to increase the fidelity of the trajectories that may be more suitable for a target reference orbit. 
This subsection demonstrates the proposed relative guidance and control strategy in ER3BP and BCR4BP to orbits where the perturbation from the CR3BP solution is particularly salient.  

\subsubsection{ER3BP $L_2$ South 4:1 Halo Orbit (sidereal resonance)}

\begin{figure}[ht!]
     \centering
     \begin{subfigure}[b]{0.5\textwidth}
         \centering
         \includegraphics[width=\textwidth]{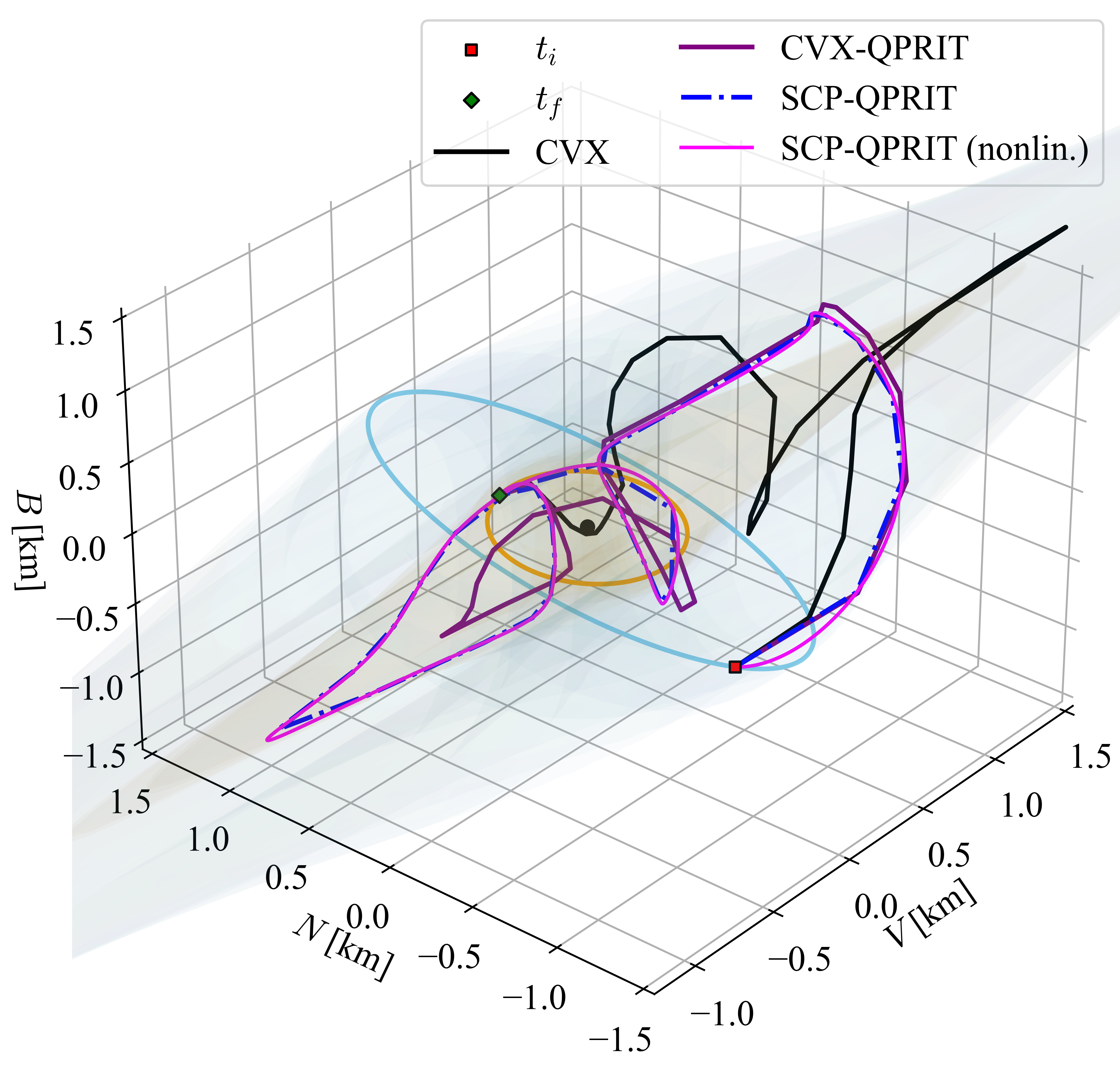}
         \caption{VNB} frame
         \label{fig:OptTraj_er3bp_tnw}
     \end{subfigure} 
     \begin{subfigure}[b]{0.35\textwidth}
         \centering
         \includegraphics[width=\textwidth]{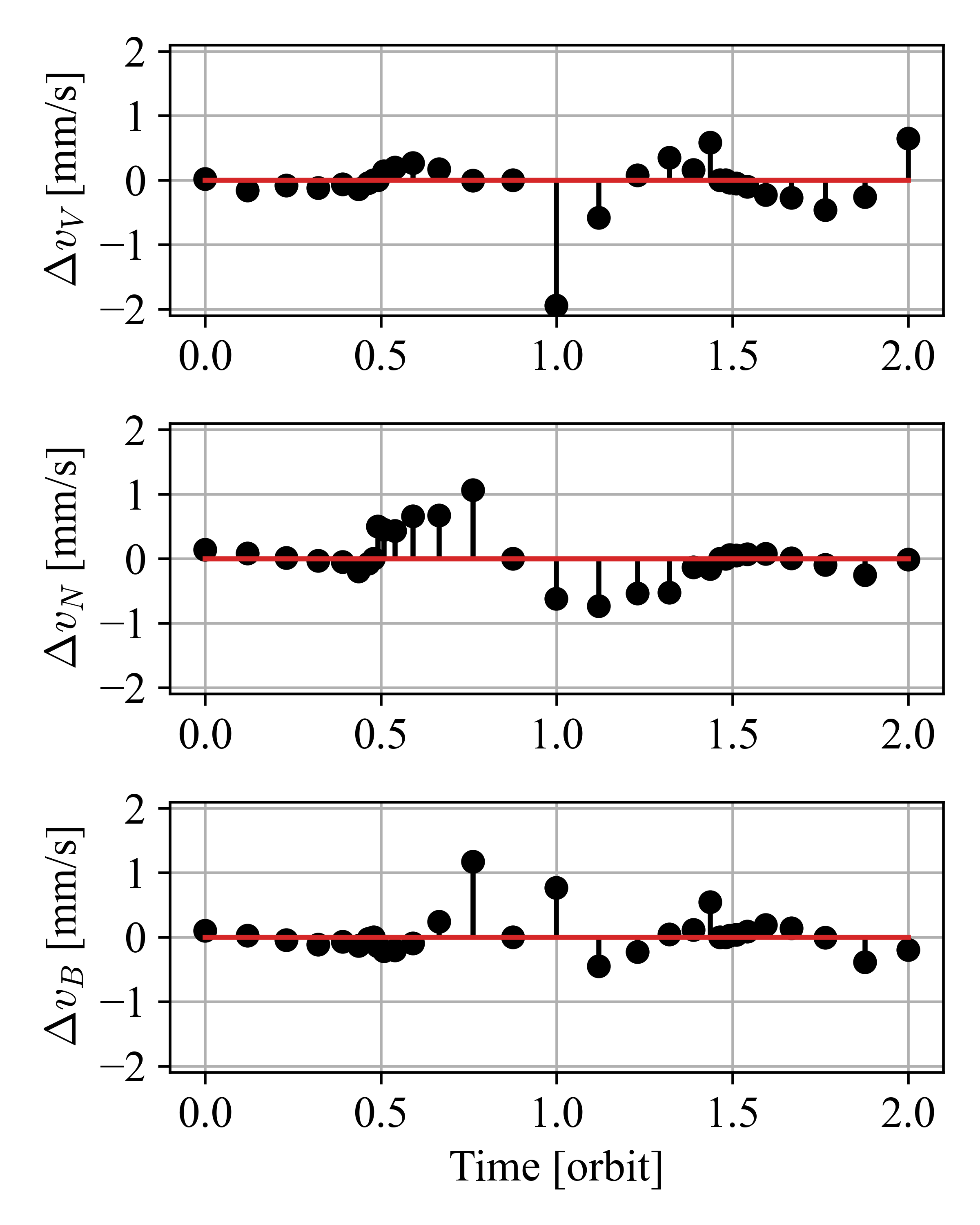}
         \caption{Control History (VNB} frame)
         \label{fig:cntrl_hist_QPRIT_er3bp}
     \end{subfigure}  \\
     \centering
     \begin{subfigure}[b]{0.35\textwidth}
         \centering
         \includegraphics[width=\textwidth]{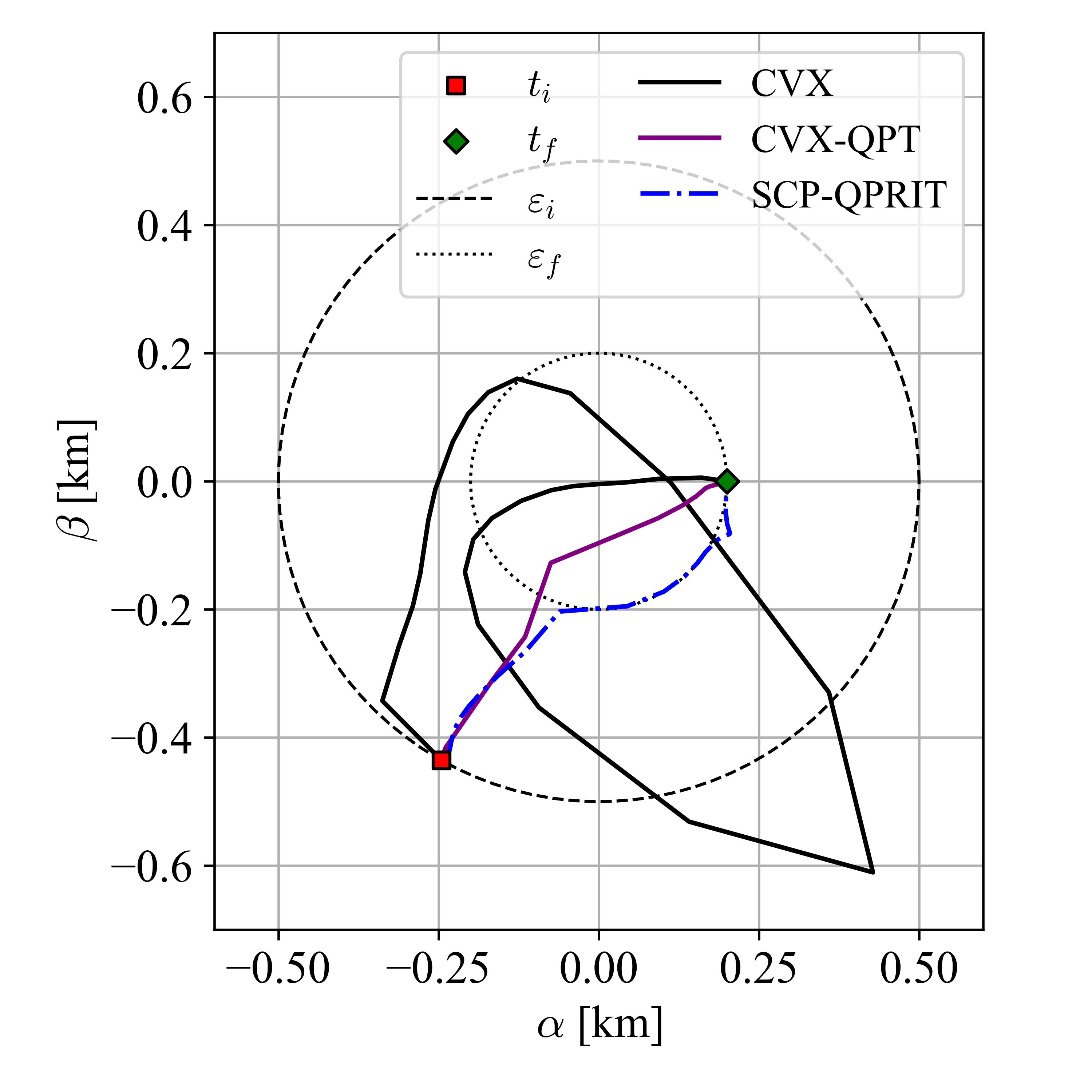}
         \caption{$\alpha-\beta$ plane}
         \label{fig:traj_ltc1_er3bp}
     \end{subfigure}
     \begin{subfigure}[b]{0.35\textwidth}
         \centering
         \includegraphics[width=\textwidth]{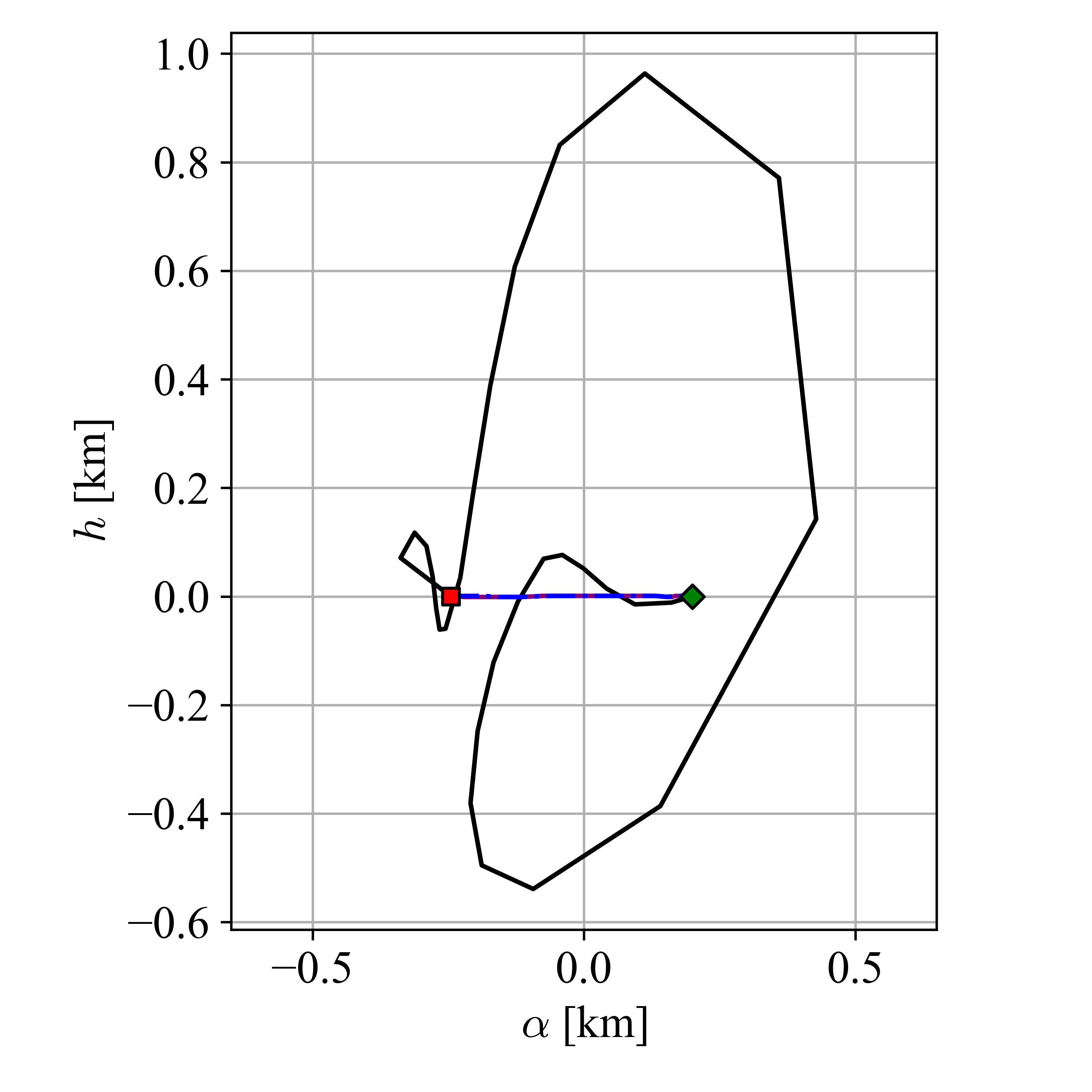}
         \caption{$\alpha-h$ plane}
         \label{fig:traj_ltc2_er3bp}
     \end{subfigure}  
    \caption{Optimized trajectory in the ER3BP.}
    \label{fig:opt_traj_qpt_er3bp} 
\end{figure}

In ER3BP, periodic orbits only permit sidereal resonance due to the eccentricity of the Earth-Moon orbit.
In this case study, the $L_2$ South 4:1 NRHO is chosen as the reference orbit, as shown in Fig. \ref{fig:4_1_er3bp}.
The optimized passively-safe solution (SCPP-QPRIT) is presented in Fig. \ref{fig:opt_traj_qpt_er3bp}, where Fig. \ref{fig:OptTraj_er3bp_tnw} is the trajectory shown in the VNB frame, Fig. \ref{fig:cntrl_hist_QPRIT_er3bp} is the control history. 
Furthermore, Figs. \ref{fig:traj_ltc1_er3bp} and \ref{fig:traj_ltc2_er3bp} are the trajectories plotted in the LTC. 
Similarly to the previous subsection, Table \ref{tab:exp_results2} summarizes the optimization results in CVX, CVX-QPRIT, and SCP-QPRIT, respectively.

\begin{table}[ht]
\centering
\caption{Experiment results of open-loop trajectories in ER3BP and BCR4BP (Sec. \ref{sec:mid_fidelity_demo})}
\label{tab:exp_results2}
\scalebox{0.8}{
    \begin{tabular}{c ccc c ccc}
    \hline\hline
    & \multicolumn{3}{c}{ER3BP} & & \multicolumn{3}{c}{BCR4BP} \\
    \cline{2-4} \cline{6-8}
    \multicolumn{1}{c}{Safety Strategy} 
    & CVX & CVX-QPRIT & SCP-QPRIT & & CVX & CVX-QPRIT & SCP-QPRIT \\
    \hline
    STM Computation [s]  & \multicolumn{3}{c}{87.23} &&  \multicolumn{3}{c}{28.46}  \\
    Opt. runtime [s] & 0.131 & 0.170 & 2.590 && 0.051 & 0.200 & 1.378 \\
    One-rev. passive safety &  & & \checkmark && & & \checkmark \\
    Fuel cost [mm/s]  & 8.192 & 13.330 & 13.770 && 2.051 & 5.559 & 5.637 \\
    \makecell{Nonlinear error \\ (terminal position) [m]} & 0.038  & 0.666 & 0.0164 && 0.870 & 0.667 & 0.657\\
    \hline\hline
    \end{tabular}
    }
\end{table}

As shown in Fig.~\ref{fig:OptTraj_er3bp_tnw}, the shape and orientation of the osculating invariant circles at the initial and terminal epoch are different because of the 4:1 resonance of the reference, while the transfer time is 2 revolutions. 
Still, the SCP-QPRIT trajectory successfully performs safe reconfiguration using the quasi-periodic structure, compared to the CVX solution that has a total cost of $\Delta v_{tot} = 8.19$ mm/s. 
Particularly, the solution of CVX-QPRIT results in the trajectory that stays on QPRIT surfaces for the entire transfer, although the deputy intrudes inside the terminal QPRIT, as shown in its LTC history (cf. Fig.~\ref{fig:traj_ltc1_er3bp}) and also smaller minimum separation as shown in the VNB frame (cf. Fig.~\ref{fig:OptTraj_er3bp_tnw}).
The total costs of CVX-QPRIT and SCP-QPRIT are $\Delta v_{tot} = 13.33$ mm/s and 13.77 mm/s, respectively, where passive safety is achieved via a marginal increase in transfer cost. 
When propagated under the nonlinear ER3BP dynamics, the terminal positional error between the linearized trajectory and the open-loop nonlinear trajectory is 0.0164 m, validating the accuracy of the linearized relative dynamics model in the ER3BP. 

\subsubsection{BCR4BP $L_2$ South 3:1 Halo Orbit (synodic resonance)}

For the demonstration in the BCR4BP, the $L_2$ South 3:1 synodic resonance NRHO is chosen as a reference orbit. 
As shown in Fig. \ref{fig:3_1_bcr4bp}, this orbit has a significant perturbation from the frequency-matching CR3BP solution \cite{boudad2020dynamics}. 
Note that only synodic resonance orbits are permitted for the reference orbits in the BCR4BP due to the Sun's revolution around the Earth-Moon barycenter in the synodic frame. 

\begin{figure}[ht!]
     \centering
     \begin{subfigure}[b]{0.5\textwidth}
         \centering
         \includegraphics[width=\textwidth]{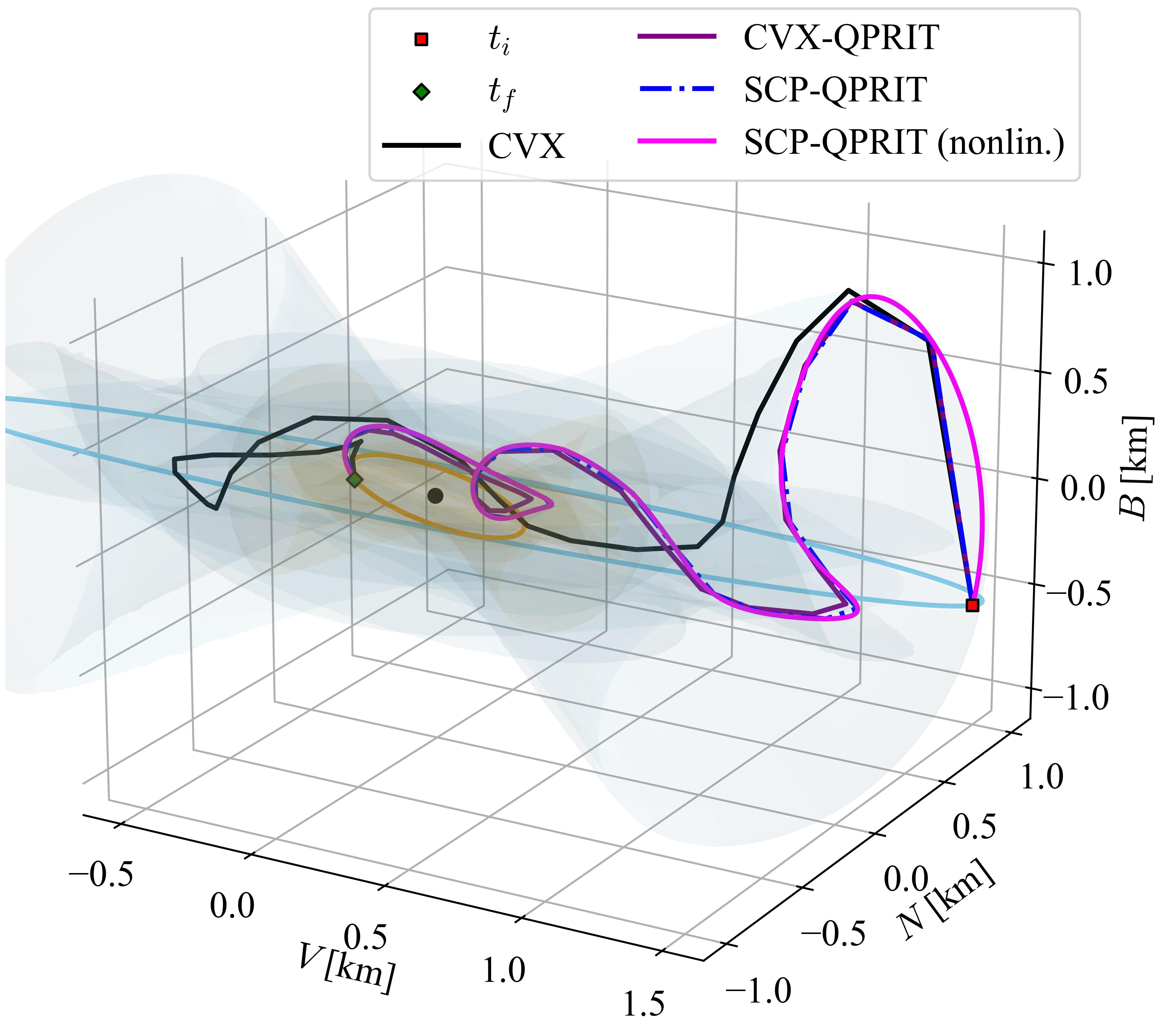}
         \caption{VNB} frame
         \label{fig:OptTraj_bcr4bp_tnw}
     \end{subfigure}
     \begin{subfigure}[b]{0.35\textwidth}
         \centering
         \includegraphics[width=\textwidth]{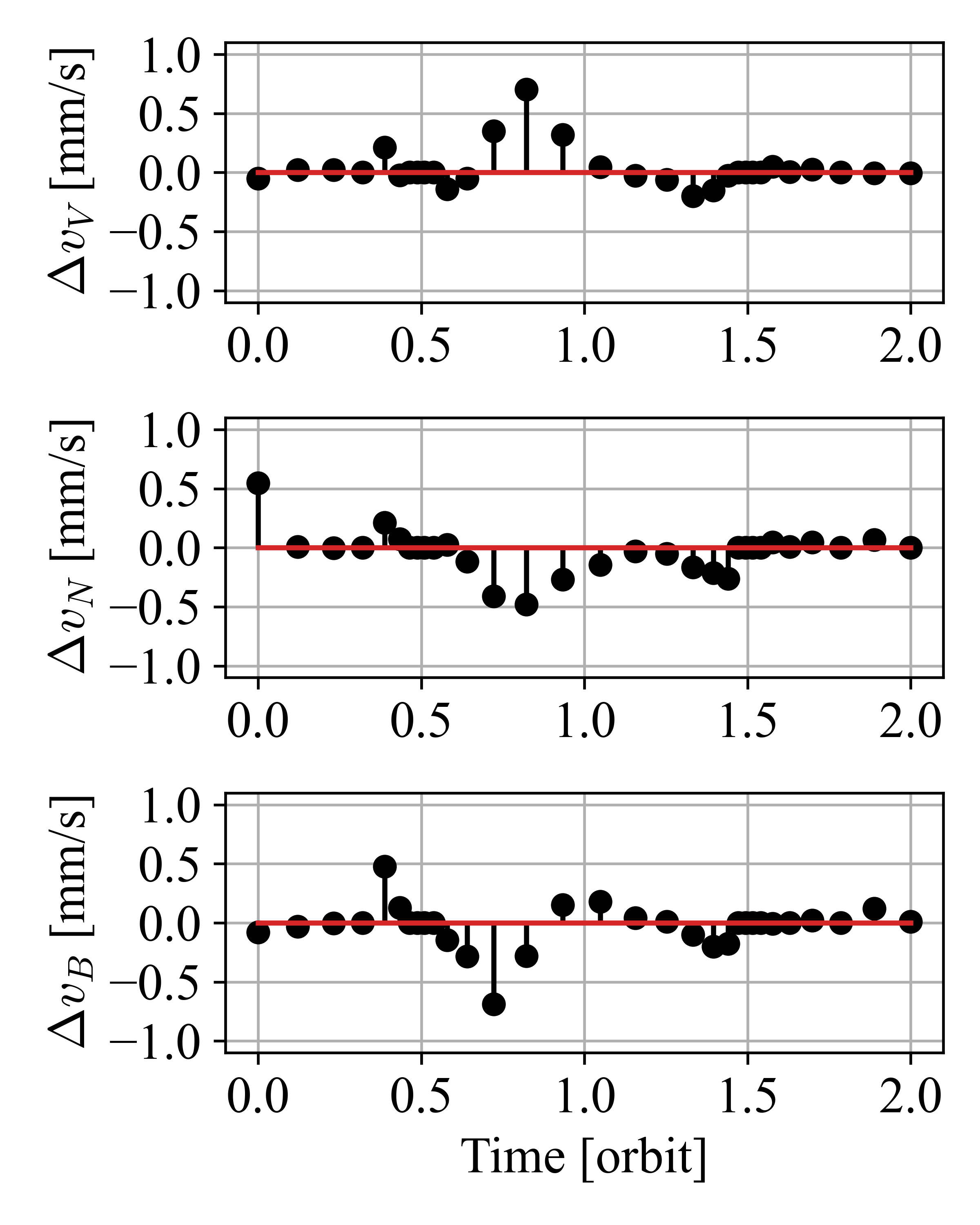}
         \caption{Control History (VNB} frame)
         \label{fig:cntrl_hist_QPRIT_bcr4bp}
     \end{subfigure}  \\
     \centering
     \begin{subfigure}[b]{0.4\textwidth}
         \centering
         \includegraphics[width=\textwidth]{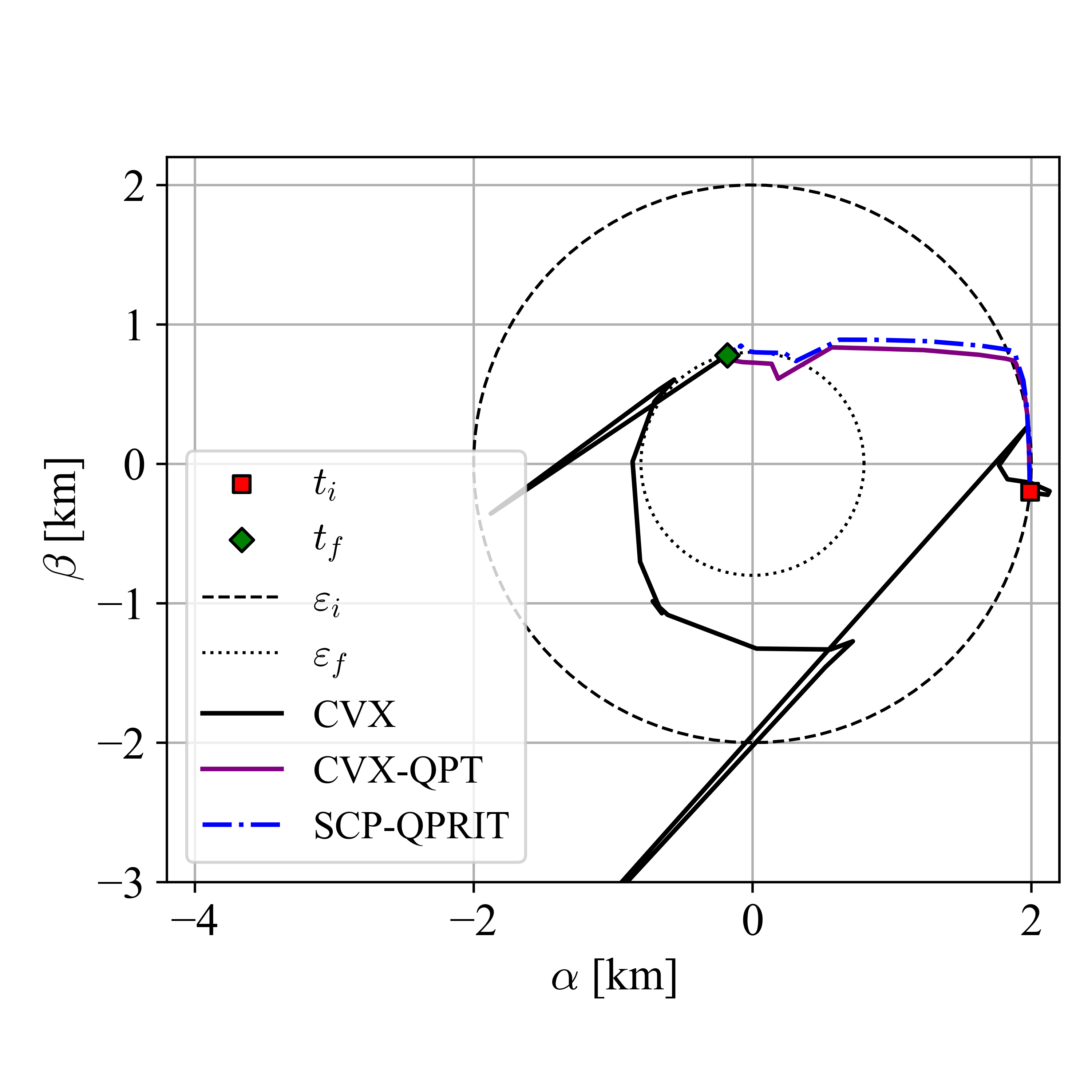}
         \caption{$\alpha-\beta$ plane}
         \label{fig:traj_ltc1_bcr4bp}
     \end{subfigure}
     \begin{subfigure}[b]{0.35\textwidth}
         \centering
         \includegraphics[width=\textwidth]{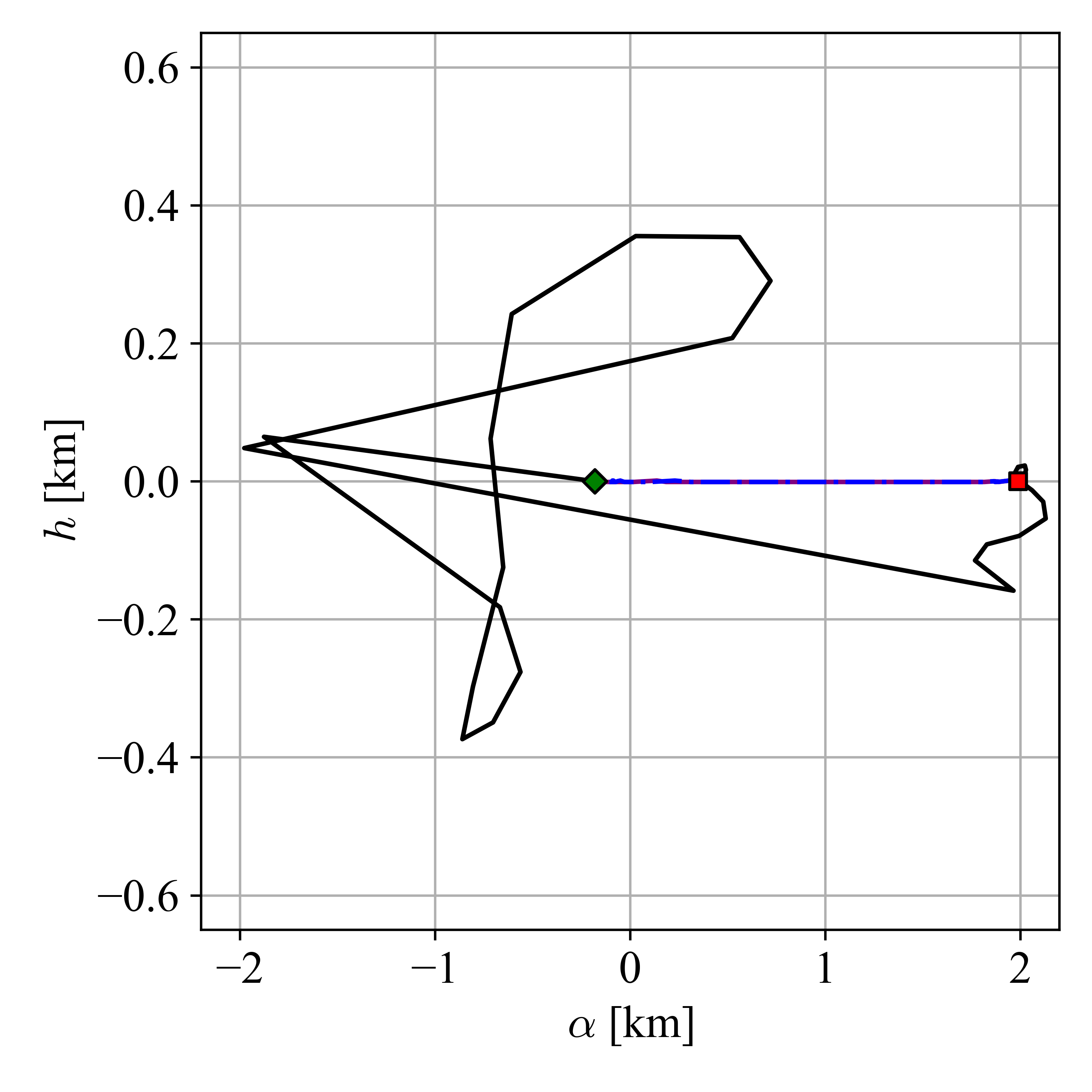}
         \caption{$\alpha-h$ plane}
         \label{fig:traj_ltc2_bcr4bp}
     \end{subfigure}  
    \caption{Optimal passively-safe solution in the BCR4BP.}
    \label{fig:opt_traj_qpt_bcr4bp} 
\end{figure}

The optimization results are presented in Fig. \ref{fig:opt_traj_qpt_bcr4bp} in the same fashion to Fig. \ref{fig:opt_traj_qpt_er3bp}.
First, as shown in Fig. \ref{fig:OptTraj_bcr4bp_tnw}, the shape of the QPRIT is remarkably different from the previous two cases. 
The perturbation of the reference orbit leads to the highly eccentric invariant curve. 
The optimized transfer cost of the CVX, CVX-QPRIT, and SCP-QPRIT is $\Delta v_{tot}$ = 2.051 mm/s,  5.559 mm/s, and  5.637 mm/s, respectively.  
It is observed that the SCP-QPRIT solution stays along the surface of the torus (cf. Figs. \ref{fig:OptTraj_bcr4bp_tnw} and \ref{fig:traj_ltc2_bcr4bp}) before reducing $\varepsilon$. 
The CVX-QPRIT solution almost satisfies passive safety, where only slight modification is needed in the following SCP to ensure the constraint in Eq. \ref{eq:scp_qpt_constr2}.
Another key observation is the high sensitivity of the plan OCP solution in the LTC space. This is due to the rapid rotation of the invariance circle, where the slight change in the Cartesian states leads to a spike in the LTC when the change is orthogonal to the change of LTC over time. 
This sensitivity is suppressed by simply imposing the constraint in Eq. \ref{eq:scp_qpt_constr1}.
The terminal nonlinear error between the SCP-QPRIT solution under the linearized dynamics and the nonlinear dynamics of BCR4BP is 0.657m, which is sufficiently low for the open-loop guidance and control of this formation reconfiguration problem. 

\subsection{N-spacecraft Swarm Reconfiguration in a Full-Ephemeris Dynamics Model} \label{sec:hifi_demo}

\begin{figure}[ht!]
    \centering    
    \includegraphics[width=0.8\textwidth]{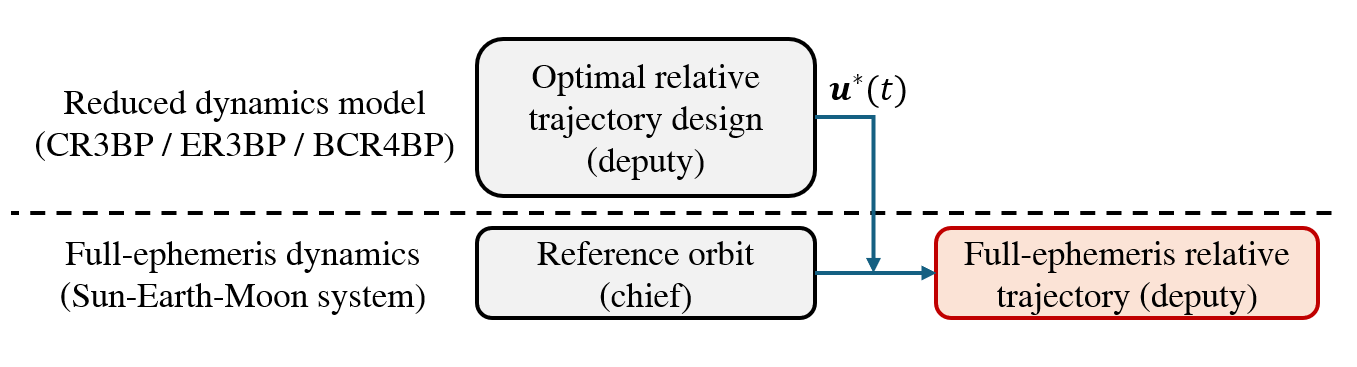}
    \caption{Validation process in the full-ephemeris model}
    \label{fig:validation_full_ephem}
\end{figure}

As the final case study, the proposed method is validated in a five-spacecraft swarm reconfiguration in the $L_2$ South 9:2 synodic NRHO under the full-ephemeris dynamics model, where exact periodic solutions are not attainable. 
The full-ephemeris model is based on the JPL SPICE toolkit \footnote{\url{https://naif.jpl.nasa.gov/naif/toolkit.html}} and considers the point-mass model of the Sun-Earth-Moon system, where the initial epoch is set to 2025-01-01 00:00:00 UTC.
Other perturbations, such as spherical harmonics or solar radiation pressure, are not considered. 
The validation process is illustrated in Fig. \ref{fig:validation_full_ephem}.
Because the LTC (and an oscillatory mode) is well-defined only if the reference orbit is perfectly periodic, an optimal relative guidance trajectory is generated within the simplified dynamics model, which permits periodic solutions. 
In this case study, the trajectory is optimized in the relative dynamics of CR3BP. 
In contrast, the chief orbit is separately obtained by converging a trajectory in the full-ephemeris model. 
In this case study, the 9-revolutions of NRHO is continued using multiple-shooting with the Newton-Raphson method from the CR3BP solution \cite{zimovan2023baseline}. 
The full-ephemeris NRHO is shown in Fig. \ref{fig:9_2_nrho}.
The optimal relative control profile obtained in the first step is then propagated in the full-ephemeris model along with the chief orbit that is continued in the full-ephemeris model. 
Note that the same dynamics model that is used to continue the full-ephemeris NRHO is also used for the propagation of relative motion. 
By performing this validation process, the station-keeping of the reference orbit around the Moon is effectively neglected from both the chief and deputy trajectories. 
This enables the extraction of the error in relative dynamics models between the simplified dynamics model (CR3BP) and the full-ephemeris dynamics. 
Since the propagation period in this case study is up to a few revolutions, this can serve as a simple validation scheme of the relative orbital control based on the proposed method. 

Because of the simplifying approximation of the CR3BP, the dynamics modeling error causes a large deviation when the feed-forward control is propagated in the full-ephemeris dynamics model. 
A shrinking-horizon MPC framework is introduced to address this issue, where the nonconvex trajectory optimization is re-solved via the SCP for the remaining horizon at predefined epochs $t_{\text{MPC}}$. 
The corresponding time step indices $k_{\text{MPC}}$ are listed in Table \ref{tab:sim_setup}, where the OCP is resolved twice during the transfer at $k_{\text{MPC}}=$ 13 and 25.
Throughout the transfer, a history of STM, control input matrix, and eigenvector $\boldsymbol{w}(t)$ is carried over (i.e., the LTC is not re-initialized).
In this case study, navigation uncertainty is not explicitly considered, and each agent is assumed to have perfect knowledge, while more sophisticated event-based MPC may be considered using the state covariance information \cite{Guffanti2023Visors}. 
However, to account for the actuation error, unmodeled acceleration (e.g., differential in spherical harmonics and solar radiation pressure), a linear process noise is added to the dynamics propagation, modeled as $Q(t_{k+1} - t_k)$, where $Q$ represents the noise matrix (cf. Table \ref{tab:sim_setup}). 

\begin{figure}[ht!]
     \centering
     \begin{subfigure}[b]{0.38\textwidth}
         \centering
         \includegraphics[width=\textwidth]{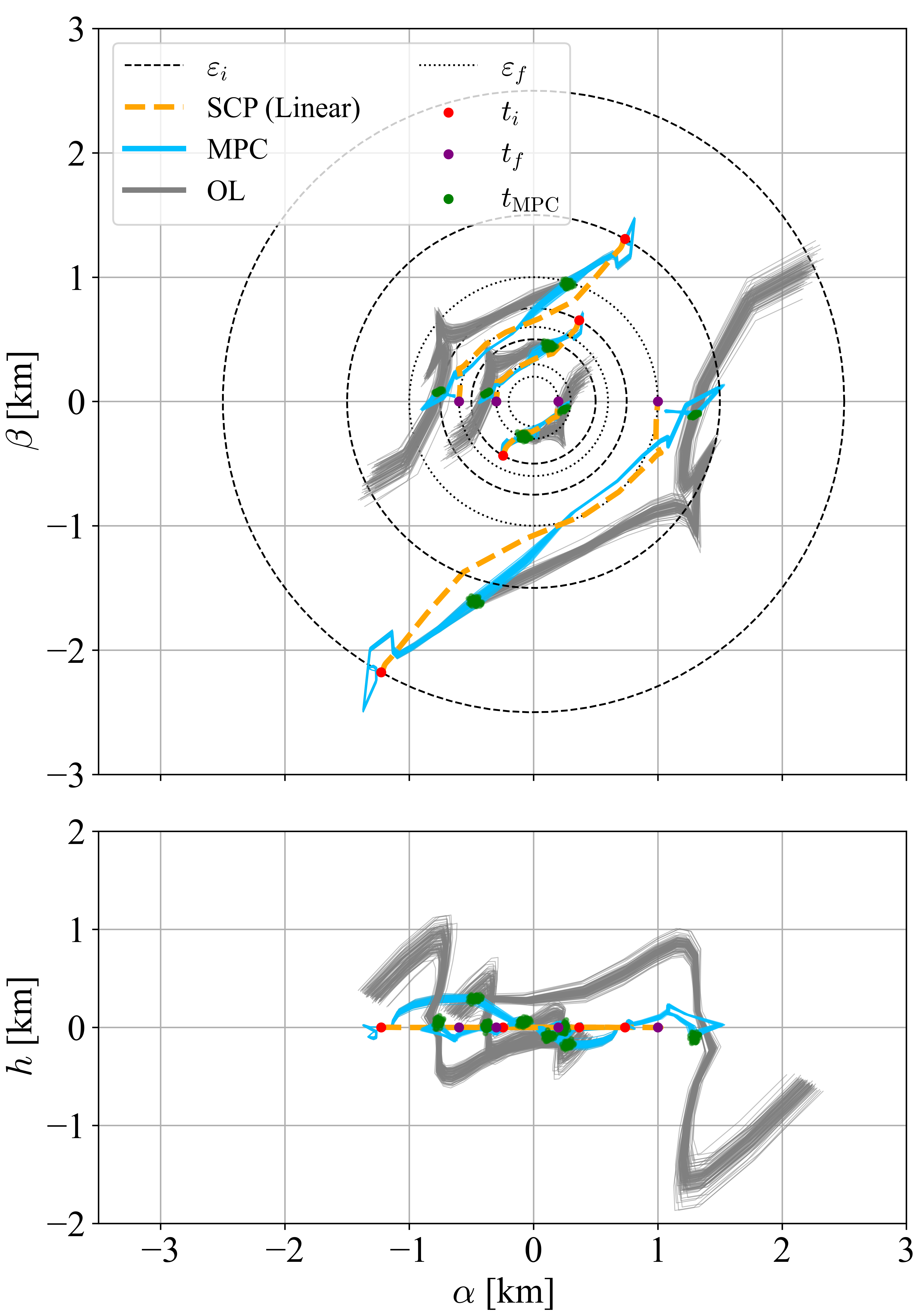}
         \caption{LTC}
         \label{fig:mpc_ltc}
     \end{subfigure}  
      \begin{subfigure}[b]{0.6\textwidth}
     \centering
     \includegraphics[width=\textwidth]{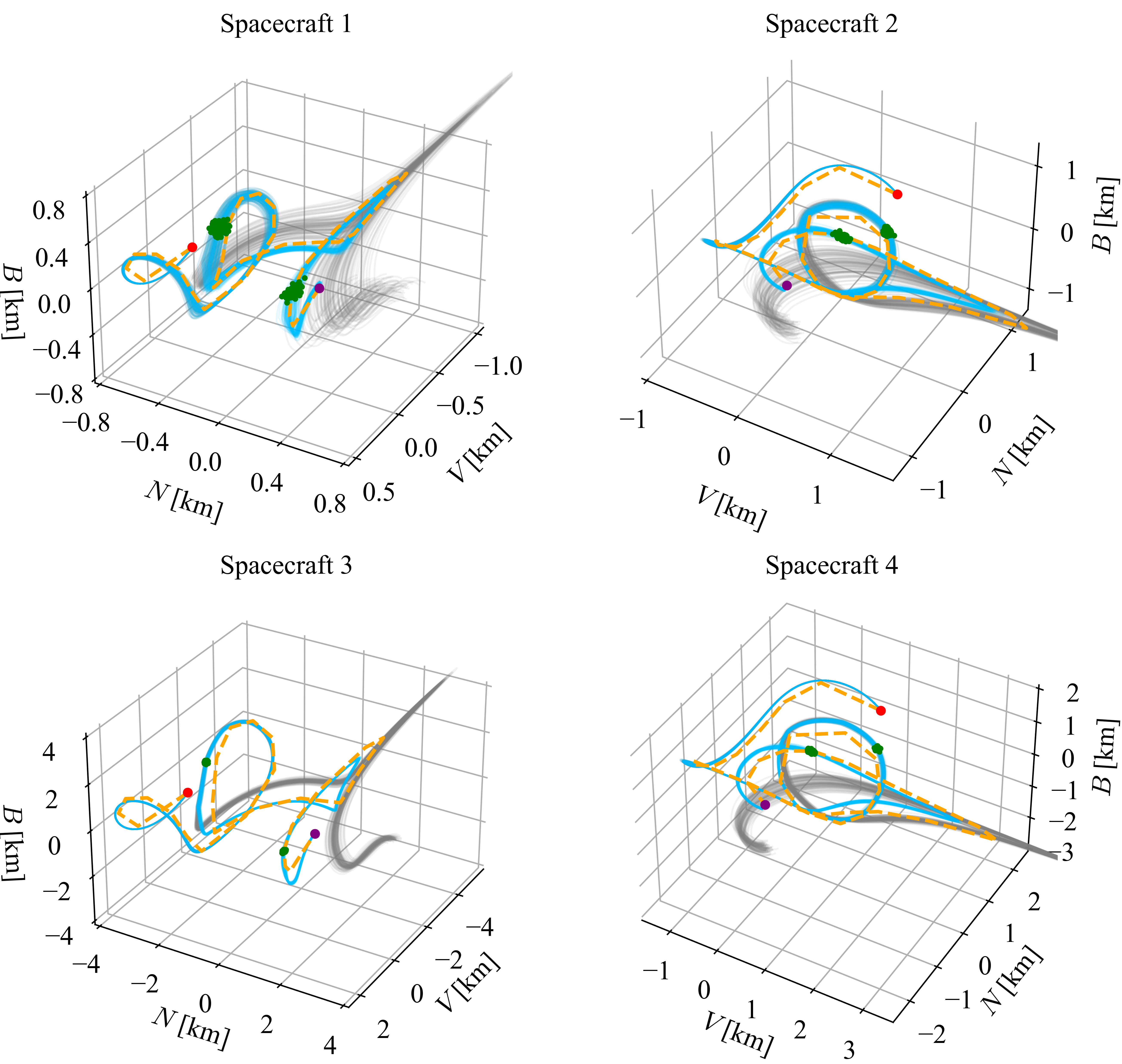}
     \caption{Trajectories}
     \label{fig:mpc_traj}
     \end{subfigure}
    \caption{Trajectories from 100 Monte-Carlo simulations of swarm reconfiguration in the full‑ephemeris dynamics model under Open‑Loop (OL) control and MPC.}
    \label{fig:mpc_sim} 
\end{figure}

A hundred Monte-Carlo simulations are conducted for both open-loop control and closed-loop control via MPC, where the ensemble of trajectories is shown in Fig. \ref{fig:mpc_sim}. 
Note that Fig. \ref{fig:mpc_ltc} represents the simulated trajectories in the LTC of the reference NRHO in the CR3BP. 
Although this coordinate frame does not exist in the full-ephemeris model, it represents the deviation from the expected behaviors of the trajectory with respect to the ideal local eigenspace. 
Furthermore, Table \ref{tab:mpc_result} presents the mean and standard deviation of the terminal position, velocity, and LTC of each deputy (labeled as SC1 to SC4), respectively. 
The leftmost column also displays the total fuel usage $\mathcal{J}_j$ for each individual spacecraft.

\begin{table}[ht!]
    \centering
    \scalebox{0.75}{
    \begin{tabular}{cc c c cc c cccccc}
    \hline \hline
    & & \multicolumn{1}{c}{Cost} & &  \multicolumn{9}{c}{Terminal State errors} \\
    \cline{3-3} \cline{5-13}
    & & $\mathcal{J}_j$ [mm/s] &  & $\rho_{\text{VNB}}$ [m] & $\dot{\rho}_{\text{VNB}}$ [mm/s] & &  $h$ [m] & $\alpha$ [m] &  $\beta$ [m]  & $\dot{h}$ [mm/s]  & $\dot{\alpha}$ [mm/s] & $\dot{\beta}$ [mm/s] \\
    \hline
    \multirow{2}{*}{SC1} & OL  & 8.33      &     & 367.9 $\pm$ 98.3     & 2.540 $\pm$ 0.675    && 66.3 $\pm$ 36.0 & 199.0 $\pm$ 49.9 & 188.0 $\pm$ 75.6 & 1.39 $\pm$ 0.78 & 0.89 $\pm$ 0.37 & 0.34 $\pm$ 0.26 \\
     & MPC & $9.15 \pm 0.84$ && 56.2  $\pm$ 7.0     & 0.652 $\pm$ 0.131    && 50.4 $\pm$ 7.6  & 11.9  $\pm$ 3.9  & 12.7  $\pm$ 4.9  & 0.07 $\pm$ 0.05 & 0.27 $\pm$ 0.08 & 1.06 $\pm$ 0.12 \\
    \hline
    \multirow{2}{*}{SC2} & OL  & 12.99        &   & 604.3 $\pm$ 104.0    & 4.120 $\pm$ 0.760    && 126.0 $\pm$ 34.2 & 327.0 $\pm$ 49.2 & 307.0 $\pm$ 90.0 & 2.14 $\pm$ 0.98 & 1.45 $\pm$ 0.44 & 0.43 $\pm$ 0.29 \\
     & MPC & $13.80 \pm 0.83$ && 73.7  $\pm$ 6.9 & 0.879 $\pm$ 0.111    && 63.8  $\pm$ 7.4  & 17.8  $\pm$ 4.0  & 20.5  $\pm$ 5.0  & 0.09 $\pm$ 0.07 & 0.41 $\pm$ 0.08 & 1.35 $\pm$ 0.14 \\
    \hline
    \multirow{2}{*}{SC3} & OL  & 47.66   &       & 2196.1 $\pm$ 94.4    & 14.800 $\pm$ 0.676   && 572.0 $\pm$ 35.2 & 1200.0 $\pm$ 47.4 & 1060.0 $\pm$ 78.0 & 6.60 $\pm$ 0.93 & 5.14 $\pm$ 0.39 & 2.71 $\pm$ 0.40 \\
     & MPC & $52.52 \pm 1.17$ && 199.0  $\pm$ 7.6  & 2.680  $\pm$ 0.155   && 150.0 $\pm$ 8.3  & 64.1  $\pm$ 4.3  & 76.6  $\pm$ 5.1  & 0.62 $\pm$ 0.09 & 1.60 $\pm$ 0.10 & 3.15 $\pm$ 0.18 \\
    \hline
    \multirow{2}{*}{SC4} & OL  & 27.71   &       & 1280.0 $\pm$ 96.8    & 8.670 $\pm$ 0.759    && 315.0 $\pm$ 39.3 & 698.0 $\pm$ 46.6 & 628.0 $\pm$ 81.3 & 4.03 $\pm$ 0.93 & 3.03 $\pm$ 0.39 & 1.34 $\pm$ 0.47 \\
     & MPC & $30.22 \pm 1.11$ & & 123.0  $\pm$ 7.9  & 1.550 $\pm$ 0.150    && 97.0  $\pm$ 8.5  & 37.2  $\pm$ 4.0  & 44.1  $\pm$ 5.2  & 0.35 $\pm$ 0.09 & 0.89 $\pm$ 0.10 & 2.08 $\pm$ 0.17 \\
     \hline \hline 
    \end{tabular}
    }
    \caption{Statistical Performance of the 100 Monte-Carlo Simulations for Open-Loop (OL) trajectories and MPC trajectories.}
    \label{tab:mpc_result}
\end{table}

It is evident from both the figure and the table that the open-loop control leads to a severe deviation of trajectories from the nominal trajectory designed in the CR3BP. 
The trajectories quickly diverge from the surface of the QPRIT (i.e., $\alpha-\beta$ plane) and struggle to reach the terminal state, particularly for the reconfiguration from large QRPITs (SC3 and SC4). 
In contrast, the MPC scheme successfully replans the reconfiguration trajectories and contracts its dispersion in both position and velocity, compared to the open-loop control, resulting in a significant improvement of the terminal state error with about a 10 \% increase in individual transfer cost from the nominal solution. 
Dispersions and errors in $h$ and $\dot{h}$ are particularly mitigated, indicating that swarms stay on the surface of the SPRIT and therefore maintain passive safety. 

\begin{figure}
    \centering
    \includegraphics[width=0.7\linewidth]{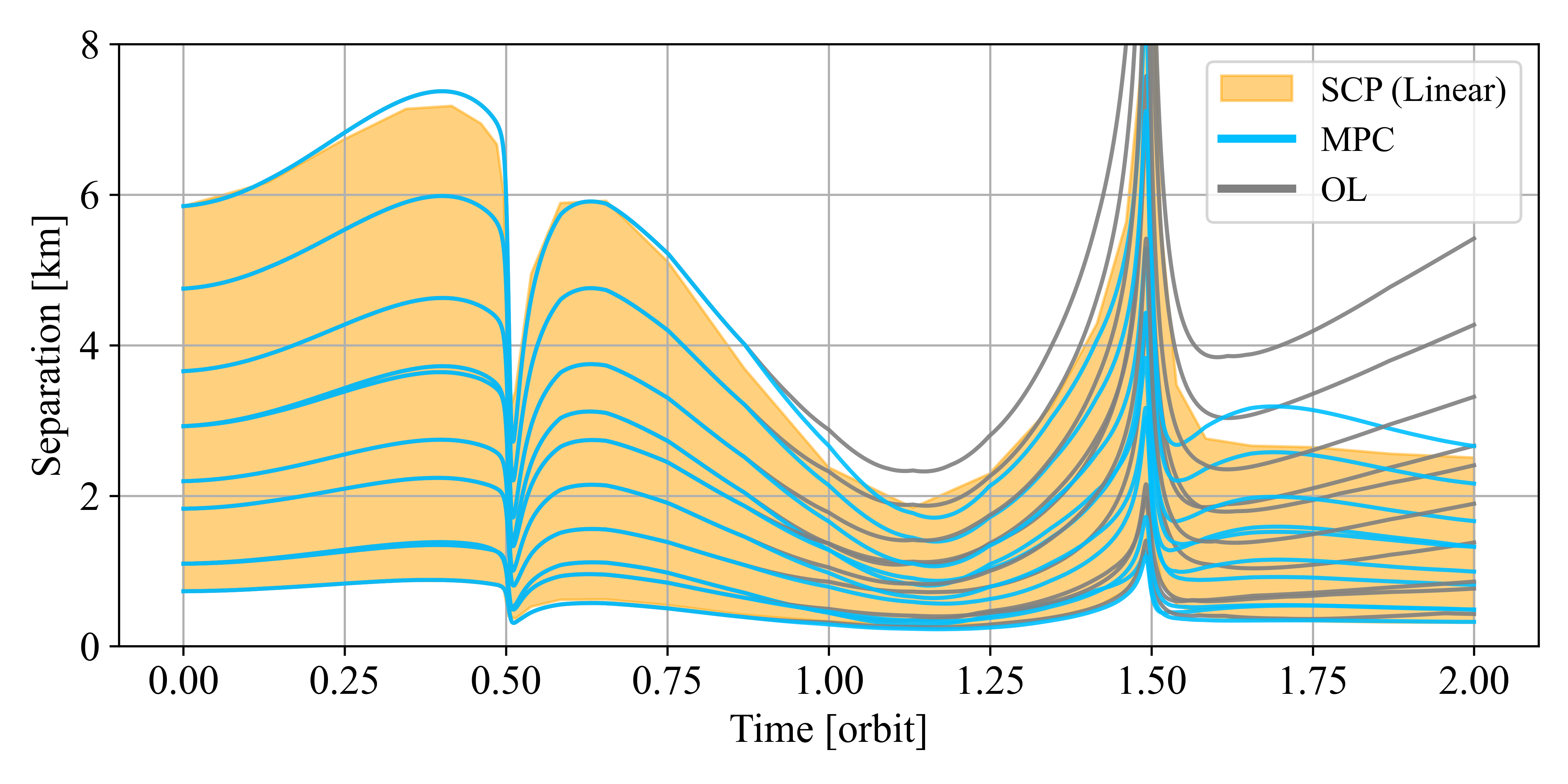}
    \caption{Comparison of inter‑spacecraft separation: the nominal bounds in the SCP solution and the worst‑case separations between each agent from 100 Monte‑Carlo simulations in the full-ephemeris model (Open-Loop (OL) and MPC trajectories).}
    \label{fig:minsep_mpc}
\end{figure}

Finally, Fig. \ref{fig:minsep_mpc} presents the worst separation between each agent in the 100 rollouts. 
The plot is overlaid with the upper and lower bounds of the inter-spacecraft separation in the CR3BP (i.e., open-loop nominal solution of the SCP-QPRIT). 
It is confirmed that, particularly when tracked by the MPC scheme, the inter-spacecraft separation in the full-epheemris model can be well approximated by the nominal trajectories optimized in the low-fidelity dynamics model, even under the process noise. 

As a final note, it is observed from the extensive testing that MPC trigger epochs are highly sensitive to the chosen full‑ephemeris reference orbit.
Even minor perturbations to the full‑ephemeris NRHO can produce different trajectories in the CR3BP LTC space, complicating swarm‑level implementation.
Therefore, in future research, a concurrent control scheme of station-keeping of reference orbit and around the Moon and formation reconfiguration between spacecraft should be addressed.

\section{Conclusion}

Safety of bounded relative motion in Restricted Multi-Body Problems (RMBPs) is one of the largest concerns in the practical deployment of spacecraft swarms in non-Keplerian cislunar orbits.
This paper presents a novel passively safe relative reconfiguration strategy for RMBPs.
The main contribution is the formulation of the optimal control problem in the Local Toroidal Coordinates (LTC), a time-varying frame that projects the relative motion onto the oscillatory eigenspace that enables the bounded quasi-periodic relative motion.
In the LTC, safety constraints drawn from the quasi-periodic structure are significantly simplified and decentralized between each agent, enabling the large-scale swarm reconfiguration. 
To support this method, relative motion of the PRMBPs is resolved in the Velocity/Normal/Binormal(VNB) frame, a co-moving frame that aligns with the spacecraft velocity. 
The VNB frame is advantageous to discuss the relative motion in highly elliptic orbits, such as Near Rectilinear Halo Orbits (NRHOs), as it slows the natural relative motion in the radial direction.  
Furthermore, geometric relationships of the bounded quasi-periodic relative motion in the RMBPs and the bounded periodic relative motion in the Perturbed Restricted Two-Body Problem (PR2BP) are analyzed through their local eigensystems. 
The proposed reconfiguration strategy is demonstrated in various cislunar NRHOs in the Circular Restricted Three-Body Problem, Elliptic Three-Body Problem, and Bicircular Restricted Four-Body Problem. 
Furthermore, the deployment of the guidance strategy in the full-ephemeris model is successfully performed through model predictive control. 
This paper builds a practical stepping stone for the deployment of large-scale spacecraft swarming operations in a highly perturbed dynamics in cislunar space. 
Future research includes the integration of station-keeping of a reference orbit and relative motion control, as well as a streamlined swarm design optimization method based on obaservability analysis. 

\section*{Appendix}

\subsection{Pseudo-Potential Functions in ER3BP and BCR4BP}

\subsubsection{ER3BP}
The pseudo-potential of the ER3BP is given by \cite{hiday1992optimal}
\begin{align}
    \Upsilon_{\text{ER3BP}} = \frac{1}{2}\omega_z^2(x^2 + y^2) + \frac{\mu}{r_{ec}} + \frac{1-\mu}{r_{mc}}.
\end{align}
The distance between the Earth and the Moon, the angular velocity, and its time derivative are given by 
\begin{align}
    \omega_z = \frac{\sqrt{1-e^2}}{(1-e \cos{E})^2}, \quad [\dot{\omega}_z]_{\mathcal{B}} = \frac{-2e\sqrt{1-e^2} \sin{E}}{(1-e \cos{E})^4}, \quad r_{em} = 1 - e \cos E,
\end{align}
where $E$ is the eccentric anomaly of the Moon. 
The evolution of the eccentric anomaly is given based on the progression of the mean anomaly as 
\begin{align}
    \tilde{M} = E - e \sin E, \quad \tilde{M}(t_1) = \tilde{M}(t_0) + (t_1 - t_0).
\end{align}

\subsubsection{BCR4BP}

The pseudo-potential of the BCR4BP is given by \cite{boudad2020dynamics, takubo2023optimization}
\begin{align}
    \Upsilon_{\text{BCR4BP}} = \frac{1}{2}\omega_z^2(x^2 + y^2) + \frac{\mu}{r_{ec}} + \frac{1-\mu}{r_{mc}} + \frac{\mu_s}{r_{sc}} - \frac{\mu_s}{(r_{b_1s}/r_{em})^3}(x_sx + y_sy + z_sz),
\end{align}
with constant $\omega_z=1$ (i.e., $\dot{\omega}=0$). 

\subsection{Derivation of the Angular Velocoity and Acceleration of the VNB Frame to the Moon Synodic Frame}

The derivation of the angular velocity of the VNB frame with respect to the Moon synodic frame  $\boldsymbol{\omega}_{v/m}$ follows the procedure elaborated in Ref. \cite{franzini_relative_2019}.
The derivation first requires a time derivative of basis vectors of the VNB frame expressed in the Moon synodic frame: $\left\{[\dot{\hat{\boldsymbol{\imath}}}]_{\mathcal{M}}, [\dot{\hat{\boldsymbol{\jmath}}}]_{\mathcal{M}}, [\dot{\hat{\boldsymbol{k}}}]_{\mathcal{M}}\right\} $. 

The time derivative of $\boldsymbol{\hat{\imath}}$ is computed as follows: 

\begin{subequations}
\label{eq:i_deriv}
\begin{align}
    [\dot{\hat{\boldsymbol{\imath}}}]_{\mathcal{M}} 
    & = \dfrac{d}{dt} \left( \dfrac{[\dot{\boldsymbol{r}}]_\mathcal{M}}{v} \right) 
    = \dfrac{1}{v} ( [\ddot{\boldsymbol{r}}]_\mathcal{M} - \hat{\boldsymbol{\imath}} \cdot (\hat{\boldsymbol{\imath}} \cdot [\dot{\boldsymbol{r}}]_\mathcal{M} )) 
    \quad 
    \left( \because \dot{v} = \dfrac{[\dot{\boldsymbol{r}}]_\mathcal{M} \cdot [\ddot{\boldsymbol{r}}]_\mathcal{M}}{v} = \hat{\boldsymbol{\imath}} \cdot  [\ddot{\boldsymbol{r}}]_\mathcal{M} \right) \\
    & = \dfrac{1}{v} ( ([\ddot{\boldsymbol{r}}]_\mathcal{M} \cdot \hat{\boldsymbol{\jmath}}) \hat{\boldsymbol{\jmath}} + ([\ddot{\boldsymbol{r}}]_\mathcal{M} \cdot \hat{\boldsymbol{k}}) \hat{\boldsymbol{k}} ) \\
    & = \dfrac{1}{hv} ( [\ddot{\boldsymbol{r}}]_\mathcal{M} \cdot \boldsymbol{h} ) \hat{\boldsymbol{\jmath}} + \dfrac{1}{v} ([\ddot{\boldsymbol{r}}]_\mathcal{M} \cdot \hat{\boldsymbol{k}}) \hat{\boldsymbol{k}} .
\end{align}
\end{subequations}
Using the definition in Eq. \eqref{eq:tnw_basis}, the time derivative of $\boldsymbol{\hat{\jmath}}$ is computed as follows:
\small
\begin{subequations}
\label{eq:j_deriv}
\begin{align}
    [\dot{\hat{\boldsymbol{\jmath}}}]_{\mathcal{M}} 
    & = \dfrac{1}{h^2} ([\dot{\boldsymbol{h}}]_\mathcal{M} h - \boldsymbol{h} \dot{h}) \\
    & = \dfrac{1}{h} ([\dot{\boldsymbol{h}}]_\mathcal{M} - \dot{h}{\hat{\boldsymbol{\jmath}}}) \label{eq:j_deriv2} \\
    & = \dfrac{1}{h} ( [\dot{\boldsymbol{h}}]_\mathcal{M}  \cdot {\hat{\boldsymbol{\imath}}}) {\hat{\boldsymbol{\imath}}} + \dfrac{1}{h} ( [\dot{\boldsymbol{h}}]_\mathcal{M}  \cdot {\hat{\boldsymbol{k}}}) {\hat{\boldsymbol{k}}} \label{eq:j_deriv3} \\
    & = - \dfrac{1}{hv} ( [\ddot{\boldsymbol{r}}]_\mathcal{M} \cdot \boldsymbol{h} ) \hat{\boldsymbol{\imath}} + \dfrac{1}{h} ( [\dot{\boldsymbol{h}}]_\mathcal{M}  \cdot {\hat{\boldsymbol{k}}}) {\hat{\boldsymbol{k}}} 
\end{align}
\end{subequations}
\normalsize
From Eq. \ref{eq:j_deriv2} to \ref{eq:j_deriv3}, the following equality is used, noting that $  [\dot{\boldsymbol{h}}]_\mathcal{V} = \dot{h} \hat{\boldsymbol{\jmath}} $:
\small
\begin{subequations}
\begin{align}
    \dot{h} 
    & = [\dot{\boldsymbol{h}}]_\mathcal{V} \cdot \hat{\boldsymbol{\jmath}} 
    = [\dot{\boldsymbol{h}}]_\mathcal{M} \cdot \hat{\boldsymbol{\jmath}} - (\boldsymbol{\omega}_{v/m} \times \boldsymbol{h}) \cdot \hat{\boldsymbol{\jmath}} 
    = [\dot{\boldsymbol{h}}]_\mathcal{M} \cdot \hat{\boldsymbol{\jmath}}. 
    \quad \left(\because (\boldsymbol{\omega}_{v/m} \times \boldsymbol{h}) \cdot \hat{\boldsymbol{\jmath}} = \dfrac{1}{h}(\boldsymbol{\omega}_{v/m} \times \boldsymbol{h}) \cdot \boldsymbol{h} = 0 \right)
\end{align}
\end{subequations}
\normalsize
Also note that $[\dot{\boldsymbol{h}}]_\mathcal{M} = \boldsymbol{r} \times [\ddot{\boldsymbol{r}}]_\mathcal{M}$.
Finally, the time derivative of $\hat{\boldsymbol{k}}$ is computed as follows:
\begin{align}\label{eq:k_deriv}
    [\dot{\hat{\boldsymbol{k}}}]_{\mathcal{M}} 
    & = [\dot{\hat{\boldsymbol{\imath}}}]_{\mathcal{M}} \times \hat{\boldsymbol{\jmath}} + \hat{\boldsymbol{\imath}} \times [\dot{\hat{\boldsymbol{\jmath}}}]_{\mathcal{M}} = - \dfrac{1}{v} ([\ddot{\boldsymbol{r}}]_\mathcal{M} \cdot \hat{\boldsymbol{k}}) \hat{\boldsymbol{\imath}} - \dfrac{1}{h} ( [\dot{\boldsymbol{h}}]_\mathcal{M}  \cdot {\hat{\boldsymbol{k}}}) {\hat{\boldsymbol{\jmath}}}.
\end{align}
The time derivative of unit vectors of the VNB frame expressed in the Moon synodic frame can also be expressed using $\boldsymbol{\omega}_{v/m}$ as follows:
\begin{align}
    [\dot{\hat{\boldsymbol{\imath}}}]_\mathcal{M} = \boldsymbol{\omega}_{v/m} \times \hat{\boldsymbol{\imath}}, \quad 
    [\dot{\hat{\boldsymbol{\jmath}}}]_\mathcal{M} = \boldsymbol{\omega}_{v/m} \times \hat{\boldsymbol{\jmath}}, \quad 
    [\dot{\hat{\boldsymbol{k}}}]_\mathcal{M} = \boldsymbol{\omega}_{v/m} \times \hat{\boldsymbol{k}}
\end{align}
Taking a cross product with each unit vector from the left side, the following expression is obtained:
\small
\begin{align}
    \hat{\boldsymbol{\imath}} \times [\dot{\hat{\boldsymbol{\imath}}}]_\mathcal{M} = \boldsymbol{\omega}_{v/m} - ( \boldsymbol{\omega}_{v/m} \cdot \hat{\boldsymbol{\imath}} ) \hat{\boldsymbol{\imath}}, 
    \quad 
    \hat{\boldsymbol{\jmath}} \times [\dot{\hat{\boldsymbol{\jmath}}}]_\mathcal{M} = \boldsymbol{\omega}_{v/m} - ( \boldsymbol{\omega}_{v/m} \cdot \hat{\boldsymbol{\jmath}} ) \hat{\boldsymbol{\jmath}}, 
    \quad 
    \hat{\boldsymbol{k}} \times [\dot{\hat{\boldsymbol{k}}}]_\mathcal{M} = \boldsymbol{\omega}_{v/m} - ( \boldsymbol{\omega}_{v/m} \cdot \hat{\boldsymbol{k}} ) \hat{\boldsymbol{k}}.
\end{align}
\normalsize
By summing up the three equations, the following equation is derived: 
\begin{align}
    \boldsymbol{\omega}_{v/m} = \dfrac{1}{2} \left( \hat{\boldsymbol{\imath}} \times [\dot{\hat{\boldsymbol{\imath}}}]_\mathcal{M} + \hat{\boldsymbol{\jmath}} \times [\dot{\hat{\boldsymbol{\jmath}}}]_\mathcal{M} + \hat{\boldsymbol{k}} \times [\dot{\hat{\boldsymbol{k}}}]_\mathcal{M} \right)
\end{align}
With Eqs. \eqref{eq:i_deriv}, \eqref{eq:j_deriv}, and \eqref{eq:k_deriv}, the angular velocity $\boldsymbol{\omega}_{v/m}$ is expressed in the VNB frame as follows: 
\small
\begin{subequations}
\begin{align}
    \boldsymbol{\omega}_{v/m} 
    & = {\omega}^{x}_{v/m} \hat{\boldsymbol{\imath}} + {\omega}^{y}_{v/m} \hat{\boldsymbol{\jmath}} + {\omega}^{z}_{v/m} \hat{\boldsymbol{k}} \\
    & = \dfrac{1}{h} ( [\dot{\boldsymbol{h}}]_\mathcal{M} \cdot {\hat{\boldsymbol{k}}})  \hat{\boldsymbol{\imath}} 
    - \dfrac{1}{v} ([\ddot{\boldsymbol{r}}]_\mathcal{M} \cdot \hat{\boldsymbol{k}})  \hat{\boldsymbol{\jmath}} 
    + \dfrac{1}{hv} ( [\ddot{\boldsymbol{r}}]_\mathcal{M} \cdot \boldsymbol{h} ) \hat{\boldsymbol{k}} \\
    & = \dfrac{1}{h} ( ( \boldsymbol{r} \times [\ddot{\boldsymbol{r}}]_{\mathcal{M}}) \cdot {\hat{\boldsymbol{k}}})  \hat{\boldsymbol{\imath}} 
    - \dfrac{1}{v} ([\ddot{\boldsymbol{r}}]_\mathcal{M} \cdot \hat{\boldsymbol{k}})  \hat{\boldsymbol{\jmath}} 
    + \dfrac{1}{hv} ( [\ddot{\boldsymbol{r}}]_\mathcal{M} \cdot \boldsymbol{h} ) \hat{\boldsymbol{k}} 
\end{align}
\end{subequations}
\normalsize
Each component of $[\dot{\boldsymbol{\omega}}_{v/m}]_\mathcal{V}$ is computed by simply taking the component-wise time derivative of the above expression, which results in  
\small
\begin{subequations}
\begin{align}
    \dot{\omega}^{x}_{v/m} & = 
    - \dfrac{\dot{h}}{h^2} ((\boldsymbol{r} \times [\ddot{\boldsymbol{r}}]_{\mathcal{M}}) \cdot \hat{\boldsymbol{k}} ) 
    + \dfrac{1}{h} \left(  
    ([\dot{\boldsymbol{r}}]_{\mathcal{M}} \times [\ddot{\boldsymbol{r}}]_{\mathcal{M}} ) \cdot \hat{\boldsymbol{k}} 
    + (\boldsymbol{r} \times [\dddot{\boldsymbol{r}}]_{\mathcal{M}} ) \cdot \hat{\boldsymbol{k}} 
    + (\boldsymbol{r} \times [\ddot{\boldsymbol{r}}]_{\mathcal{M}}) \cdot [\dot{\hat{\boldsymbol{k}}}]_{\mathcal{M}}  \right) \\
    \dot{\omega}^{y}_{v/m} & = \dfrac{\dot{v}}{v^2} ([\ddot{\boldsymbol{r}}]_{\mathcal{M}} \cdot \hat{\boldsymbol{k}}) 
    - 
    \frac{1}{v} ([\dddot{\boldsymbol{r}}]_{\mathcal{M}} \cdot \hat{\boldsymbol{k}} +[\ddot{\boldsymbol{r}}]_{\mathcal{M}} \cdot [\dot{\hat{\boldsymbol{k}}}]_{\mathcal{M}}  )
    \\ 
    \dot{\omega}^{z}_{v/m} & = \dfrac{-\dot{h}v - h\dot{v}}{h^2 v^2} ([\ddot{\boldsymbol{r}}]_{\mathcal{M}} \cdot \boldsymbol{h}) 
    +
    \dfrac{1}{hv} ( [\ddot{\boldsymbol{r}}]_{\mathcal{M}} \cdot \boldsymbol{h}).
\end{align}
\end{subequations}
\normalsize

The acceleration and the jerk of the chief in CR3BP, ER3BP, and BCR4BP within the Moon synodic frame are provided in Appendix  \ref{app:acc_jerk}. 

\subsection{Acceleration and Jerks in the Moon Synodic Frame} \label{app:acc_jerk}

\subsubsection{CR3BP}

For the CR3BP, the acceleration and jerk of the chief in the Moon synodic frame are calculated as follows:
\small
\begin{subequations}
\begin{align} 
    [\ddot{\boldsymbol{r}}]_{\mathcal{M}} 
    =& -2 \boldsymbol{\omega}_{m / i} \times[\dot{\boldsymbol{r}}]_{\mathcal{M}} - \boldsymbol{\omega}_{m / i} \times\left(\boldsymbol{\omega}_{m / i} \times \boldsymbol{r}\right) - \mu \frac{\boldsymbol{r}}{r^3} - (1-\mu)\left(\frac{ \boldsymbol{r}_{em} + \boldsymbol{r} }{\left\|\boldsymbol{r}_{em} + \boldsymbol{r}\right\|^3} - \frac{\boldsymbol{r}_{em}}{r_{em}^3}\right) \label{eq:rddot_M_cr3bp} \\
    [\dddot{\boldsymbol{r}}]_{\mathcal{M}} 
    =& - 2 \boldsymbol{\omega}_{m / i} \times[\ddot{\boldsymbol{r}}]_{\mathcal{M}} 
    - \boldsymbol{\omega}_{m / i} \times \left(\boldsymbol{\omega}_{m / i} \times[\dot{\boldsymbol{r}}]_{\mathcal{M}}\right) 
    - \mu \frac{\partial}{\partial \boldsymbol{r}}\left[\frac{\boldsymbol{r}}{r^3}\right] {[\dot{\boldsymbol{r}}]}_{\mathcal{M}} 
    - (1-\mu) \frac{\partial}{\partial \boldsymbol{r}}\left[\frac{\boldsymbol{r}_{em} + \boldsymbol{r}}{\left\|\boldsymbol{r}_{em} + \boldsymbol{r} \right\|^3}\right] [\dot{\boldsymbol{r}}]_{\mathcal{M}},
\end{align}
\end{subequations}
\normalsize
where $\boldsymbol{r}_{em}$ is the position vector from Earth to the Moon, and for a vector $\boldsymbol{\zeta}$, 
\small
\begin{align}
    \dfrac{\partial}{\partial \boldsymbol{\zeta}} \left[ \dfrac{\boldsymbol{\zeta}}{\zeta^3}\right] = \dfrac{1}{\zeta^3} \left( \boldsymbol{I}_3 - 3 \dfrac{\boldsymbol{\zeta}\boldsymbol{\zeta}^\top}{\zeta^2} \right) .
\end{align}
\normalsize

\subsubsection{ER3BP}

Similarly, the acceleration and jerk of the chief in the ER3BP within the Moon synodic frame are calculated as follows: 
\small
\begin{subequations}
\begin{align} 
    [\ddot{\boldsymbol{r}}]_{\mathcal{M}} 
    =& -2 \boldsymbol{\omega}_{m / i} \times[\dot{\boldsymbol{r}}]_{\mathcal{M}} 
    - [\dot{\boldsymbol{\omega}}_{m/i}]_{\mathcal{M}} \times \boldsymbol{r}
    - \boldsymbol{\omega}_{m / i} \times\left(\boldsymbol{\omega}_{m / i} \times \boldsymbol{r}\right)
    - \mu \frac{\boldsymbol{r}}{r^3} - (1-\mu)\left(\frac{ \boldsymbol{r}_{em} 
    + \boldsymbol{r} }{\left\|\boldsymbol{r}_{em} + \boldsymbol{r}\right\|^3} - \frac{\boldsymbol{r}_{em}}{r_{em}^3}\right) \label{eq:rddot_M_er3bp} \\
    [\dddot{\boldsymbol{r}}]_{\mathcal{M}} 
    =& - 2 \boldsymbol{\omega}_{m / i} \times[\ddot{\boldsymbol{r}}]_{\mathcal{M}} 
    - 3 [\dot{\boldsymbol{\omega}}_{m / i}]_{\mathcal{M}} \times[\ddot{\boldsymbol{r}}]_{\mathcal{M}} 
    - [\ddot{\boldsymbol{\omega}}_{m / i}]_{\mathcal{M}} \times \boldsymbol{r} \\
    & 
    - [\dot{\boldsymbol{\omega}}_{m / i}]_{\mathcal{M}} \times \left(\boldsymbol{\omega}_{m / i} \times \boldsymbol{r} \right) 
    - \boldsymbol{\omega}_{m / i} \times \left(\dot{\boldsymbol{\omega}}_{m / i} \times \boldsymbol{r}_{\mathcal{M}}\right) 
    - \boldsymbol{\omega}_{m / i} \times \left(\boldsymbol{\omega}_{m / i} \times[\dot{\boldsymbol{r}}]_{\mathcal{M}}\right) \\
    & 
    - \mu \frac{\partial}{\partial \boldsymbol{r}}\left[\frac{\boldsymbol{r}}{r^3}\right] {[\dot{\boldsymbol{r}}]}_{\mathcal{M}} 
    - (1-\mu) \left( \frac{\partial}{\partial \boldsymbol{r}}\left[\frac{\boldsymbol{r}_{em} + \boldsymbol{r}}{\left\|\boldsymbol{r}_{em} + \boldsymbol{r} \right\|^3}\right] ([\dot{\boldsymbol{r}}]_{\mathcal{M}} + [\dot{\boldsymbol{r}}_{em}]_{\mathcal{M}}) - \frac{\partial}{\partial \boldsymbol{r}_{em}} \left[\frac{\boldsymbol{r}_{em}}{r_{em}^3}\right] [\dot{\boldsymbol{r}}_{em}]_{\mathcal{M}} \right),
\end{align}
\end{subequations}
\normalsize

\subsubsection{BCR4BP}

Finally, the acceleration and jerk of the chief in the BCR4BP within the Moon synodic frame are calculated as follows \cite{Romero2023cislunar}:
\small
\begin{subequations}
\begin{align} 
    [\ddot{\boldsymbol{r}}]_{\mathcal{M}} 
    =& -2 \boldsymbol{\omega}_{m / i} \times[\dot{\boldsymbol{r}}]_{\mathcal{M}} - \boldsymbol{\omega}_{m / i} \times\left(\boldsymbol{\omega}_{m / i} \times \boldsymbol{r}\right) 
    - \mu \frac{\boldsymbol{r}}{r^3} 
    - (1-\mu)\left(\frac{ \boldsymbol{r}_{em} + \boldsymbol{r} }{\left\|\boldsymbol{r}_{em} + \boldsymbol{r}\right\|^3} - \frac{\boldsymbol{r}_{em}}{r_{em}^3}\right) 
    - \mu_S \left(\frac{ \boldsymbol{r}_{sm} + \boldsymbol{r} }{\left\|\boldsymbol{r}_{sm} + \boldsymbol{r}\right\|^3} - \frac{\boldsymbol{r}_{sb_1}}{r_{sb_1}^3}\right)
    \label{eq:rddot_M_bcr4bp} \\
    [\dddot{\boldsymbol{r}}]_{\mathcal{M}} 
    =& - 2 \boldsymbol{\omega}_{m / i} \times[\ddot{\boldsymbol{r}}]_{\mathcal{M}} 
    - \boldsymbol{\omega}_{m / i} \times \left(\boldsymbol{\omega}_{m / i} \times[\dot{\boldsymbol{r}}]_{\mathcal{M}}\right) 
    - \mu \frac{\partial}{\partial \boldsymbol{r}}\left[\frac{\boldsymbol{r}}{r^3}\right] {[\dot{\boldsymbol{r}}]}_{\mathcal{M}} 
    - (1-\mu) \frac{\partial}{\partial \boldsymbol{r}}\left[\frac{\boldsymbol{r}_{em} + \boldsymbol{r}}{\left\|\boldsymbol{r}_{em} + \boldsymbol{r} \right\|^3}\right] [\dot{\boldsymbol{r}}]_{\mathcal{M}} \\
    & -\mu_S \left( \frac{\partial}{\partial \boldsymbol{r}}\left[\frac{\boldsymbol{r}_{sm} + \boldsymbol{r}}{\left\|\boldsymbol{r}_{sm} + \boldsymbol{r} \right\|^3}\right] ([\dot{\boldsymbol{r}}]_{\mathcal{M}} + [\dot{\boldsymbol{r}}_{sm}]_{\mathcal{M}}) - \frac{\partial}{\partial \boldsymbol{r}_{sb_1}} \left[\frac{\boldsymbol{r}_{sb_1}}{r_{sb_1}^3}\right] [\dot{\boldsymbol{r}}_{sb_1}]_{\mathcal{M}} \right),
\end{align}
\end{subequations}
\normalsize

\section*{Acknowledgments}
This work is supported by Blue Origin (SPO \#299266) as an Associate Member and Co-Founder of the Stanford’s Center of AEroSpace Autonomy Research (CAESAR), and National Science Foundation Graduate Research Fellowship Program (No. DGE-2146755). 
This article solely reflects the opinions and conclusions of its authors and not any sponsors.
Yuji Takubo acknowledges financial support from the Ezoe Memorial Recruit Foundation. 
Yuji Takubo also thanks Yuri Shimane for helpful discussions on continuing NRHO in full-ephemeris dynamics.
AI tools, including ChatGPT o1 and Claude 3.7 Sonnet, were used to assist in minimizing grammatical mistakes.

\bibliography{91_references}

\end{document}